\documentclass[11pt]{article}

\usepackage{contract,latexsym,amsmath,amsfonts,amsthm,amsbsy,amscd,graphics,
epsf}
\theoremstyle{plain}
\newtheorem{lemma}{Lemma}[section]
\newtheorem{corollary}[lemma]{Corollary}
\newtheorem{theorem}[lemma]{Theorem}
\newtheorem{conjecture}[lemma]{Conjecture}
\newtheorem{definition}{Definition}[section]

\theoremstyle{definition}
\newtheorem{example}{Example}[section]
\newtheorem{remark}{Remark}[section]

\numberwithin{equation}{section}

\newcounter{ll}
\renewcommand{\thell}{$(\roman{ll})$}
\newenvironment{enum}{\begin{list}{\thell}{\usecounter{ll} \topsep 2mm
	\itemsep 2mm \parsep 0mm \labelwidth 7mm \labelsep 2mm}}{\end{list}}

\newcommand{\tab}[2]{{\rm Tab}(#1,#2)}
\newcommand{\ttab}[2]{{\rm LRT}(#1,#2)}

\newcommand{\qbin}[2]{\genfrac{[}{]}{0pt}{}{#1}{#2}}
\newcommand{\bin}[2]{\genfrac{(}{)}{0pt}{}{#1}{#2}}

\def\Integer{\mathbb{Z}}
\def\Int{\Integer_{\geq 0}}

\def\bs{\boldsymbol}
\def\vL{\bs L}
\def\ve{{\bs e}}

\def\vh{\bs h}
\def\vn{\bs n}
\def\vm{\bs m}
\def\vQ{\bs Q}
\def\vk{\bs k}
\def\vrho{\bs \rho}
\def\l{\ell}
\def\lb{\overline{\l}}
\def\la{\lambda}
\def\lm{\lambda\mu}
\def\La{\Lambda}
\def\Om{\Omega}
\def\Omp{\Om_{\rm p}}
\def\Sb{\tilde{K}}
\def\A{{\cal A}}
\def\O{{\cal O}}
\def\C{{\cal C}}
\def\Cb{\overline{\C}}
\def\Cp{\C_{\rm p}}
\def\Cbp{\Cb_{\rm p}}
\def\Z{{\cal Z}}
\def\Zb{\overline{\Z}}
\def\P{{\cal P}}
\def\Pb{\overline{\P}}
\def\D{{\cal D}}
\def\U{{\cal U}}
\def\Dp{\D_{\rm p}}
\def\Up{\U_{\rm p}}
\def\T{{\cal T}}
\def\R{{\cal R}}
\def\B{{\cal B}}
\def\X{{\cal X}}
\def\W{{\cal W}}
\def\S{{\cal S}}
\def\G{{\cal G}}
\def\Pt{\tilde{P}}
\def\pt{\tilde{p}}
\def\wt{\tilde{w}}
\def\Id{{\rm Id}}
\def\content{{\rm content}}
\def\shape{{\rm shape}}
\def\height{{\rm height}}
\def\width{{\rm width}}
\def\rank{{\rm rank}}
\def\co{{\rm co}}
\def\col{{\rm col}}

\begin{document}

\title{Inhomogeneous lattice paths, generalized Kostka polynomials\\
and A$_{n-1}$ supernomials}

\author{
{\Large Anne Schilling}
\thanks{e-mail: {\tt schillin@wins.uva.nl}}
\quad and\quad
{\Large S. Ole Warnaar}
\thanks{e-mail: {\tt warnaar@wins.uva.nl}}\\
\mbox{} \\
{\em Instituut voor Theoretische Fysica, Universiteit van Amsterdam},\\
{\em Valckenierstraat 65, 1018 XE Amsterdam, The Netherlands}}
 
\date{\Large February, 1998}
\maketitle

\begin{abstract}
Inhomogeneous lattice paths are introduced as ordered sequences of rectangular
Young tableaux thereby generalizing recent work on the Kostka polynomials by 
Nakayashiki and Yamada and by Lascoux, Leclerc and Thibon. Motivated by these
works and by Kashiwara's theory of crystal bases we define a statistic on 
paths yielding two novel classes of polynomials. One of these provides a 
generalization of the Kostka polynomials while the other, which we name the
A$_{n-1}$ supernomial, is a $q$-deformation of the expansion coefficients of 
products of Schur polynomials. Many well-known results for Kostka polynomials 
are extended leading to representations of our polynomials in terms of a 
charge statistic on Littlewood--Richardson tableaux and in terms of fermionic 
configuration sums. Several identities for the generalized Kostka polynomials 
and the A$_{n-1}$ supernomials are proven or conjectured. Finally, a 
connection between the supernomials and Bailey's lemma is made.
\end{abstract}

\section{Introduction}

Lattice paths play an important r\^{o}le in combinatorics and exactly solvable 
lattice models of statistical mechanics. In particular, the one-dimensional
configuration sums necessary for the calculation of order parameters of 
lattice models are generating functions of lattice paths (see for 
example~\cite{ABF84,DJKMO89,JMO88}). 
A classic example of lattice paths is given by sequences of upward and 
downward steps.
The number of such paths consisting of $\la_1$ steps up and $\la_2$ steps down
is given by the binomial coefficient $\tbinom{\la_1+\la_2}{\la_1}
=(\la_1+\la_2)!/\la_1!\la_2!$ which is the expansion coefficient of 
\begin{equation*}
(x_1+x_2)^L=\sum_{\substack{\la_1,\la_2\\ \la_1+\la_2=L}} x_1^{\la_1} 
x_2^{\la_2}\bin{L}{\la_1}.
\end{equation*}
An important $q$-deformation of the binomial is the $q$-binomial 
\begin{equation}\label{qbin}
\qbin{\la_1+\la_2}{\la_1}=
\begin{cases}
\frac{(q^{\la_2+1})_{\la_1}}{(q)_{\la_1}} & \text{for $\la_1,\la_2\in\Int$,}\\
0 & \text{otherwise,}
\end{cases}
\end{equation}
where $(x)_n=(1-x)(1-xq)(1-xq^2)\cdots (1-xq^{n-1})$. 
The $q$-binomial can be interpreted as the generating function of all 
paths with $\la_1$ steps up and $\la_2$ steps down where each path is 
weighted as follows. 
Let $p_1,\ldots,p_{\la_1+\la_2}$ denote the steps of the path where we 
label a step up by 1 and a step down by 2.
Then the weight of the path is given by $\sum_{i=1}^{\la_1+\la_2-1} i 
\chi(p_i<p_{i+1})$ where $\chi({\rm true})=1$ and $\chi({\rm false})=0$.

Other $q$-functions have occurred, such as a $q$-deformation of the 
trinomial coefficients~\cite{AB87,A94} in the expansion of 
$(x_1^2+x_1x_2+x_2^2)^L$ or, more generally, of the $(N+1)$-nomial 
coefficients~\cite{K95,S96,W97} in the expansion of $h_N^L$ where 
$h_N$ is the complete symmetric polynomial in the variables $x_1$ and $x_2$ of 
degree $N$.
In a study of Rogers--Ramanujan-type identities the following generalizations 
of the multinomial coefficients were introduced~\cite{SW97} 
\begin{equation}\label{supern}
h_1^{L_1} \dots h_N^{L_N} = \sum_{\substack{\la_1,\la_2\\\la_1+\la_2=\ell_N}}
x_1^{\la_1} x_2^{\la_2}\binom{\vL}{\la_1-\frac{1}{2}\ell_N},
\end{equation}
where $\vL=\sum_{i=1}^N L_i\ve_i\in\Int^N$ with $\ve_i$ the 
$i$th unit vector in $\Integer^N$ and $\ell_i=\sum_{j=1}^N \min\{i,j\}L_j$.
Since equation~\eqref{supern} reduces to the definition of the $(i+1)$-nomial 
coefficient when $\vL=L\ve_i$ (up to a shift in the lower index), the 
expansion coefficient in~\eqref{supern} was coined (A$_1$) supernomial.
In ref.~\cite{SW97} it was shown that many Rogers--Ramanujan-type identities 
admit bounded analogues involving the following $q$-deformation of the 
supernomial
\begin{equation}\label{qsuper}
\qbin{\vL}{a}=\sum_{j_1+\cdots+j_N=a+\frac{\ell_N}{2}}
q^{\sum_{k=1}^{N-1}(\ell_{k+1}-\ell_k-j_{k+1})j_k}
\qbin{L_N}{j_N}\qbin{L_{N-1}+j_N}{j_{N-1}}\cdots\qbin{L_1+j_2}{j_1}
\end{equation}
for $\vL\in\Int^N$ and $a+\frac{1}{2}\ell_N=0,1,\ldots,\ell_N$.
However, the question whether~\eqref{qsuper} or the 
Rogers--Ramanujan-type identities involving~\eqref{qsuper} can be interpreted 
as generating functions of weighted lattice paths remained unanswered.
Incidentally, the polynomials in equation \eqref{qsuper} have occurred in 
Butler's study~\cite{B87}-\cite{B94} of finite abelian groups.

In a seemingly unrelated development, Nakayashiki and Yamada~\cite{NY95}
introduced the notion of ``inhomogeneous'' lattice paths by considering
paths in which each of the elementary steps $p_i$ can be chosen from a 
different set $\B_i$.
The main result of their work is a new combinatorial representation of the 
Kostka polynomial as the generating function of inhomogeneous paths 
where either all $\B_i$ are sets of fully symmetric (one-row) Young
tableaux or all $\B_i$ are sets of fully antisymmetric (one-column) Young 
tableaux. An equivalent description of the Kostka polynomials, formulated 
in terms of the plactic monoid, was found by Lascoux, Leclerc and 
Thibon~\cite{LLT95}.

The purpose of this paper is to elucidate the connection of the work
of Nakayashiki and Yamada and of Lascoux, Leclerc and Thibon on the 
Kostka polynomials with that of ref.~\cite{SW97} on supernomials and to
extend all of them.
In particular, we introduce inhomogeneous lattice paths based on Young
tableaux with mixed symmetries, or more precisely, on Young tableaux
of rectangular shape. Motivated by the theory of crystal bases~\cite{K91}
we assign weights to these paths and relate their generating functions to 
$q$-deformations of the (A$_{n-1}$) supernomials defined through products of 
Schur polynomials, in the spirit of equation~\eqref{supern}.
By imposing suitable restrictions on the inhomogeneous lattice paths,
we obtain new polynomials that include the Kostka polynomials as a special
case. For these generalized Kostka polynomials we derive several extensions
of classical results such as a Lascoux--Sch\"utzenberger-type
representation~\cite{LS78} in terms of a charge statistic, and a 
representation akin to that of Kirillov and Reshetikhin~\cite{KR88} based on 
rigged configurations.
We furthermore prove and conjecture several identities involving the 
A$_{n-1}$ supernomials and generalized Kostka polynomials and briefly 
comment on a Bailey-type lemma~\cite{B49} for ``antisymmetric'' supernomials.
The A$_{n-1}$ supernomials include polynomials previously studied 
in~\cite{B87}-\cite{B94},\cite{HKKOTY98}.

The rest of this paper is organized as follows. 
Section~\ref{sec_review} serves to set  the notation used throughout the
paper and to review some basic definitions and properties of Young tableaux, 
words and Kostka polynomials. In section~\ref{sec_SK} inhomogeneous lattice 
paths based on rectangular Young tableaux are introduced 
(def~\ref{def_paths}). A statistic on such 
paths originating from crystal-base theory is used to define A$_{n-1}$ 
supernomials (def~\ref{def_super}) and generalized Kostka polynomials 
(def~\ref{def_gKp}) as generating functions of
inhomogeneous paths. We furthermore map the paths underlying the generalized 
Kostka polynomials onto Littlewood--Richardson (LR) tableaux 
(def~\ref{def_Stab}).
In section~\ref{sec_in_cyc} an initial cyclage and charge statistic on
LR tableaux is defined (def~\ref{def_charge}) which enables us in the 
subsequent section to prove a Lascoux--Sch\"utzenberger-type representation 
for generalized Kostka polynomials (cor~\ref{cor_SK}).
Section~\ref{sec_poset} deals with more general $\la$-(co)cyclages on 
LR tableaux, showing that these cyclages impose a ranked poset structure 
on the set of LR tableaux (thm~\ref{theo_poset}). These results are used 
in section~\ref{sec_propSK} to prove a duality formula for the generalized
Kostka polynomials (thm~\ref{theo_Kinv}) and recurrence relations for 
the A$_{n-1}$ supernomials and the generalized Kostka polynomials 
(thm~\ref{theo_rec}). 
In section~\ref{sec_fermi} the recurrences are employed to obtain a
Kirillov--Reshetikhin-type expression for the generalized Kostka polynomials
(thm~\ref{theo_KF}).
We finally conclude in section~\ref{sec_dis} with some conjectured polynomial 
identities and with a Bailey-like lemma involving the A$_{n-1}$ supernomials.

\section{Young tableaux, words and Kostka polynomials}
\label{sec_review}

This section reviews some definitions and properties of Young tableaux, 
words and Kostka polynomials and sets out the notation and terminology
used throughout the paper. For more details the reader may consult
refs.~\cite{B94,F97,M95}.

Throughout, we denote by $|A|$ the cardinality of a set $A$ and we
define $|\mu|=\sum_i \mu_i$ for an array of numbers $\mu=(\mu_1,\mu_2,\ldots)$.

\subsection{Young tableaux and words}

We begin by recalling some definitions regarding partitions.
A partition $\la=(\la_1,\la_2,\ldots)$ is a weakly decreasing sequence of 
non-negative integers such that only finitely many $\la_i\neq 0$.
We write $\la\vdash n$ if $|\la|=n$.
Partitions which differ only by a string of zeros are identified.
Each partition can be depicted by a Young diagram, which (adopting the
``French'' convention) is a collection of 
boxes with left-adjusted rows of decreasing length from bottom to top. If
$\la=(\la_1,\la_2,\ldots)$ is a partition then the corresponding Young diagram
has $\la_i$ boxes in the $i$th row from the bottom. For example
\begin{equation*}
\epsfxsize=2 cm \epsffile{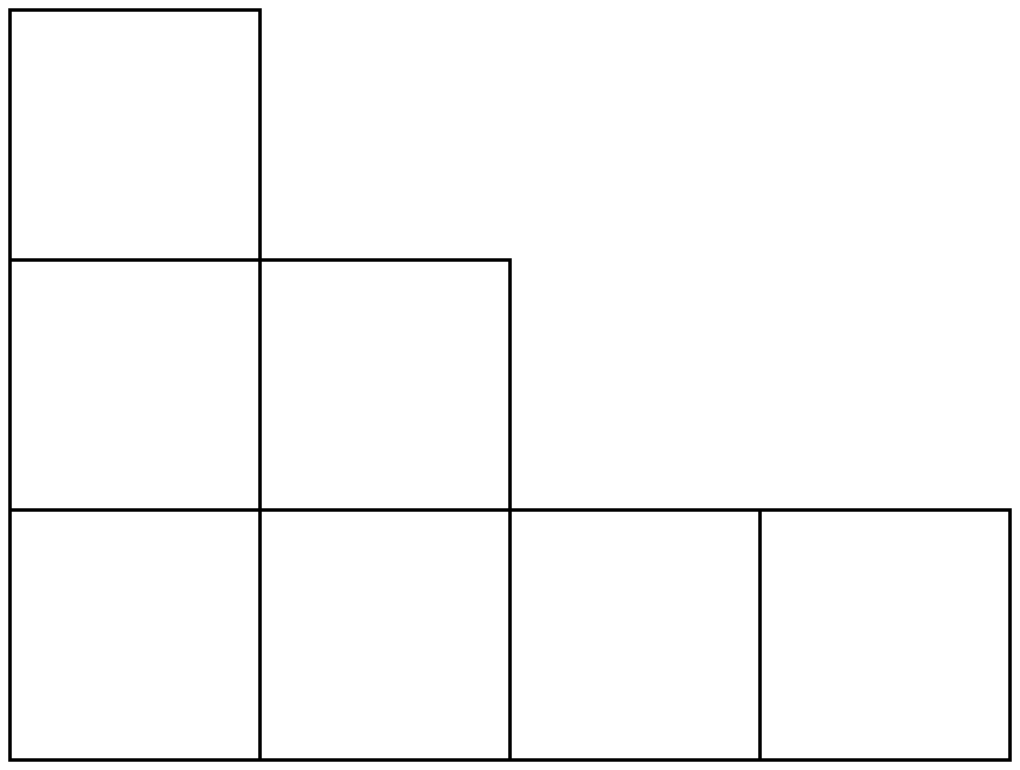}
\end{equation*}
is the Young diagram corresponding to the partition $(4,2,1)$.
The nonzero elements $\la_i$ are called the parts of $\la$. The height of 
$\la$ is the number of parts and its width equals the largest part. 
At times it is convenient to denote a partition $\la$ with $L_i$ parts
equal to $i$ by $\la=(1^{L_1}2^{L_2}\cdots)$.
The partition $\la^{\top}$ is the partition corresponding to the transposed 
diagram of $\la$ obtained by reflecting along the diagonal, i.e., if
$\la=(1^{L_1}2^{L_2}\cdots N^{L_N})$ then $(\la^{\top})_i=L_i+\cdots +L_N$.
The addition $\la+\mu$ of the partitions $\la$ and $\mu$ is defined by the 
addition of their parts $(\la+\mu)_i=\la_i+\mu_i$. The parts of the 
partition $\la\cap \mu$ are given by $(\la\cap\mu)_i=\min\{\la_i,\mu_i\}$,
and $\la/\mu$ denotes the skew shape obtained by removing the boxes
of $\mu$ from $\la$. 
By $\la\geq \mu$ in dominance order we mean $\la_1+\cdots+\la_i\geq 
\mu_1+\cdots+\mu_i$ for all $i$. The set of ``rectangular'' 
partitions (i.e., partitions with rectangular diagram) is denoted by $\R$.

In this paper we will often encounter arrays of rectangular partitions.
For such an array $\mu=(\mu_1,\ldots,\mu_L)\in\R^L$, define
$\mu^{\star}=(\mu_1^{\top},\ldots,\mu_L^{\top})$ and 
$|\mu|=|\mu_1|+\cdots+|\mu_L|$. When all components of
$\mu$ have height 1 ordered according to decreasing width, i.e., 
$\mu_i=(k_i)$ with $k_1\geq \cdots\geq k_L$, one can identify $\mu$ with the 
partition $(k_1,\ldots,k_L)$. Notice, however, that 
$\mu^{\star}\neq \mu^{\top}$ in this case. 
There is the following partial order on $\R^L$ modulo reordering. 
Define $\la^{(a)}$ as the 
partition obtained from $\la\in\R^L$ by putting the widths of all components 
of $\la$ of height $a$ in decreasing order. Then $\la\geq \mu$ for 
$\la,\mu\in\R^L$ if $\la^{(a)}\geq \mu^{(a)}$ for all $a$ by the dominance 
order on partitions.

Next we consider Young tableaux.
Let $X=\{x_1<x_2<\cdots <x_n\}$ be a totally ordered alphabet of 
non-commutative indeterminates. A Young tableau over $X$ is a filling of a 
Young diagram such that each row is weakly increasing from left to right and
 each column is strictly increasing from bottom to top. The Young diagram 
(or, equivalently, partition) underlying a Young tableau $T$ is called the 
shape of $T$ and the height of $T$ is the height of its shape. The content of 
a Young tableau $T$ is an array $\mu=(\mu_1,\ldots,\mu_n)$ where $\mu_i$ is 
the number of boxes filled with $x_i$. The set of all Young tableaux
of shape $\la$ and content $\mu$ is denoted by $\tab{\la}{\mu}$. 
It is clear that $\tab{\la}{\mu}=\emptyset$ unless $|\la|=|\mu|$.

Young tableaux can also be represented by words over the alphabet $X$.
Let $\X$ be the free monoid generated by $X$. 
By the Schensted bumping algorithm~\cite{S61}
one can associate a Young tableau to each word $w\in\X$ denoted by $[w]$. 
Knuth~\cite{K70} introduced equivalence relations on words generated by
\begin{equation}
\label{knuth}
\begin{aligned}
zxy &\equiv xzy \qquad (x\leq y < z),\\
yxz &\equiv yzx \qquad (x < y \leq z),
\end{aligned}
\end{equation}
for $x,y,z\in X$ and showed that $[w]=[w']$ if and only if $w \equiv w'$. 
The word $w_T$ obtained from a Young tableau $T$ by reading its entries
successively from left to right down the page is called a word in 
row-representation or, for short, a row-word. 
Since $T=[w_T]$, the Schensted and row-reading algorithms provide an
one-to-one correspondence between the plactic monoid $\X/\equiv$ and 
the set of Young tableaux over $X$.
Using the above correspondence, we say that a word has shape $\la$ and 
content $\mu$ if the corresponding Young tableau $[w]$ is in $\tab{\la}{\mu}$.
Furthermore, the product of two Young tableaux $S$ and $T$
is defined as $S\cdot T=[w_S w_T]$ where $w_S w_T$ is the word
formed by juxtaposing the row-words $w_S$ and $w_T$. 

Finally, the following definitions for words are needed.
A word $w=w_1 w_2 \cdots w_k$ with $w_i\in X$ is called a Yamanouchi 
word if, for all $1\leq i\leq k$, the sequence $w_k \cdots w_i$ contains at 
least as many $x_1$ as $x_2$, at least as many $x_2$ as $x_3$ and so on.
We call a word balanced if all of the letters in $X$ occur an equal number of
times. Let $w$ be a word on the two-letter alphabet $\{x<y\}$ and recursively
connect all pairs $yx$ in $w$ as in the following example:
\begin{equation*}
w=\contract{y}{\contract{y}{\contract{y}{}{x}\contract{y}{}{x}}{x}}{x}y.
\end{equation*}
A pair $\contract{y}{\cdots}{x}$ is called an inverted pair. 
All letters of $w$ which do not belong to an inverted pair are 
called non-inverted, and the subword of $w$ containing its non-inverted
letters is of the form $x^r y^s$.

\subsection{Kostka polynomials}
\label{sec_Kostkaintro}

Throughout this paper $x^{\la}:=x_1^{\la_1}\cdots x_n^{\la_n}$ where
$x_1,\dots,x_n$ are commutative variables (not to be confused with the
{\em non}commutative letters in the alphabet $X$) 
and $\la=(\la_1,\dots,\la_n)$.
The Schur polynomial $s_{\la}$ in the variables $x_1,\dots,x_n$ is defined as
\begin{equation}\label{schur}
s_{\la}(x)=\sum_{T\in\tab{\la}{\cdot}}x^T,
\end{equation}
where $x^T:=x^{\content(T)}$.
The Kostka polynomials $K_{\lm}(q)$ arise as the connection coefficients
between the Schur and Hall--Littlewood polynomials~\cite{M95},
\begin{equation}\label{HL}
s_{\la}(x)=\sum_{\mu\vdash |\la|} K_{\lm}(q)P_{\mu}(x;q).
\end{equation}
Here $\la$ and $\mu$ are partitions and $K_{\lm}(q)\neq 0$ if and only if
$|\la|=|\mu|$ and $\la\ge \mu$.

A combinatorial interpretation of the Kostka polynomials was obtained by
Lascoux and Sch\"utzenberger~\cite{LS78}, who showed that
\begin{equation}\label{LSKostka}
K_{\lm}(q)=\sum_{T\in\tab{\la}{\mu}} q^{c(T)},
\end{equation}
where $c(T)$ is the charge of a Young tableau defined below.

Let $T\in\tab{\cdot}{\mu}$ be a Young tableau of content $\mu$ over 
$X=\{x_1<x_2<\cdots<x_n\}$ and let 
$T_{\min}=T_{\min}(\mu):= \left[x_1^{\mu_1} \cdots x_n^{\mu_n}\right]$ 
be the one-row tableau of content $\mu$. When $T\neq T_{\min}$ and $w_T=x_i u$ 
the initial cyclage $\C$ on $T$ is defined as $\C(T)=[u x_i]$. The cocharge 
$\co(T)$ of $T\in\tab{\cdot}{\mu}$ 
is the number of times one has to apply $\C$ to obtain $T_{\min}$. 
The charge of $T$ is defined as $c(T)=\|\mu\|-\co(T)$ where $\|\mu\|$ is 
the cocharge of the Young tableau 
$T_{\max}:=\left[x_n^{\mu_n}\cdots x_1^{\mu_1}\right]$, given by
$\|\mu\|=\sum_{i<j}\min\{\mu_i,\mu_j\}$.

To illustrate the above definitions take, for example,
$T=[x_3x_2x_1^2x_2]$. Then
\begin{equation*}
\begin{aligned}
\C(T)&=[x_2x_1^2x_2x_3],\\
\C^2(T)&=[x_1^2x_2x_3x_2]=[x_3x_1^2x_2^2],\\
\C^3(T)&=[x_1^2x_2^2x_3],
\end{aligned}
\end{equation*}
so that $\co(T)=3$ and $c(T)=1$ in this example.

Another combinatorial description of the Kostka polynomials, in terms of
rigged configurations, is due to Kirillov and Reshetikhin~\cite{KR88}
and provides an explicit formula for calculating the Kostka polynomials as
\begin{equation}\label{KRKostka}
K_{\lm}(q)=\sum_{\alpha} q^{C(\alpha)} \prod_{a,i\geq 1} 
\qbin{P_i^{(a)}(\alpha)+\alpha_i^{(a)}-\alpha_{i+1}^{(a)}}
{\alpha_i^{(a)}-\alpha_{i+1}^{(a)}}.
\end{equation}
The summation is over sequences $\alpha=(\alpha^{(0)},\alpha^{(1)},\dots)$ of
partitions such that $\alpha^{(0)}=\mu^{\top}$ and 
$|\alpha^{(a)}|=\la_{a+1}+\la_{a+2}+\cdots$. Furthermore
\begin{align}
P_i^{(a)}(\alpha)&=
\sum_{k=1}^i(\alpha_k^{(a-1)}-2\alpha_k^{(a)}+\alpha_k^{(a+1)})\\
\intertext{and}
C(\alpha)&=\sum_{a,i\geq 1}\binom{\alpha_i^{(a-1)}-\alpha_i^{(a)}}{2},
\end{align}
where $\tbinom{a}{2}=a(a-1)/2$ for $a\in\Integer$.
Expressions of the type~\eqref{KRKostka} are often called 
fermionic as they can be interpreted as the partition 
function for a system of quasi-particles with fractional statistics
obeying Pauli's exclusion principle~\cite{KKMM93a,KKMM93b}.

In section~\ref{sec_KP} a third combinatorial representation of the 
Kostka polynomials as the generating function of paths will
be discussed. This representation is due to Lascoux, Leclerc and 
Thibon~\cite{LLT95} and Nakayashiki and Yamada~\cite{NY95} and is the 
starting point for our generalized Kostka polynomials.
As we will see in subsequent sections, these generalized Kostka polynomials 
also admit representations stemming from equations~\eqref{LSKostka} 
and~\eqref{KRKostka}.

\section{A$_{n-1}$ supernomials and generalized Kostka polynomials}
\label{sec_SK}

This section deals with paths defined as ordered sequences of rectangular 
Young tableaux. Assigning weights to the paths, we consider the generating 
functions over two different sets of paths called unrestricted and 
classically restricted. These are treated in 
sections~\ref{sec_super} and~\ref{sec_KP}, respectively.
As will be shown in section~\ref{sec_propSK}, the generating functions 
over the set of unrestricted paths are A$_{n-1}$ generalizations 
of the A$_1$ supernomials \eqref{qsuper}.
The generating functions over the set of classically restricted paths lead to 
generalizations of the Kostka polynomials.

\subsection{Unrestricted paths and A$_{n-1}$ supernomials}
\label{sec_super}

Denote by $\B_{\la}$ the set $\tab{\la}{\cdot}$ of
Young tableaux of shape $\la$ over the alphabet $\{1,2,\cdots,n\}$. 
An element of $\B_{\la}$ is called a step and
an ordered sequence of $L$ steps is a path of length $L$
denoted by $p_L\otimes\cdots\otimes p_1$. 
We treat here only paths with rectangular steps $p_i$, i.e., $p_i\in\B_{\mu_i}$
for $\mu_i\in\R$. Let us however emphasize that the steps in a path can
have different shapes indicated by the subscript $i$ on $\mu_i$.
Paths with this property are called inhomogeneous~\cite{NY95}.

The reason for the tensor product notation for paths (treated
here as ordered sequences of steps only) is for notational
convenience, but is motivated by the relation to the theory of crystal
bases~\cite{K91}. In this setting $\B_{(i^a)}$ is usually labelled by
$\B_{i\La_a}$ where $\La_a$ are the fundamental weights
of A$_{n-1}$. The set $\B_{i\La_a}$ is called a perfect crystal and 
parametrizes a basis of the irreducible highest weight module of A$_{n-1}$
with highest weight $i\La_a$~\cite{KN94}. There exist crystal bases for all
integrable highest weight modules and they are compatible with the tensor
product structure.

\begin{definition}[Unrestricted paths]
\label{def_paths}
For fixed integers $n\geq 2$ and $L\geq 0$ let $\la\in \Int^n$ and 
$\mu=(\mu_1,\ldots,\mu_L) \in \R^L$. The set of paths $\P_{\lm}$ is 
defined as
\begin{equation}
\P_{\lm}=\{p_L\otimes\cdots\otimes p_1|~p_i\in\B_{\mu_i} \text{ and }
\sum_{i=1}^L \content(p_i)=\la \}.
\end{equation}
\end{definition}

To each path $P\in \P_{\lm}$ we assign an energy $h(P)\in \Int$ as
\begin{equation}
\label{h}
h(P)=\sum_{i=1}^{L-1} i h(p_{i+1}\otimes p_i),
\end{equation}
where $h(p\otimes p')$ for the steps $p\in\B_{\nu}$ and 
$p'\in\B_{\nu'}$ is defined as the number of boxes in 
the product $p\cdot p'$ that lie outside the Young diagram 
$\nu+\nu'$ or, more formally, as
\begin{equation}
\label{h_expl}
h(p\otimes p')=|\nu+\nu'|-|\shape(p\cdot p')\cap (\nu+\nu')|.
\end{equation}

\begin{example}
\begin{equation*}
\text{Let}\quad
P=p_2\otimes p_1=
\thickspace \raisebox{-0.4cm}{\epsfxsize=1 cm \epsffile{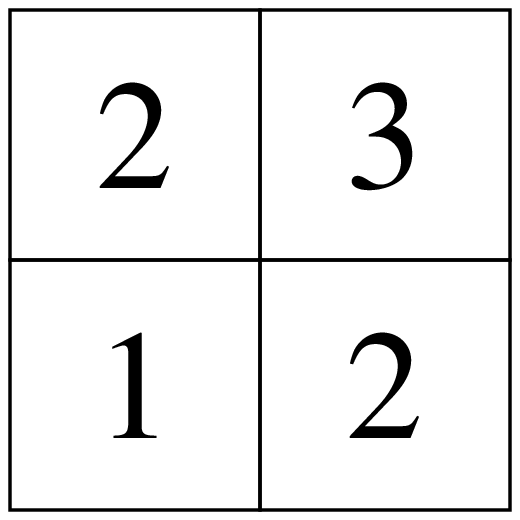}}
\thickspace\otimes\thickspace 
\raisebox{-0.4cm}{\epsfxsize=1 cm \epsffile{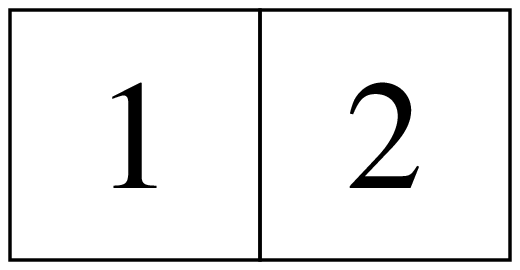}}\; .
\quad\text{Then}\quad
p_2 \cdot p_1=\thickspace
\raisebox{-0.4cm}{\epsfxsize=1.5 cm \epsffile{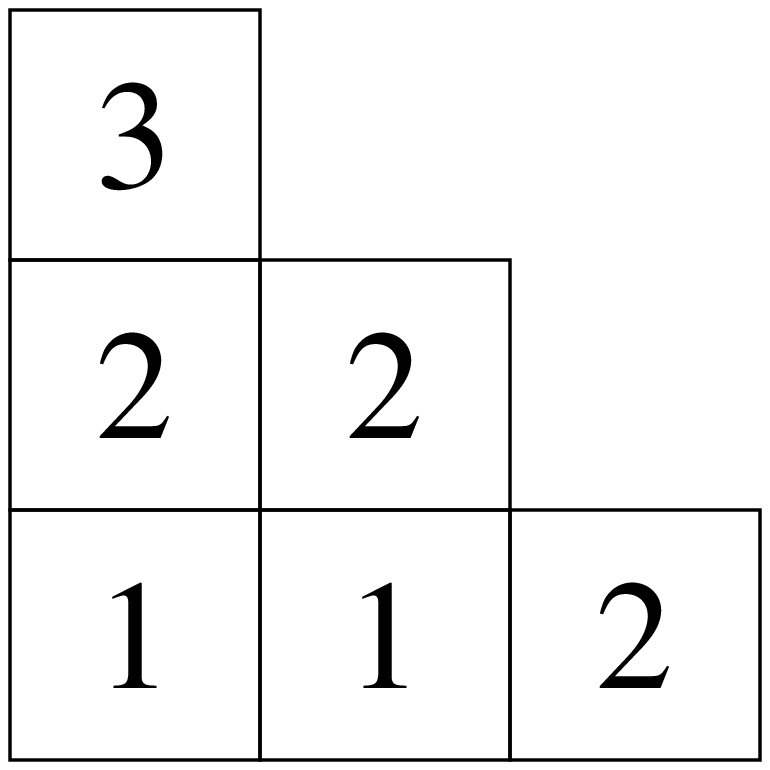}}
\end{equation*}
and $\shape(p_2\cdot p_1)=(3,2,1)$. Hence $h(P)=|(4,2)|-|(3,2,1)\cap (4,2)|
=6-5=1$.
\end{example}

The cardinality $S_{\lm}$ of $\P_{\lm}$ does not depend on the 
ordering of $\mu$, i.e.,
\begin{equation}
S_{\lm}=S_{\la\tilde{\mu}},
\end{equation}
where $\tilde{\mu}$ is a permutation of $\mu$. In general the generating 
function of $\P_{\lm}$ with paths weighted by the energy function $h$ does not
have this symmetry. To obtain a weight function such that the
resulting generating function does respect this symmetry we introduce an 
isomorphism $\sigma: \B_{\alpha} \otimes \B_{\alpha'} \to \B_{\alpha'} 
\otimes \B_{\alpha}$ for $\alpha,\alpha' \in \R$ between two successive steps. 
Let $p\otimes p' \in \B_{\alpha} \otimes \B_{\alpha'}$. Then
$\sigma(p \otimes p')=\tilde{p}' \otimes \tilde{p}$, where
$\tilde{p}'$ and $\tilde{p}$ are the unique Young tableaux of
shape $\alpha'$ and $\alpha$, respectively, which satisfy
\begin{equation}
\label{iso}
p \cdot p'=\tilde{p}' \cdot \tilde{p}.
\end{equation}
The uniqueness of the Young tableaux $\tilde{p}'$ and $\tilde{p}$ is ensured
since the Littlewood--Richardson coefficients have the symmetry 
$c_{\alpha\alpha'}^{\beta}=c_{\alpha'\alpha}^{\beta}$ and 
for rectangular shapes $\alpha$ and $\alpha'$ obey 
$c_{\alpha\alpha'}^{\beta}\leq 1$.
Notice that $\sigma$ is the identity if $p$ and $p'$ have the same shape.

\begin{definition}[Isomorphism]
\label{def_iso}
For a path $P=p_L \otimes \cdots \otimes p_1\in\P_{\lm}$ we define the 
isomorphism $\sigma_i$ as
\begin{equation}
\sigma_i(P)=p_L\otimes \cdots \otimes \sigma(p_{i+1} \otimes p_i)
\otimes \cdots \otimes p_1.
\end{equation}
\end{definition}

The group generated by the isomorphisms $\sigma_i$ is the symmetric group,
i.e., $\sigma_i^2=\Id$, $\sigma_i \sigma_{i+1} \sigma_i= \sigma_{i+1} \sigma_i
\sigma_{i+1}$ and $\sigma_i \sigma_j= \sigma_j \sigma_i$ for $|i-j|\geq 2$.
The proof of the braiding relation is non-trivial (see~\cite{S98a,S98c}).

\begin{definition}[Orbit]
\label{def_orbit}
The set $\O_P$ is the orbit of the path $P\in\P_{\lm}$ under the group 
generated by the isomorphisms $\sigma_i$.
\end{definition}

The weight of a path $P$ is now given by the mean of the energy function
$h$ over the orbit of $P$.

\begin{definition}[Weight]
\label{def_weight}
For $P\in\P_{\lm}$, the weight function $H:\P_{\lm}\to \Int$ is defined as
\begin{equation}
\label{H}
H(P)=\frac{1}{|\O_P|} \sum_{P'\in \O_P} h(P').
\end{equation}
\end{definition}

It is not obvious from~\eqref{H} that the weight $H(P)$ of a path $P$
is indeed integer. This will follow from theorem~\ref{theo_Hco}.
Before we continue to define the generating functions over the set of
paths $\P_{\lm}$, some remarks on the relation of our definitions to
lattice paths of exactly solvable lattice models and the theory of crystal
bases are in order.

\begin{remark}
\label{rem_ref}
For homogeneous paths, i.e., $P\in\P_{\lm}$ with $\mu_1=\cdots =\mu_L$, 
the weight simplifies to $H(P)=h(P)$ which is the weight function
of configuration sums of A$^{(1)}_{n-1}$ solvable lattice models.
For example, for $p,p'\in\B_{(N)}$, the energy function $h(p\otimes p')$
coincides with the one of refs.~\cite{JMO88,DJKMO89} (and references 
therein) given by
\begin{equation}
\label{h_sym}
h(p\otimes p')=\max_{\tau\in S_N}\bigl\{
\sum_{i=1}^N \chi(p_i> p'_{\tau_i})\bigr\}.
\end{equation}
Here $p_i,p'_i\in\{1,2,\ldots,n\}$ are the letters in $p=[p_1\cdots p_N]$ and
$p'=[p'_1\cdots p'_N]$, $S_N$ is the permutation group on $1,2,\ldots,N$, 
$\chi({\rm true})=1$ and $\chi({\rm false})=0$.
An alternative combinatorial expression of~\eqref{h_sym} in terms of 
so-called nonmovable tableaux is given in~\cite{KKN97}.
When $p,p'\in\B_{(1^N)}$, our energy function reduces to 
$h(p\otimes p')=\min_{\tau\in S_N}\bigl\{\sum_{i=1}^N \chi(p_i> p'_{\tau_i})
\bigr\}$ of ref.~\cite{O97}.

Nakayashiki and Yamada~\cite{NY95} defined weight functions on
inhomogeneous paths when either $\mu$ or $\mu^{\star}$ is a partition, i.e.,
when $|\mu_1|\geq \cdots \geq |\mu_L|$ and either
$\height(\mu_i)=1$ for all $i$ or $\width(\mu_i)=1$ for all $i$.
Their isomorphism, defined in terms of graphical rules (Rule~3.10 and~3.11 
of ref.~\cite{NY95}), is a special case of the isomorphism $\sigma$.
The expression for $H(P)$ that they give is quite different to that of 
equation~\eqref{H} even though it is the same function for the subset of 
paths they consider. For example when $\height(\mu_i)=1$ for all $i$, 
$H$ of ref.~\cite{NY95} is, in our normalization, given by
\begin{equation*}
H(P)=\sum_{i=2}^L\sum_{j=1}^{i-1} h(p_i\otimes p_j^{(i-1)}),
\end{equation*}
where $P=p_L\otimes\cdots \otimes p_1\in\P_{\lm}$, $p_i^{(i)}=p_i$ and 
$p_j^{(i)}=p'_i$ with $P'=\sigma_{i-1}\circ\sigma_{i-2}\circ\cdots\circ
\sigma_j(P)$ for $j<i$.

Lascoux, Leclerc and Thibon~\cite{LLT95} defined a weight function $b(T)$ 
for Young tableaux $T$ as the mean over certain orbits very similar in spirit 
to equation~\eqref{H} (see theorem 5.1 in ref.~\cite{LLT95}). In fact, when 
$\height(\mu_i)=1$ for all $i$, each path $P\in\P_{\cdot\mu}$ can 
be mapped to a Young tableau $T\in\tab{\cdot}{\mu}$ (by the virtue of
the map $\omega$ of equation~\eqref{omega} below, i.e., $T=[\omega(P)]$),
and in this case one finds that $H(P)=\|\mu\|-b(T)$ where we recall that
$\|\mu\|=\sum_{i<j}\min\{\mu_i,\mu_j\}$.
\end{remark}

\begin{remark}
Kashiwara~\cite{K91} defined raising and lowering operators $f_i$ and $e_i$ 
$(0\leq i\leq n-1)$ acting on elements of a crystal $\B_{i\La_a}$.
Set $\B_k=\B_{i_k\La_{a_k}}$ $(k=1,2)$. Then for $p_1\in\B_1$ and $p_2\in\B_2$
the lowering operators act on the tensor product $p_2\otimes p_1$ as follows
\begin{equation*}
e_i(p_2\otimes p_1)=\begin{cases}
p_2\otimes e_i p_1 & \text{if $\varphi_i(p_1)\geq \varepsilon_i(p_2)$},\\
e_i p_2\otimes p_1 & \text{if $\varphi_i(p_1)< \varepsilon_i(p_2)$},
\end{cases}
\end{equation*}
where $\varepsilon_i(b)=\max\{k|e_i^k(b)\neq 0\}$ and 
$\varphi_i(b)=\max\{k|f_i^k(b)\neq 0\}$. The action of $f_i$ on
a tensor product is defined in a similar way.
(Note that to conform with the rest of this paper the order of the tensor 
product is inverted in comparison to the usual definitions).
Let $e_i(p_2\otimes p_1)\neq 0$. Up to an additive constant
the energy function used in crystal theory is recursively defined as 
\begin{equation}
\label{h_crystal}
E(e_i(p_2\otimes p_1))=\begin{cases}
E(p_2\otimes p_1)+1 & 
 \text{if $i=0$ and $\varphi_0(p_1)\geq \varepsilon_0(p_2)$,
 $\varphi_0(\tilde{p}_2)\geq \varepsilon_0(\tilde{p}_1)$},\\
E(p_2\otimes p_1)-1 & 
 \text{if $i=0$ and $\varphi_0(p_1)< \varepsilon_0(p_2)$,
 $\varphi_0(\tilde{p}_2)<\varepsilon_0(\tilde{p}_1)$},\\
E(p_2\otimes p_1) & \text{otherwise}. \end{cases}
\end{equation}
Here $p_2\otimes p_1\mapsto \tilde{p}_1\otimes \tilde{p}_2$ with
$p_1,\tilde{p}_1\in\B_1$ and $p_2,\tilde{p}_2\in\B_2$ is an isomorphism
obeying certain conditions~\cite{KKMMNN92a,KKMMNN92b,K94}. 
The isomorphism $\sigma$ 
defined through~\eqref{iso} yields the isomorphism of crystal theory. 
Up to a sign the energy $h(p_2\otimes p_1)$ as defined in~\eqref{h_expl}
provides an explicit expression for the recursively defined energy 
of~\eqref{h_crystal}, i.e., $E(p_2\otimes p_1)=-h(p_2\otimes p_1)$. We do 
not prove these statements in this paper.
\end{remark}

Let us now define the A$_{n-1}$ supernomial as the generating function of the 
set of paths $\P_{\lm}$ weighted by $H$ of definition~\ref{def_weight}.

\begin{definition}[Supernomials]
\label{def_super}
Let $\la \in \Int^n$ and $\mu\in \R^L$. Then the supernomial
$S_{\lm}(q)$ is defined as
\begin{equation}
\label{super}
S_{\lm}(q)= \sum_{P\in \P_{\lm}} q^{H(P)}.
\end{equation}
\end{definition}

Since $H(P)=H(P')$ for $P' \in \O_P$ it is clear that
\begin{equation}
\label{step_order}
S_{\lm}(q)=S_{\la\tilde{\mu}}(q),
\end{equation}
where $\tilde{\mu}$ is a permutation of $\mu$.

To conclude this section we comment on the origin of the
terminology supernomial as first introduced in ref.~\cite{SW97}.
Recalling definition~\ref{def_paths} 
of the set of paths $\P_{\lm}$ and equation~\eqref{schur} for
the Schur polynomial, one finds that 
\begin{equation}\label{expansion}
s_{\mu_1}(x)\cdots s_{\mu_L}(x) 
= \sum_{\la\vdash |\mu|}S_{\lm} x^{\la},
\end{equation}
where $S_{\lm}:=S_{\lm}(1)=|\P_{\lm}|$.
For homogeneous paths, i.e., $\mu_1=\dots=\mu_L$,
this is the usual definition for various kinds of multinomial coefficients.
The supernomials (or, more precisely, $q$-supernomial coefficients) 
can thus be viewed as $q$-deformations of generalized multinomial 
coefficients.

\subsection{Classically restricted paths and generalized Kostka polynomials}
\label{sec_KP}

Analogous to the previous section we now introduce classically restricted
paths and their generating function. To describe the set of classically
restricted paths we first map paths onto words and then specify the
restriction on these words.

For our purposes it will be most convenient to label the alphabet underlying
the words associated to paths as
\begin{equation}
\label{alph_X}
X^a=\{ x_1^{(1)}< x_1^{(2)} < \cdots < x_1^{(a_1)}<
       x_2^{(1)}< \cdots < x_2^{(a_2)}<
       \cdots<
       x_L^{(1)}< \cdots < x_L^{(a_L)} \}
\end{equation}
for some fixed integers $1\leq a_i\leq n$. As before, $\X^a$ denotes the free
monoid generated by $X^a$.
The $i$-subword of a word $w \in \X^a$ is the subword of $w$ consisting of the 
letters $x_i^{(j)}$ $(1\leq j\leq a_i)$ only.
More generally, the $(i_1,\ldots,i_{\ell})$-subword of $w$ is the subword
consisting of the letters with subscripts $i_1,\ldots,i_{\ell}$ only.

We are interested in the following subset of $\X^a$
\begin{equation}
\label{W}
\W=\{ w\in \X^a | \text{ each $i$-subword of $w$ is a balanced Yamanouchi 
word} \}.
\end{equation}
By the Schensted bumping algorithm each word $w\in\W$
corresponds to a Young tableau $[w]$ over the alphabet $X^a$. It is an 
easy matter to show that if $w\in\W$ is Knuth equivalent to $w'\in\X^a$, 
then $w'\in \W$. Hence it makes sense to consider the set of Young tableaux 
over $X^a$ corresponding to $\W/\equiv$. 
Instead of labelling the content of such a tableau 
by $\mu=(\mu_1^{(1)},\ldots,\mu_L^{(a_L)})$ (where $\mu_i^{(j)}$ is the
number of $x_i^{(j)}$'s), we set $\mu=(\mu_1,\ldots,\mu_L)$
where $\mu_i=(\mu_i^{(1)},\ldots,\mu_i^{(a_i)})$ so that 
$\mu\in\R^L$ with $\height(\mu_i)=a_i$. The set of words $w\in\W$
with $\content([w])=\mu$ is denoted $\W_{\mu}$.
\begin{definition}[Littlewood--Richardson tableaux]
\label{def_Stab}
Let $L\geq 0$ and $n\geq 2$ be integers, $\la$ a partition and $\mu\in\R^L$.
Then the set of LR tableaux of shape $\la$ and content $\mu$ is defined as
\begin{equation}
\ttab{\la}{\mu}=\{ T | w_T\in\W_{\mu}\text{ and } \shape(T)=\la \}.
\end{equation}
\end{definition}
The set of LR tableaux $\ttab{\la}{\mu}$ reduces to the set $\tab{\la}{\mu}$ 
of Young tableaux over $X=\{x_1<\cdots<x_L\}$ when $a_1=a_2=\cdots=a_L=1$ 
in~\eqref{alph_X}.

We now define a map
\begin{equation}
\label{omega}
\omega: \P_{\lm} \to \W_{\mu}
\end{equation}
in the following way. Let $P=p_L\otimes\cdots\otimes p_1\in\P_{\lm}$
and let $(j,k)$ denote the $k$th row of $p_j$. Set $P'=P$, $w$ to the empty 
word and carry out the following procedure $|\mu|$ times:
\begin{quote}
If $(j,k)$ labels the position of the rightmost, maximal entry in $P'$,
obtain a new $P'$ by removing this maximal entry from $P'$ and
append $x_j^{(k)}$ to $w$.
\end{quote}\noindent
The resulting word $w$ defines $\omega(P)$. Equivalently, $[\omega(P)]$
is the column insertion recording tableaux of $\col(p_L)\ldots\col(p_1)$
where $\col(T)$ is the column word of the tableau $T$.

The word $\omega(P)$ obtained via the above procedure is indeed in
$\W_{\mu}$. The Yamanouchi condition is guaranteed in each intermediate
$w$ in the construction thanks to the fact that all steps $p_j$ are
Young tableaux. Since all $p_j$ have rectangular shape $\omega(P)$
is a balanced Yamanouchi word.
Note that the integer $a_j$ in the alphabet~\eqref{alph_X} used in
definition~\eqref{W} of $\W$ is exactly the height of step $p_j$. 

\begin{example}
\label{ex_omega}
To illustrate the map $\omega$ take for example
\begin{equation*}
P=\; \raisebox{-0.4cm}{\epsfxsize=1 cm \epsffile{fig_ex1.ps}}\; \otimes
  \; \raisebox{-0.4cm}{\epsfxsize=1 cm \epsffile{fig_ex2.ps}}\; \otimes
  \; \raisebox{-0.4cm}{\epsfxsize=0.5 cm \epsffile{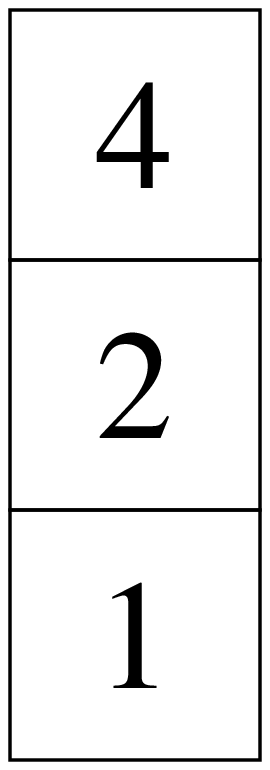}}\; .
\end{equation*}
Then $\omega(P)=x_1^{(3)} x_3^{(2)} x_1^{(2)} x_2^{(1)} x_3^{(1)}
x_3^{(2)} x_1^{(1)} x_2^{(1)} x_3^{(1)}$ and the corresponding 
LR tableau (identifying $x_j^{(k)}$ with $j^{(k)}$) is
\begin{equation}
\label{gyt}
\epsfxsize=2 cm \epsffile{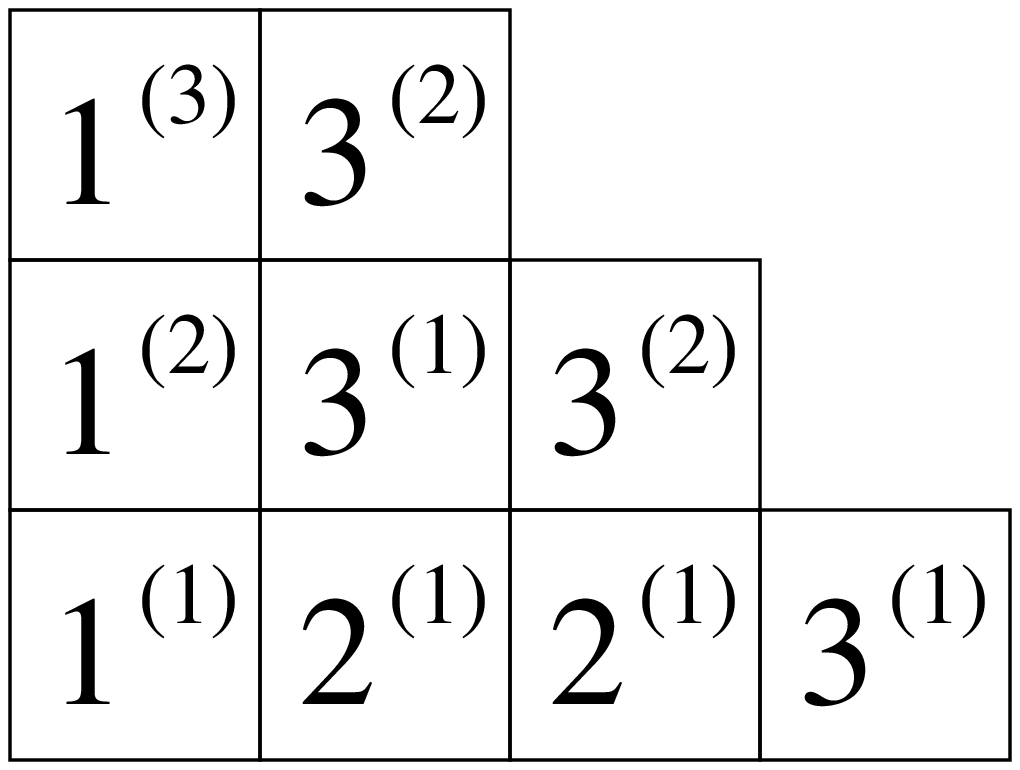}\; .
\end{equation}
\end{example}

The map $\omega$ has the following properties which are needed later.
\begin{lemma} \label{lem_om}
Let $P=p_L\otimes\cdots\otimes p_1\in\P_{\lm}$. Then:
\begin{enum}
\item \label{i_1}
The shape of $p_{i+1}\cdot p_i$ is the same as the shape of the
$(i,i+1)$-subword of $\omega(P)$.
\item \label{i_2}
Let $w=w_1\ldots w_{|\mu|}:=\omega(P)$ and
$w'=w'_1\ldots w'_{|\mu|}:=\omega(\sigma_i(P))$. 
Then $\shape(w_{\ell}\ldots w_{\ell'})=
\shape(w'_{\ell}\ldots w'_{\ell'})$ for all
$1\leq \ell\leq \ell'\leq |\mu|$. 
In particular,
$\omega(P)$ and $\omega(\sigma_i(P))$ have the same shape, and if
$\omega(P)$ is in row-representation then so is $\omega(\sigma_i(P))$.
\end{enum}
\end{lemma}

\begin{proof}
\ref{i_1}~Let $L(w,k)$ be the largest possible sum of the lengths of $k$ 
disjoint increasing sequences extracted from the word $w$. If $w$ has shape
$\nu$ then $L(w,k)=\nu_1+\cdots+\nu_k$ 
(see for example lemma 1 on page 32 of ref.~\cite{F97}). 
Since $\shape(p_{i+1}\cdot p_i)=\shape(w_{p_{i+1}} w_{p_i})$
\ref{i_1} is proven if we can show that 
$L(w_{p_{i+1}} w_{p_i},k)=L(\omega(p_{i+1}\otimes p_i),k)$ for all $k$.
Let $s_1,\ldots,s_k$ be disjoint increasing sequences of 
$w_{p_{i+1}} w_{p_i}$ such that the sum of their lengths is 
$L(w_{p_{i+1}} w_{p_i},k)$. Now $\omega$ successively maps the 
rightmost largest element of $P$ to $x_{\ell}^{(m)}$ if its position is in 
step $\ell$ and row $m$. Interpreting $s_1,\ldots,s_k$ as decreasing sequences
from right to left we see that the image of $s_j$ under $\omega$ is an
increasing sequence from left to right in $\omega(p_{i+1}\otimes p_i)$.
Conversely, each increasing sequence of $\omega(p_{i+1}\otimes p_i)$ 
is a decreasing sequence from right to left of $w_{p_{i+1}} w_{p_i}$ 
which proves $L(w_{p_{i+1}} w_{p_i},k)=L(\omega(p_{i+1}\otimes p_i),k)$.

\ref{i_2}~Since $\sigma_i$ only changes the $(i,i+1)$-subword of $\omega(P)$ 
it suffices to prove~\ref{i_2} for paths $P=p\otimes p'$ of length two. 
Set $\tilde{p}' \otimes \tilde{p}:=\sigma(p \otimes p')$
so that $p\cdot p'=\tilde{p}'\cdot \tilde{p}$. In particular $p\cdot p'$ and 
$\tilde{p}'\cdot \tilde{p}$ have the same shape.
Hence from~\ref{i_1} we conclude that $\omega(p \otimes p')$
and $\omega(\sigma(p\otimes p'))$ have the same shape as well.
Denoted by $w_p^{(\ell,\ell')}w_{p'}^{(\ell,\ell')}$ and 
$w_{\tilde{p}'}^{(\ell,\ell')}w_{\tilde{p}}^{(\ell,\ell')}$
the words obtained from $w_p w_{p'}$ and $w_{\tilde{p}'}w_{\tilde{p}}$
after successively removing the $\ell-1$ rightmost biggest letters and 
the $|\mu|-\ell'$ leftmost smallest letters, respectively.  
Then the arguments of point~\ref{i_1} still go through, that is, 
$L(w_p^{(\ell,\ell')}w_{p'}^{(\ell,\ell')},k)=L(w_{\ell}\ldots w_{\ell'},k)$ 
and $L(w_{\tilde{p}'}^{(\ell,\ell')}w_{\tilde{p}}^{(\ell,\ell')},k)=
L(w_{\ell}'\ldots w_{\ell'}',k)$
for all $k$. Since $w_p^{(\ell,\ell')}w_{p'}^{(\ell,\ell')}\equiv 
w_{\tilde{p}'}^{(\ell,\ell')}w_{\tilde{p}}^{(\ell,\ell')}$, it follows that
$\shape(w_{\ell}\ldots w_{\ell'})=\shape(w'_{\ell}\ldots w'_{\ell'})$.
{}From this one can immediately deduce that $\omega(P)$ is in 
row-representation if and only if $\omega(\sigma_i(P))$ is in 
row-representation.
\end{proof}

\begin{definition}[Classically restricted paths]
\label{def_rps}
Let $\la$ be a partition such that $\height(\la)\leq n$ and let $\mu\in\R^L$. 
The set of classically restricted paths $\Pb_{\lm}$ is defined as
\begin{equation}
\Pb_{\lm}= \{ P \in \P_{\lm} |\, \shape(\omega(P))=\la\}.
\end{equation}
\end{definition}

Since each path $P\in\P_{\lm}$ contains $\la_i$ boxes filled with $i$, 
the condition that $\omega(P)$ has shape $\la$ implies that $\omega(P)$ is in 
row representation. Hence $\Pb_{\lm}$ is isomorphic to $\ttab{\la}{\mu}$.

Let us now introduce the restricted analogue of the supernomials of
definition~\ref{def_super}.

\begin{definition}
For $\la$ a partition with $\height(\la)\leq n$ and $\mu \in \R^L$ define
\begin{equation}
\label{Sb}
\Sb_{\lm}(q)=\sum_{P \in \Pb_{\lm}} q^{H(P)}.
\end{equation}
\end{definition}

When $\mu$ is a partition, i.e., $\mu\in\R^L$ such that its components $\mu_i$ 
are one-row partitions of decreasing size, $\Sb_{\lm}(q)$ reduces to the
cocharge Kostka polynomial. This follows from the work of Nakayashiki and 
Yamada~\cite{NY95} and Lascoux, Leclerc and Thibon~\cite{LLT95} and the 
relation between the weight~\eqref{H} and their statistics as explained in 
remark~\ref{rem_ref}, and is our motivation for the 
definition of generalized Kostka polynomials for all $\mu\in\R^L$.
\begin{definition}[Generalized Kostka polynomials]
\label{def_gKp}
For $\la$ a partition with $\height(\la)\leq n$ and $\mu \in \R^L$, the 
generalized Kostka polynomial $K_{\lm}(q)$ is defined as
\begin{equation}
\label{gK}
K_{\lm}(q)=q^{\|\mu\|}\Sb_{\lm}(1/q),
\end{equation}
where $\|\mu\|=\sum_{i<j}|\mu_i\cap\mu_j|$.
\end{definition}

Lemma~\ref{lem_om} ensures that if $P'\in\O_P$ with $P\in\Pb_{\lm}$ then
$\shape(\omega(P'))=\la$. Therefore $P'\in\Pb_{\la\tilde{\mu}}$ for some
permutation $\tilde{\mu}$ of $\mu$, and hence also
$\Sb_{\la\tilde{\mu}}(q)=\Sb_{\lm}(q)$ in analogy with
equation~\eqref{step_order}. We also find that
\begin{equation}
\label{step_orderK}
K_{\lm}(q)=K_{\la\tilde{\mu}}(q)
\end{equation}
since $\|\mu\|=\|\tilde{\mu}\|$.

We remark that the lattice path representation of the Kostka polynomials
was used in~\cite{KMOTU96} to express the Demazure characters as a sum over 
Kostka polynomials. It also yields expressions for the 
A$_{n-1}^{(1)}/{\rm A}_{n-1}$ branching functions in terms of the Kostka 
polynomials~\cite{NY95,KKN97} and has been employed in~\cite{HKKOTY98} to 
obtain various further A$_{n-1}^{(1)}$ branching- and string function 
identities.

\section{Initial cyclage and cocharge for LR tableaux}
\label{sec_in_cyc}

In this section we define the notions of initial cyclage, cocharge and
charge for LR tableaux which ``paves the path'' for section~\ref{sec_LS} where 
an expression of the Lascoux--Sch\"utzenberger-type~\eqref{LSKostka} for the 
generalized Kostka polynomials is derived. In the case when $\ttab{\la}{\mu}$ 
coincides with $\tab{\la}{\mu}$ our definitions reduce to the usual 
definitions of the initial cyclage etc. as mentioned in 
section~\ref{sec_Kostkaintro}.

The definition of the initial cyclage for $T\in\ttab{\cdot}{\mu}$ has to be 
altered when $a\neq (1,\dots,1)$ for the alphabet $X^a$
as given in~\eqref{alph_X}. Namely, if $T=\bigl[x_i^{(a_i)} u\bigr]
\in\ttab{\cdot}{\mu}$ with $x_i^{(a_i)}u$ in row-representation
(by the Yamanouchi condition on each of the $i$-subwords the first letter
has to be $x_i^{(a_i)}$ for some $i$), then $T':=\bigl[ux_i^{(a_i)}\bigr]$ 
obtained by cycling the first letter is not in $\ttab{\cdot}{\mu}$ since the 
Yamanouchi condition on the $i$-subword of $ux_i^{(a_i)}$ is violated if
$a_i>1$. To repair this fault we define the initial cyclage for $T=[w]$ with 
$w=x_i^{(a_i)}u$ in row-representation by considering the following chain of 
transformations
\begin{equation}\label{chain}
w\to w^{(a_i)}\to w^{(a_i-1)}\to\dots\to w^{(1)}\in\W_{\mu}.
\end{equation}
Here $w^{(a_i)}:=ux_i^{(a_i)}$ and $w^{(j)}$ for $1<j\leq a_i$ is a word such 
that ({\em i}) the $k$-subword for each $k\neq i$ is a balanced Yamanouchi 
word, ({\em ii}) the $i$-subword is balanced, ({\em iii}) the subword 
consisting of the letters $x_i^{(j)}$ and $x_i^{(j-1)}$ has one
non-inverted $x_i^{(j-1)}$ and one non-inverted $ x_i^{(j)}$
and ({\em iv}) the subword consisting of the letters $x_i^{(k)}$ and
$x_i^{(k-1)}$ with $k\neq j$ has no non-inverted letters.
The transformation $w^{(j+1)} \to w^{(j)}$ is defined as follows.
Consider the subword of $w^{(j+1)}$ consisting of the letters $x_i^{(j)}$ 
and $x_i^{(j+1)}$ only. Determine all inverted pairs for this subword and
exchange the two non-inverted letters $x_i^{(j)}$ and $x_i^{(j+1)}$
(making them an inverted pair). The resulting word is $w^{(j)}$.
Clearly, this means that $w^{(1)}\in \W_{\mu}$.

\begin{definition}[Initial cyclage]\label{def_cyc}
The initial cyclage $\C$ on $T\in\ttab{\cdot}{\mu}$ is defined as 
\begin{equation}
\C(T)=\bigl[w^{(1)}\bigr],
\end{equation}
where $w^{(1)}\in\W_{\mu}$ is the last word in the chain of 
transformations~\eqref{chain}.
\end{definition}

\begin{example}
If
\begin{equation*}
T=\; \raisebox{-0.6cm}{\epsfxsize=1.5 cm \epsffile{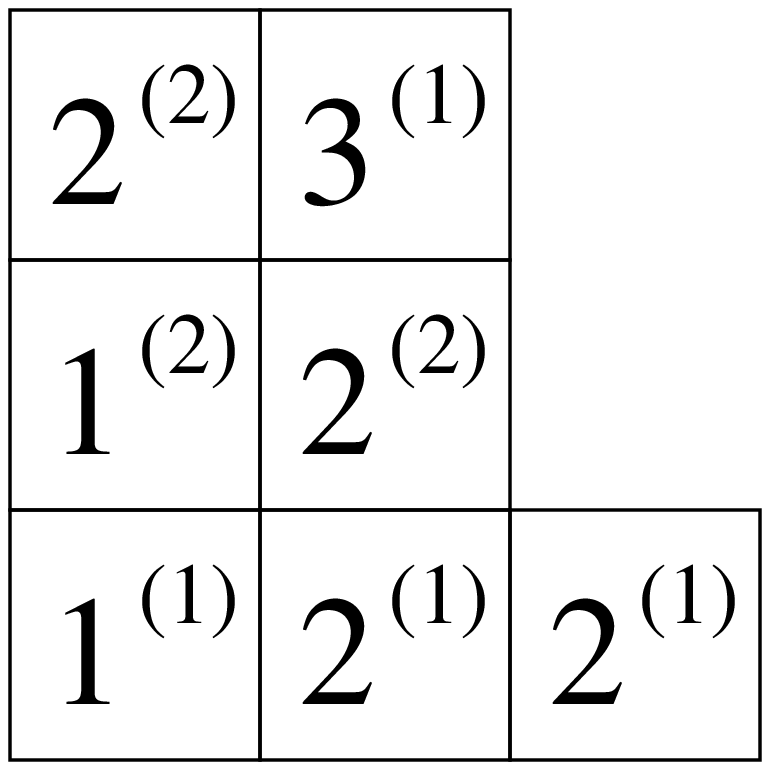}}
\quad\text{then}\quad 
\C(T)=\; \raisebox{-0.6cm}{\epsfxsize=1.5 cm \epsffile{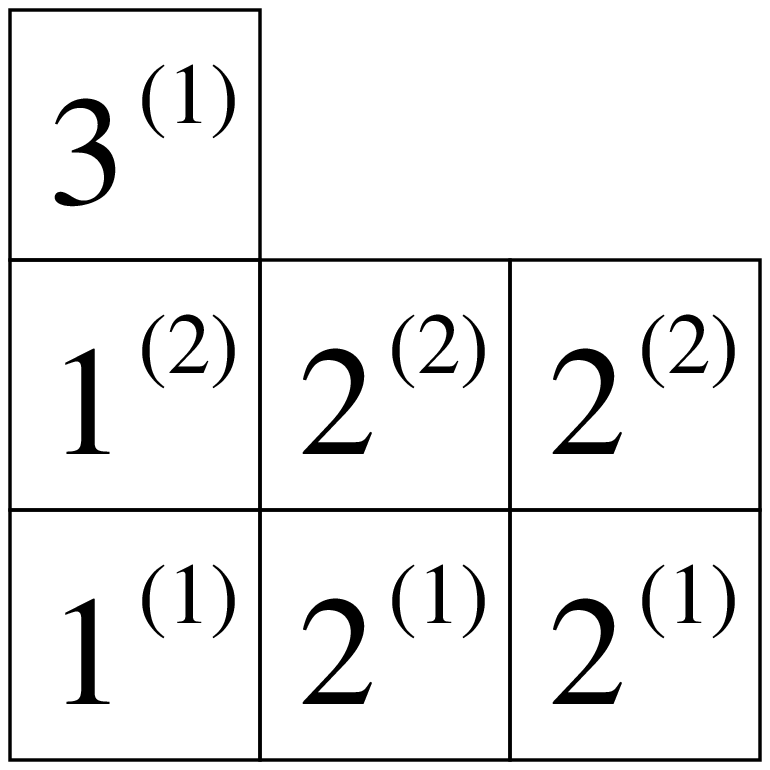}}\; ,
\end{equation*}
where the words $w^{(2)}$ and $w^{(1)}$ in the chain~\eqref{chain}
are given by
$w^{(2)}=x_3^{(1)}x_1^{(2)}x_2^{(2)}x_1^{(1)}x_2^{(1)}x_2^{(1)}x_2^{(2)}$ 
and
$w^{(1)}=x_3^{(1)}x_1^{(2)}x_2^{(2)}x_1^{(1)}x_2^{(1)}x_2^{(2)}x_2^{(1)}$.
\end{example}

For $\tab{\cdot}{\mu}$ the initial cyclage defines a partial order ranked by 
the cocharge, that is, $\co(T):=\rank(T)$. In particular, $\co(T)=\co(T')+1$
if $T'=\C(T)$ and the minimal element in this poset is
$T_{\min}=\left[x_1^{\mu_1}\cdots x_n^{\mu_n}\right]$.
For $\ttab{\cdot}{\mu}$ we would like to mimic this structure, that is, we 
wish to turn $\ttab{\cdot}{\mu}$ into a ranked poset with minimal element 
defined as the LR tableau 
$T_{\min}=\left[x_1^{\mu_1}x_2^{\mu_2}\cdots x_L^{\mu_L}\right]$,
where $x_i^{\mu_i}$ abbreviates $(x_i^{(a_i)}x_i^{(a_i-1)}\cdots x_i^{(1)})^j$ 
for $\mu_i=(j^{a_i})\in\R$.
Note that for $a\neq (1,\dots,1)$, $T_{\min}$ is no longer a one-row tableau 
but a tableau of shape $\mu_1+\mu_2+\cdots+\mu_L$ where the $j$th row is 
filled with the letters $x_i^{(j)}$ (all possible $i$) only.

Having fixed $T_{\min}$ we observe two important differences between the sets 
$\ttab{\cdot}{\mu}$ (for $a\neq (1,\dots,1)$) and $\tab{\cdot}{\mu}$.
\begin{remark}\label{rem_diff}\mbox{}
\begin{enumerate}
\item \label{rem_min}
If for $T,T_{\min}\in\tab{\cdot}{\mu}$, $w_T$ and $w_{T_{\min}}$ both start 
with the same letter then $T=T_{\min}$. Generally this is not true 
for $T,T_{\min}\in\ttab{\cdot}{\mu}$.
Indeed, the LR tableau~\eqref{gyt} of example~\ref{ex_omega} starts 
with $x_1^{(3)}$, but is not the minimal LR tableaux, which reads
\begin{equation*}
T_{\min}=\;\raisebox{-0.4cm}{\epsfxsize=2.5 cm \epsffile{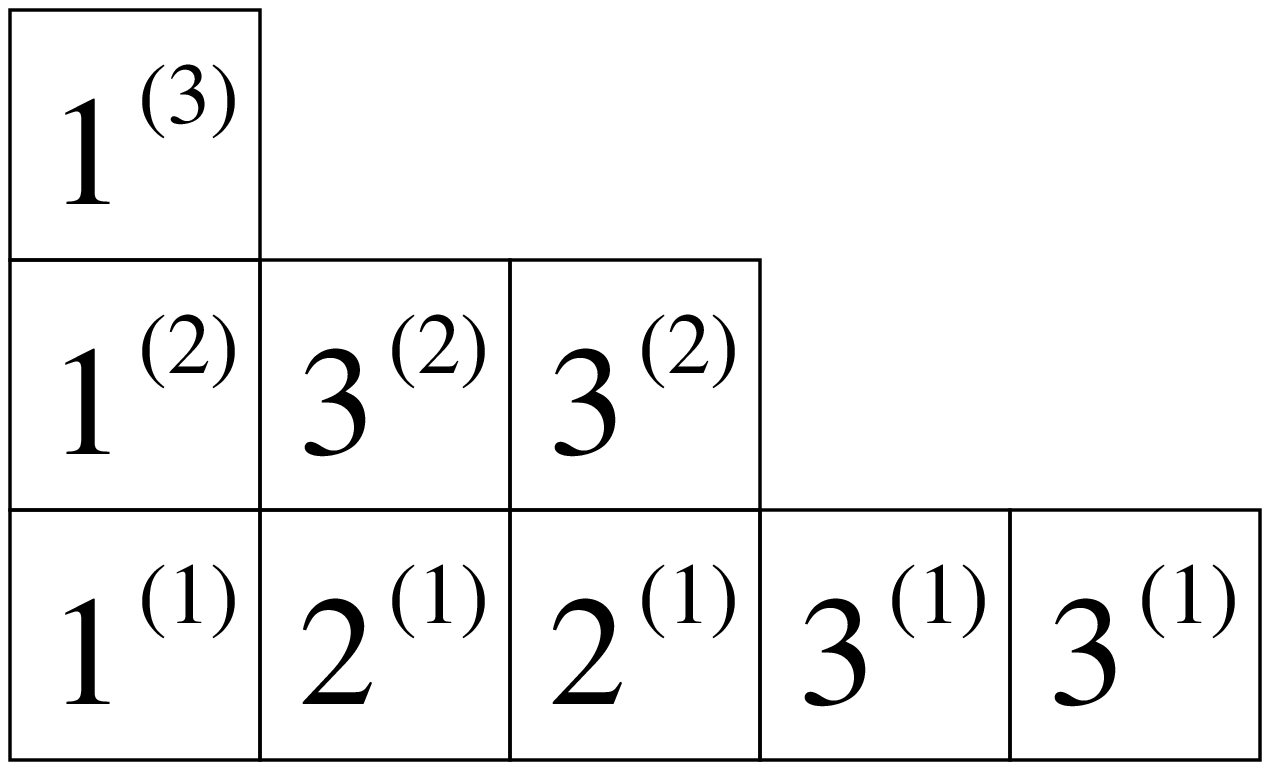}}\; .
\end{equation*}
\item
For $a=(1,\dots,1)$ the initial cyclage $\C$ has no fixed points.
In other words there is no $T\in\tab{\cdot}{\mu}$ with $T\neq T_{\min}$
such that $\C(T)=T$. For $\ttab{\cdot}{\mu}$, however, fixed points may
occur. The following LR tableau, for example, is not minimal but obeys 
$\C(T)=T$:
\begin{equation}
\label{fixpoint}
T=\;\raisebox{-0.4cm}{\epsfxsize=1.5 cm \epsffile{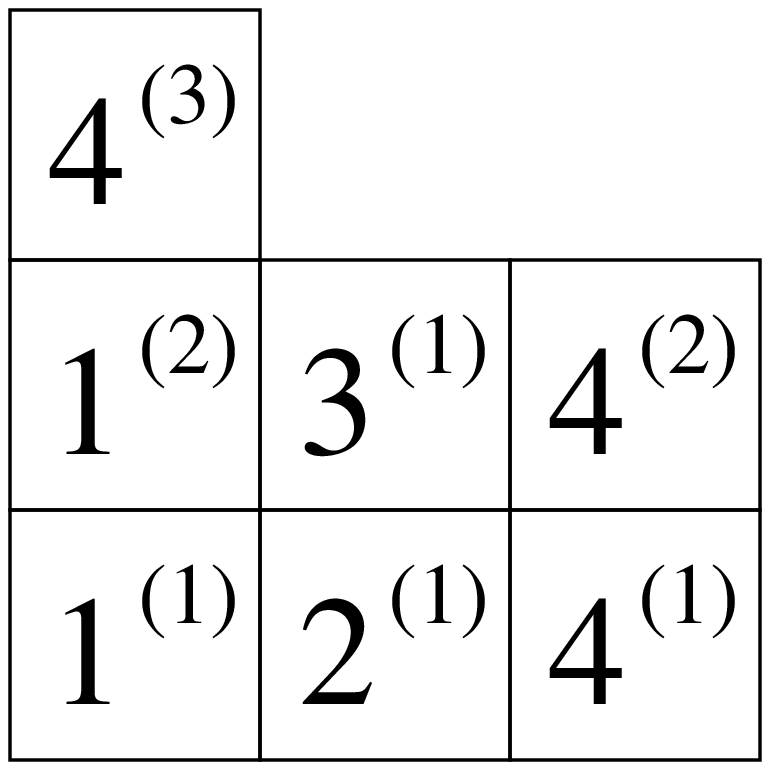}}\; .
\end{equation}
\end{enumerate}
\end{remark}

The second remark shows that the initial cyclage $\C$ does not induce 
a ranked poset structure on $\ttab{\cdot}{\mu}$ and hence needs 
further modification. For this purpose we define the ``dropping'' and 
``insertion'' operators $\D$ and $\U$, respectively. Let 
$T\in\ttab{\cdot}{\mu}$ with $\mu\in\R^L$. 
If, for some fixed $i\in\{1,\dots,L\}$ and all $j=1,\dots,\height(T)$ 
the $j$th row of $T$ contains $x_i^{(j)}$ (there may be more than one
$x_i^{(j)}$ in row $j$), then drop all the boxes containing the letters 
$x_i^{(j)}$. Repeat this operation on the reduced tableau until no more 
letters can be dropped. The final tableau defines $\D(T)$. Obviously the 
condition for dropping occurs if and only if $\height(T)=a_i$ for some $i$.
The operator $\U$ is somewhat intricate in that we only define 
$\U\circ\O\circ\D$, where $\O$ can be any content preserving operator acting 
on LR tableaux. So assume $T'=(\O\circ\D)(T)$. Then $\U$ acts on $T'$ by
reinserting all boxes that have been dropped by $\D$, inserting a box with 
filling $x_i^{(j)}$ in the $j$th row such that the conditions for a Young
tableau are satisfied (i.e., each row remains non-decreasing and each column 
strictly increasing on $X^a$). The insertion operator $\U$ is not to be 
confused with the insertion of boxes defined by the Schensted algorithm. 
$\U$ never bumps any boxes.
\begin{remark}\label{rem_U}
Note that $\U\circ\D=\Id$, but $\D\circ\U\circ\O\circ\D=\D\circ\O\circ\D$
and not $\O\circ\D$.
\end{remark}
\begin{definition}[Modified initial cyclage]
\label{def_mic}~\newline
The modified initial cyclage $\Cb:\ttab{\cdot}{\mu}\to\ttab{\cdot}{\mu}$ is 
defined as
\begin{equation}\label{mic}
\Cb=\U\circ\C\circ\D.
\end{equation}
\end{definition}
Note that $\D(T)=T$ when $\height(T)>\max\{a_1,\dots,a_L\}$, in which case
$\Cb=\C$.

Finally we are in the position to define the cocharge and charge of
an LR tableau.
\begin{definition}[Cocharge and charge]
\label{def_charge}
Let $T\in\ttab{\cdot}{\mu}$.
\begin{enum}
\item
The cocharge $\co(T)$ of $T$ is the number of times one has to 
apply $\Cb$ to obtain the minimal LR tableau $T_{\min}$.
\item
The charge is $c(T)=\|\mu\|-\co(T)$ where $\|\mu\|=\sum_{i<j}|\mu_i\cap\mu_j|$.
\end{enum}
\end{definition}

\begin{example}
For $T$ in~\eqref{fixpoint} we have
\begin{equation*}
\D(T)=\; \raisebox{-0.5cm}{\epsfxsize=0.5 cm \epsffile{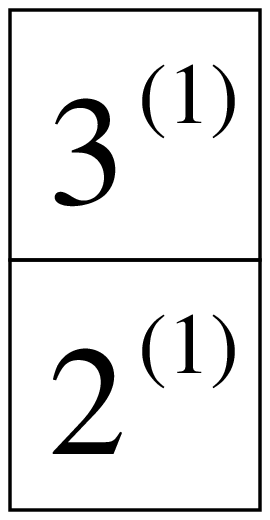}}\; ;\qquad
(\C\circ\D)(T)=\; 
 \raisebox{-0.2cm}{\epsfxsize=1 cm \epsffile{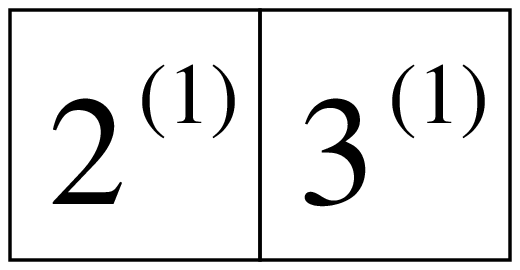}}\; ;\qquad
\Cb(T)=(\U\circ\C\circ\D)(T)=
 \raisebox{-0.5cm}{\epsfxsize=2 cm \epsffile{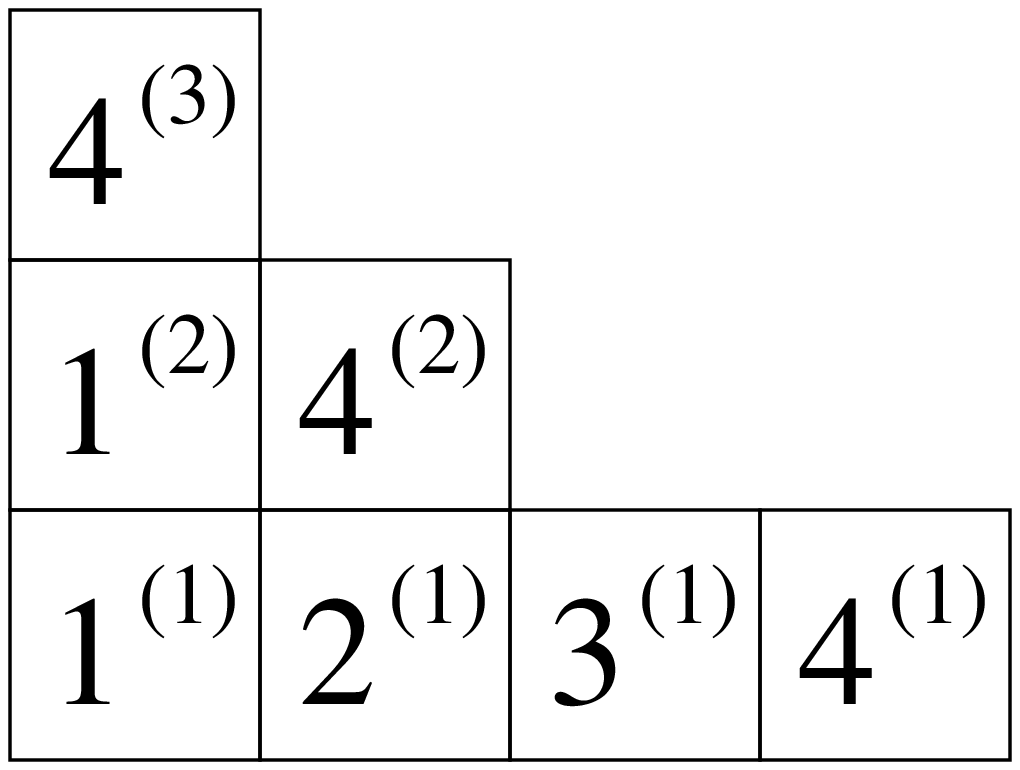}}\; .
\end{equation*}
Since $\Cb(T)=T_{\min}$ and $\|\mu\|=7$ we see that $\co(T)=1$ and $c(T)=6$.
\end{example}

Further examples of the action of the modified initial cyclage can be found
in appendix~\ref{sec_C}.

As will be shown in section~\ref{sec_poset}, the modified initial cyclage
$\Cb$ indeed turns $\ttab{\cdot}{\mu}$ into a ranked poset with
$\co:=\rank$. This implies in particular that $\co(T)$ is a bounded 
non-negative integer. In fact, we will do more and define more general 
$\la$-cyclages which induce a ranked poset structure on $\ttab{\cdot}{\mu}$.
We also show in section~\ref{sec_poset} that the charge is non-negative and 
that $\|\mu\|=\co(T_{\max})$ where 
$T_{\max}=\left[x_L^{\mu_L}\cdots x_1^{\mu_1}\right]$ has maximal cocharge.

For convenience we write $\co(w)$ instead of $\co([w])$ for $w\in\W_{\mu}$.

\section{Charge statistic representation for the generalized Kostka
polynomials}
\label{sec_LS}

There is a relation between the weight-function of definition~\ref{def_weight}
and the cocharge of definition~\ref{def_charge} which is stated in 
theorem~\ref{theo_Hco}. This relation enables us to derive an 
expression for the generalized Kostka polynomials stemming from 
the Lascoux and Sch\"utzenberger representation~\eqref{LSKostka} given in
corollary~\ref{cor_SK}.
Section~\ref{sec_proof} is devoted to the proof of theorem~\ref{theo_Hco}. 

\subsection{Relation between cocharge and weight}

To state the precise relation between cocharge and weight,
we need to introduce the anti-automorphism $\Om$ on words in $\W$. 
Recall that for every alphabet $X=\{x_1<x_2<\cdots<x_L\}$
there exists the dual alphabet $X^*=\{x_L^*<x_{L-1}^*<\cdots<x_1^*\}$.
Setting $(x_i^*)^*=x_i$, $(X^*)^*=X$. The letter $x_i^*$ which is often 
identified with $x_{L+1-i}$ is called dual to $x_i$. Under $\Om$ a word
$w=x_{i_1}x_{i_2}\ldots x_{i_k}$ in the monoid $\X$ is mapped to 
$w^*=x_{i_k}^* x_{i_{k-1}}^* \ldots x_{i_1}^*$ in $\X^*$.
Obviously, $\Om$ is an involution. For the alphabet $X^a$ of
equation~\eqref{alph_X} one may identify $(x_i^{(j)})^*$
with $x_{L+1-i}^{(a_i+1-j)}$ so that $(X^a)^*$ becomes $X^{a^*}$ with
$a^*=(a_{L+1},\ldots,a_1)$. One can easily show that if $w\in\W$ over 
$X^a$, then $\Om(w)\in\W^*$ with $\W^*=\{w\in\X^{a^*} | 
\text{ each $i$-subword of $w$ is a balanced Yamanouchi word} \}$. 
Since $\Om$ respects the Knuth equivalence relations, $\Om$ is also 
well-defined on LR tableaux by setting $\Om(T)=\Om(w_T)$.
On paths $P=p_L\otimes\cdots\otimes p_1\in\P_{\lm}$ we define 
$\Omp(P)=\Om(p_1)\otimes\cdots\otimes\Om(p_L)$.

Recall the map $w:\P_{\lm}\to\W_{\mu}$ defined in~\eqref{omega} from
paths to words. With this we can now state the following theorem.

\begin{theorem}[Weight-cocharge relation]
\label{theo_Hco}
For $n\geq 2$ and $L\geq 0$ integers let $\la\in\Int^n$ and $\mu\in\R^L$.
Then for $P\in\P_{\lm}$ the weight $H(P)$ is a non-negative integer and
\begin{equation}
\label{eqn_Hco}
H(P)=\co(\Om\circ\omega(P)).
\end{equation}
\end{theorem}

In the special case when $n=2$ and $\mu=(1^{|\la|})$ a similar
relation was noticed in~\cite{DF96}. Theorem~\ref{theo_Hco} generalizes 
theorem 5.1 of ref.~\cite{LLT95} valid when $\mu$ is a partition.
It also implies that the generalized Kostka polynomials can be 
expressed as the generating function of LR tableaux
with the charge statistic. This is summarized in the following corollary.

\begin{corollary}
\label{cor_SK}
The generalized Kostka polynomial $K_{\lm}(q)$ can be expressed as
\begin{equation}
\label{K_charge}
K_{\lm}(q)=\sum_{T\in\ttab{\la}{\mu}}q^{c(T)}.
\end{equation}
\end{corollary}

\begin{proof}
We start by rewriting the generalized cocharge Kostka polynomial 
$\Sb_{\lm}(q)$.
Recalling that $\Pb_{\lm}$ is isomorphic to $\ttab{\la}{\mu}$ it follows
from~\eqref{Sb} and~\eqref{eqn_Hco} that
\begin{equation}\label{Sb_int}
\Sb_{\lm}(q)=\sum_{T\in\ttab{\la}{\mu}}q^{\co(\Om(T))}.
\end{equation}
Since $\Om$ does not change the shape of a tableau (see for example 
ref.~\cite{F97}), but changes its content $\mu=(\mu_1,\ldots,\mu_L)$
to $(\mu_L,\ldots,\mu_1)$, and since $\Sb_{\lm}(q)=\Sb_{\la\tilde{\mu}}(q)$ 
for a permutation $\tilde{\mu}$ of $\mu$, we can drop $\Om$ in~\eqref{Sb_int}.
Recalling equation~\eqref{gK} and $c(T)=\|\mu\|-\co(T)$ completes the proof 
of~\eqref{K_charge}.
\end{proof}

\subsection{Proof of theorem~\ref{theo_Hco}}
\label{sec_proof}

The proof of theorem~\ref{theo_Hco} requires several steps. First we use the 
fact that all paths with Knuth equivalent words have the same weight. 

\begin{lemma}
\label{lem_weight}
Let $P,P'\in\P_{\cdot\mu}$ such that $\omega(P)\equiv\omega(P')$. Then
$H(P)=H(P')$.
\end{lemma}

\begin{proof}
Since $\omega(P)\equiv \omega(P')$ also $\omega(p_{i+1}\otimes p_i)\equiv
\omega(p'_{i+1}\otimes p'_i)$ (see lemma 3 on page 33 of ref.~\cite{F97}). 
In particular, they have the same shape and hence, by~\ref{i_1} of 
lemma~\ref{lem_om}, $p_{i+1}\cdot p_i$ and $p'_{i+1}\cdot p'_i$ have the 
same shape. This implies $h(p_{i+1}\otimes p_i)=h(p'_{i+1}\otimes p'_i)$ for 
all $i$ and therefore $h(P)=h(P')$. Since $\sigma_i$ only changes steps $p_i$ 
and $p_{i+1}$ and since, by~\ref{i_2} of lemma~\ref{lem_om}, $\omega(P)$ and 
$\omega(\sigma_i(P))$ have the same shape,
it follows that $\omega(\sigma_i(P))\equiv \omega(\sigma_i(P'))$. Hence,
repeating the argument, $h(\sigma_i(P))=h(\sigma_i(P'))$. This implies 
$H(P)=H(P')$.
\end{proof}

Thanks to the above lemma it suffices to prove~\eqref{eqn_Hco} for just one 
representative path $P$ for each $T\in\ttab{\cdot}{\mu}$ such that 
$[\omega(P)]=T$. Let us now find a suitable set of such paths.

Define $\P_{\mu}:=\P_{\lm}$ where $\la=(1^{|\mu|})$ and $|\mu|=\sum_i |\mu_i|$
so that for $P\in\P_{\mu}$ each letter $1,\ldots,|\mu|$ occurs exactly once. 
There is a bijection between $\P_{\mu}$ and $\W_{\mu}$. The map from 
$\P_{\mu}$ to $\W_{\mu}$ is just given by $\omega$ of equation~\eqref{omega}. 
The inverse map 
\begin{equation}
\label{omega_inv}
\omega^{-1}:\W_{\mu}\to \P_{\mu}
\end{equation}
is given as follows. Let $w=w_{|\mu|} w_{|\mu|-1}\cdots w_1$ be a word in 
$\W_{\mu}$.
Then reading $w$ from left to right, place $i$ in the rightmost empty box 
in row $k$ of step $j$ if $w_i=x_j^{(k)}$. Obviously, since $w\in\W_{\mu}$ 
the steps of the resulting path $P$ are Young tableaux of rectangular 
shapes $\mu_i$.

Denote the set of paths $P\in\P_{\mu}$ with $\omega(P)$ in 
row-representation by $\Pb_{\mu}$.
Since $\ttab{\cdot}{\mu}$ and $\W_{\mu}/\equiv$ are isomorphic, the
bijection between $\P_{\mu}$ and $\W_{\mu}$ also implies a bijection
between $\Pb_{\mu}$ and $\ttab{\cdot}{\mu}$ still denoted by $\omega$.
By lemma~\ref{lem_weight} we are thus left to prove~\eqref{eqn_Hco} for 
all $P\in\Pb_{\mu}$.

The bijection between $\ttab{\cdot}{\mu}$ and $\Pb_{\mu}$ induces a 
modified initial cyclage $\Cbp:\Pb_{\mu}\to\Pb_{\mu}$ defined as
$\Cbp:=\omega^{-1}\circ\Cb\circ\omega$. Before setting out for the proof 
of~\eqref{eqn_Hco} let us study some of the properties of this induced 
function.

\begin{figure}
\centering
\parbox{6cm}{\centering
$\raisebox{-0.4cm}{\epsfysize=1 cm \epsffile{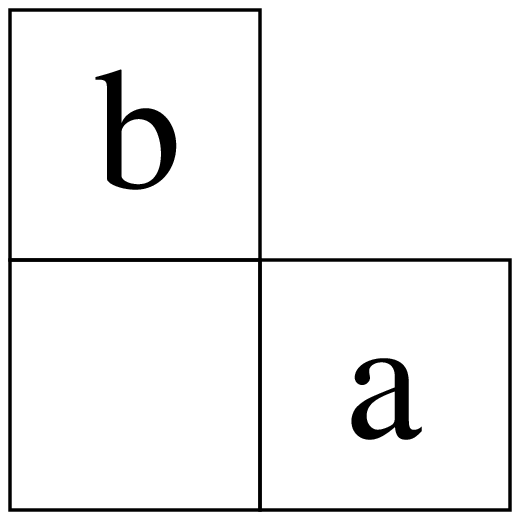}}\; \to \;
\begin{cases}
\; \epsfysize=1 cm \epsffile{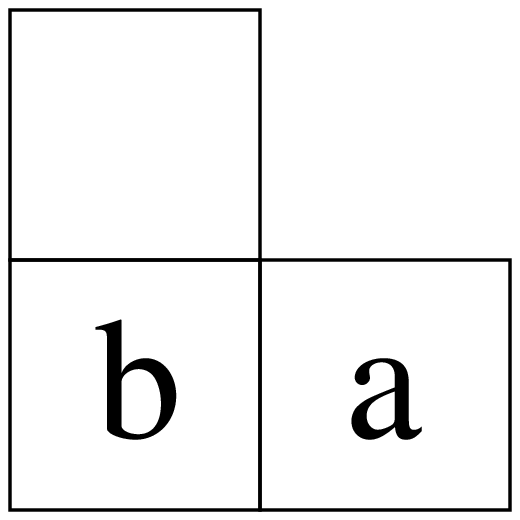} 
& \raisebox{0.3cm}{\text{if $b\leq a$}}\\[2mm]
\; \epsfysize=1 cm \epsffile{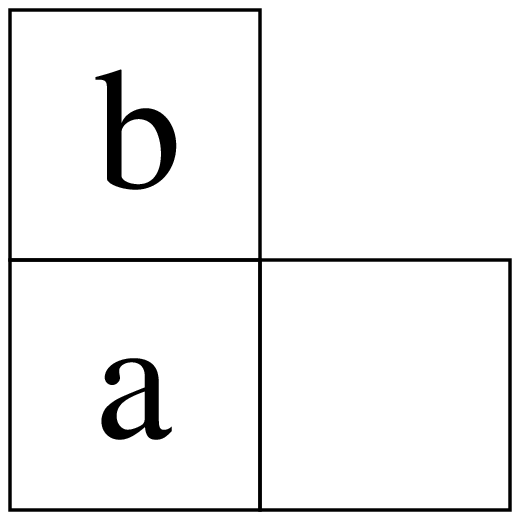}
& \raisebox{0.3cm}{\text{if $b>a$}} \end{cases}$
\caption{Sliding mechanism}\label{fig_slide}
}
\hspace{1.5cm}
\parbox{6cm}{\centering
$\raisebox{-0.4cm}{\epsfysize=1 cm \epsffile{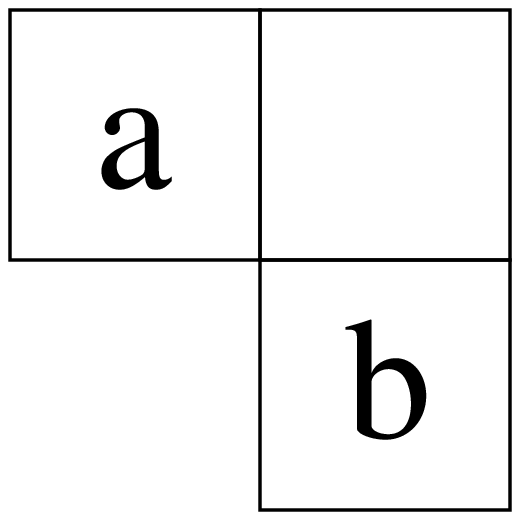}} \; \to \;
\begin{cases}
\; \epsfysize=1 cm \epsffile{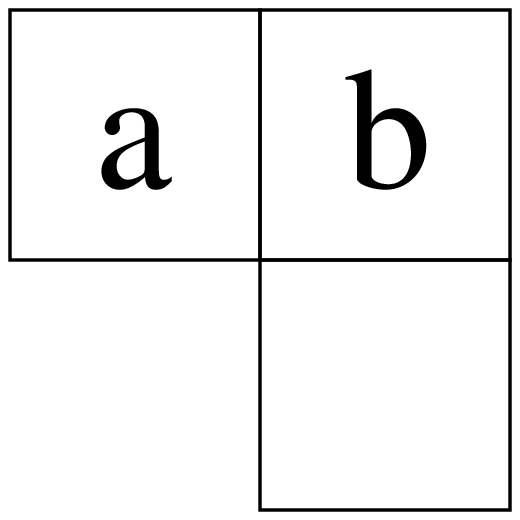} 
& \raisebox{0.3cm}{\text{if $a\leq b$}}\\[2mm]
\; \epsfysize=1 cm \epsffile{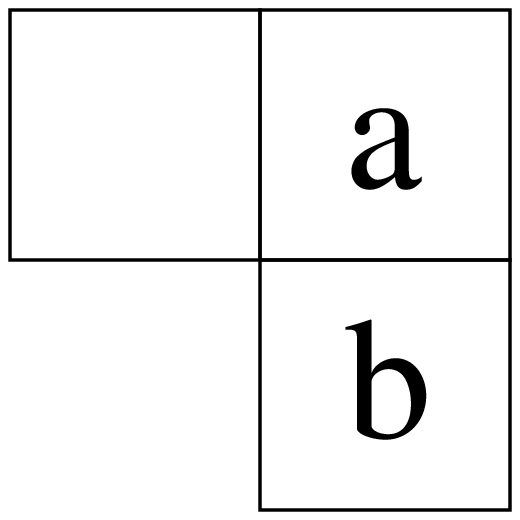}
& \raisebox{0.3cm}{\text{if $a>b$}} \end{cases}$
\caption{Inverse sliding mechanism}\label{fig_Islide}
}
\end{figure}

First consider the induced map 
$\Cp:=\S\circ\omega^{-1}\circ\C\circ\omega$ of the initial cyclage $\C$ of
definition~\ref{def_cyc}. Here $\S$ is a shift operator which decreases
the letters in each step of $P\in\P_{\mu}$ by one and hence makes $\Cp(P)$ 
a path over $\{0,1,\ldots,|\mu|-1\}$. The reason for including $\S$ in the 
definition of $\Cp$ is merely for convenience so that $\Cp$ acts only 
on one step of paths in $\Pb_{\mu}$ as will be shown in lemma~\ref{lem_Cp}.
One may always undo the effect of $\S$ by acting with $\S^{-1}$ which
adds one to each entry of a path.
To state the precise action of $\Cp$ on a path, let us briefly review 
Sch{\"u}tzenberger's (inverse) sliding mechanism~\cite{S63}.
Suppose there is an empty box with neighbours to the right and above. 
Then slide the smaller of the two neighbours in the hole; if both neighbours 
are equal choose the one above. Similarly for the inverse sliding mechanism 
consider an empty box with neighbours to the left and below. Slide the bigger 
of the two neighbours in the hole; if both are equal choose the one below.
If there is only one neighbour in either case slide this one into the
empty box. The sliding and inverse sliding mechanisms are illustrated in 
figures~\ref{fig_slide} and~\ref{fig_Islide}, respectively.

\begin{lemma} \label{lem_Cp}
Let $P=p_L\otimes\cdots\otimes p_1\in\Pb_{\mu}$ be a path over 
$\{1,2,\ldots,|\mu|\}$ and let the letter $|\mu|$ be contained in step 
$p_i$. Then
\begin{enum}
\item \label{l_1} $\Cp$ acts only on step $p_i$ of $P$, i.e.,
$\Cp(P)=p_L\otimes\cdots\otimes\Cp(p_i)\otimes\cdots\otimes p_1$, and
\item \label{l_2} $\Cp(p_i)$ is obtained by first removing $|\mu|$ from the 
top-right box of $p_i$, then using the inverse sliding mechanism to move 
the empty box to the bottom-left corner and finally inserting $0$ into the 
empty box.
\end{enum}
\end{lemma}

\begin{proof}
Since $P\in\Pb_{\mu}$ and since the largest letter $|\mu|$ occurs in step $i$,
the word $\omega(P)$ is in row-representation and of the form 
$\omega(P)=x_i^{(a_i)}u$. Let $T=[\omega(P)]$. 
In the chain of transformations~\eqref{chain} with $w=\omega(P)$
only the $i$-subword of $w$ gets changed and all letters
in $w^{(1)}$ not in the $i$-subword are shifted one position to the left.
Hence $\S\circ\omega^{-1}\circ\C(T)$ leaves all but the $i$th step in $P$
invariant which implies $(i)$.
To prove $(ii)$ observe that in row $j+1$ the empty box moves to the
left up to the point where the left neighbour is smaller than the neighbour
below. Under the map $\omega$ these two neighbours correspond to the 
two non-inverted letters $x_i^{(j)}$ and $x_i^{(j+1)}$ in $w^{(j+1)}$ 
of~\eqref{chain} used for the definition of $\C$.
\end{proof}

Some properties of the initial cyclage $\Cp$, the map $\omega$, the involution 
$\Omp$ and the isomorphism $\sigma_i$ are summarized in the following lemma.
For $P\in\P_{\lm}$ we set $h_i(P)=h(p_{i+1}\otimes p_i)$.
\begin{lemma} \label{lem_prop}
For $\la\in\Int^n$ and $\mu\in\R^L$ we have on $\P_{\lm}$
\begin{align}
\label{hO} &h_{L-i}=h_i\circ\Omp,\\
\label{sO} &\Omp\circ\sigma_i=\sigma_{L-i}\circ\Omp,\\
\intertext{and on $\Pb_{\mu}$}
\label{Omp} &\Omp=\omega^{-1}\circ\Om\circ\omega,\\
\label{sC} &[\sigma_i,\Cp]=0,\\
\label{CO} &\Cp=\Omp\circ\Cp^{-1}\circ\Omp,
\end{align}
where $\Cp^{-1}$ is defined as follows. It acts on the step with the smallest 
entry in $P\in\Pb_{\mu}$ by removing the 1, moving the empty box by the 
sliding mechanism to the top right corner and inserting $|\mu|+1$.
\end{lemma}

\begin{proof}
Let $P\in\P_{\lm}$. The energy $h_i(P)$ is 
determined by the shape of $p_{i+1}\cdot p_i$. Hence $h_i(\Omp(P))$
is determined by the shape of $\Om(p_{L-i})\cdot\Om(p_{L+1-i})
=\Om(p_{L+1-i}\cdot p_{L-i})$. But $\Om$ leaves the shape of a Young
tableau invariant (see for example ref.~\cite{F97}), yielding~\eqref{hO}.

Since the isomorphism $\sigma_i$ acts only locally on $p_{i+1}\otimes p_i$
and $\Omp$ reverses the order of the steps, it suffices to prove~\eqref{sO}
for a path of length two. Define $\tilde{p}_1\otimes\tilde{p}_2
=\sigma(p_2\otimes p_1)$ so that 
$\tilde{p}_1\cdot \tilde{p}_2=p_2\cdot p_1$. 
Acting on the last equation with $\Om$ yields 
$\Om(\tilde{p}_2)\cdot\Om(\tilde{p}_1)=\Om(p_1)\cdot\Om(p_2)$. Since $\Om$ 
does not change the shape of a Young tableau and because of the uniqueness of 
the decomposition into the product of two rectangular Young tableaux we 
conclude that $\sigma(\Om(p_1)\otimes\Om(p_2))
=\Om(\tilde{p}_2)\otimes\Om(\tilde{p}_1)$ which proves~\eqref{sO}.

Equation~\eqref{Omp} follows in a straightforward manner from
the definitions of $\omega$ and $\omega^{-1}$.

Let $P\in\Pb_{\mu}$ and let the letter $|\mu|$ be contained in 
step $p_j$ of $P$. By~\ref{l_1} of lemma~\ref{lem_Cp}, $\Cp$ acts only on 
step $p_j$, and $\sigma_i$ acts only on $p_{i+1}\otimes p_i$. Hence the
proof of \eqref{sC} reduces to showing that $[\sigma,Z_p]=0$ on
$\P_{(\mu_1,\mu_2)}$. Here $Z_p=\S\circ\omega^{-1}\circ Z\circ\omega$ and 
$Z:\W_{\mu}\to\W_{\mu}$ is defined as $Z(w)=w^{(1)}$ where
$w^{(1)}$ as given in \eqref{chain} (note that $w$ need not be in
row-representation).
Let $P=p_2\otimes p_1\in \P_{(\mu_1,\mu_2)}$ and set 
$w=w_1\ldots w_{|\mu|}:=\omega(P)$ and $\wt=\wt_1\ldots\wt_{|\mu|}:=
\omega(\sigma(P))$. 
The map $\omega:\P_{\mu}\to \W_{\mu}$ is a bijection. 
Since for a given shape $\la$ the set $\ttab{\la}{(\mu_1,\mu_2)}$ can have
at most one element, a word $w\in\W_{(\mu_1,\mu_2)}$ is uniquely specified by
$\shape(w_1\ldots w_k)$ for all $1\le k\le |\mu|$. 
Hence \eqref{sC} amounts to showing that, for 
$w'=w_1'\ldots w_{|\mu|}':=Z(w)$ and $\wt'=\wt_1'\ldots\wt_{|\mu|}':=Z(\wt)$,
$\shape(w_1'\ldots w_k')=\shape(\wt_1'\ldots \wt_k')$ for all
$1\le k\le |\mu|$. By construction, $\shape(w_1'\ldots w_k')
=\shape(w_2\ldots w_{k+1})$ and $\shape(\wt_1'\ldots \wt_k')
=\shape(\wt_2\ldots \wt_{k+1})$ for all $1\le k<|\mu|$ and by 
lemma \ref{lem_om} \ref{i_2} $\shape(w_2\ldots w_k)=
\shape(\wt_2\ldots \wt_k)$. Hence we are left to show that
$\shape(w')=\shape(\wt')$. This is can be done explicitly.
In particular, one may use that the shape of the product of two rectangular 
Young tableaux has the following form
\begin{equation}\label{sh}
\shape(p_2\cdot p_1)=\thickspace \raisebox{-0.4cm}{\epsfxsize=2 cm 
\epsffile{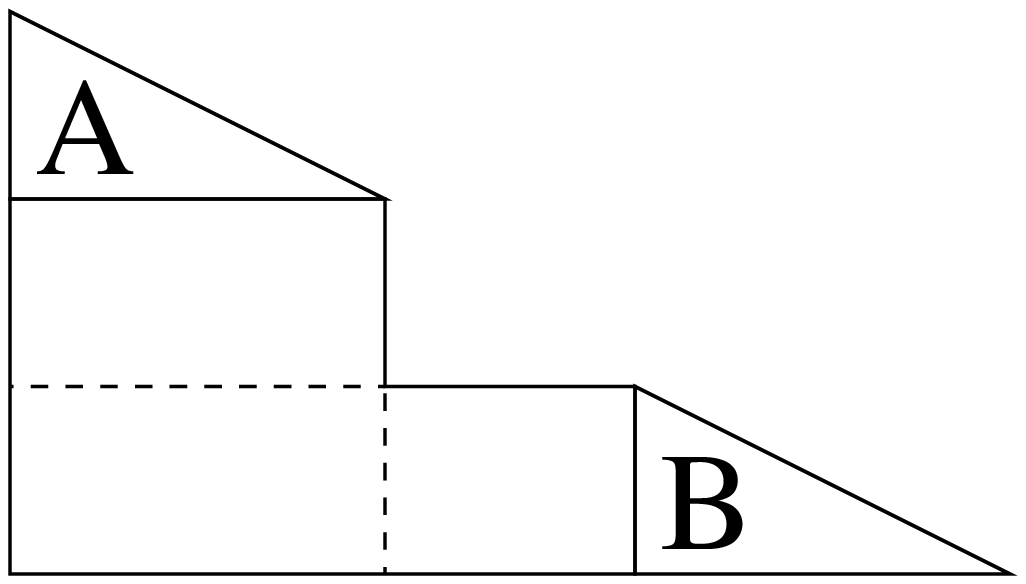}}
\end{equation}
where $A$ and $B$ are partitions and the two overlapping rectangles are the 
shapes of $p_1$ and $p_2$; one may be contained in the other.
Note that $A$ is the complement of $B$, so that knowing $A$ ($B$) fixes the 
shape. 
By lemma \ref{lem_om} \ref{i_1} also $\omega(p_2\otimes p_1)$ has 
the shape \eqref{sh}. Let $b=\shape(w_1\ldots w_{|\mu|})/
\shape(w_2\ldots w_{|\mu|})$ and
$b'=\shape(w'_1\ldots w'_{|\mu|})/\shape(w_2\ldots w_{|\mu|})$.
One may show that (i) if $b\in A$ then $b'\in B'$, (ii) if 
$b \in B$ then $b' \in A'$ and (iii) if $b\not\in A\cup B$ then
$b'\not\in A'\cup B'$. Since $b$ is the same for both $w$ and $\wt$
this implies that $\shape(w')=\shape(\wt')$. 

For the proof of~\eqref{CO} one can consider a path consisting of just
a single step thanks to~\ref{l_1} of lemma~\ref{lem_Cp}.
Suppose $p$ has $M$ boxes filled with the numbers $1,\ldots,M$.
{}From~\ref{l_2} of lemma~\ref{lem_Cp} we know that 
$\Cp$ acts by the inverse sliding mechanism and by definition
$\Cp^{-1}$ acts by the sliding mechanism.
$\Om$ acts on rectangular Young tableaux by rotation of $180^{\circ}$
and dualizing all letters. But since the inverse sliding mechanism is the
same as the sliding mechanism after rotation of $180^{\circ}$ and dualizing,
as is easily seen from figures~\ref{fig_slide} and~\ref{fig_Islide},
equation~\eqref{CO} follows.
\end{proof}

After these preliminaries we come to the heart of the proof of 
theorem~\ref{theo_Hco}. 
By lemma~\ref{lem_weight} we are left to prove equation~\eqref{eqn_Hco}
for all $P\in\Pb_{\mu}$ and by~\eqref{Omp} this is equivalent to
\begin{equation}
\label{Hco}
H'(P):=H(\Omp(P))=\co(\omega(P)).
\end{equation}
We will show that $\Cbp=\omega^{-1}\circ\Cb\circ\omega$ decreases the weight 
$H'$ of paths in $\Pb_{\mu}$ by one, i.e.,
\begin{equation}\label{HC}
H'(P)-H'(\Cbp(P))=1\quad\text{for $P\in\Pb_{\mu}$ and $P\neq P_{\min}$},
\end{equation}
where $P_{\min}:=\omega^{-1}(w_{\min})$ with 
$w_{\min}=w_{\min}(\mu)=x_1^{\mu_1}\cdots x_L^{\mu_L}$ the word corresponding 
to the minimal LR tableau $T_{\min}$.
By definition $\co(T_{\min})=0$ and one finds by direct computation that
also $H'(P_{\min})=0$. (This can be deduced from the fact that in $P_{\min}$
the number $i$ cannot be contained in a step to the left of the step 
containing $i-1$; this is also true for any $P\in\O_{P_{\min}}$ as
$P=\omega^{-1}(w_{\min}(\tilde{\mu}))$ for some permutation $\tilde{\mu}$ of 
$\mu$). The equation $H'(P_{\min})=\co(T_{\min})=0$ together with~\eqref{HC}
implies that $H'(P)$ and thus $H(P)$ are integers.
By definition $H'(P)$ is finite and non-negative. Suppose there exists
a $P\in\Pb_{\mu}$ such that $m-1<H'(P)<m$ for some integer $m$. Then we 
conclude from~\eqref{HC} that $H'(\Cbp^m(P))<0$ which contradicts the 
non-negativity of $H'$. Since $\co(T)-\co(\Cb(T))=1$ equation~\eqref{HC}
implies~\eqref{Hco} for all $P\in\Pb_{\mu}$.

Using~\eqref{h}, \eqref{H}, \eqref{hO}, \eqref{sO} and $\Omp^2=\Id$ one finds 
that
\begin{equation}
\label{Hp}
H'(P)=H(\Omp(P))=\frac{1}{|\O_P|}\sum_{P'\in\O_P}\sum_{i=1}^{L-1}
(L-i)h_i(P').
\end{equation}
Hence to show~\eqref{HC} one needs to relate the energies $h_i(P)$ and 
$h_i(\Cbp(P))$. Let us first focus on the relation between the energies of $P$ 
and $\Cp(P)$. Following ref.~\cite{LLT95} we decompose the orbit $\O_P$ of 
$P$ into chains. Let $U,V\in\O_P$ with largest entries in step $i$ and 
$i-1$, respectively. Then write $U\leadsto V$ if $\sigma_{i-1}(U)=V$ 
$(i=2,3,\ldots,L)$. Connected components of the 
resulting graph are called chains. With this notation we have the following
lemma which is proven in appendix~\ref{sec_A}.

\begin{lemma}
\label{lem_chain}
For $P\in\Pb_{\mu}$ with $\mu\in\R^L$ define the vector $\vh(P)=
(h_1(P),h_2(P),\ldots,h_{L-1}(P))$. For a chain 
$\gamma=\{P_m\leadsto P_{m-1}\leadsto \cdots \leadsto P_{\ell}\}$
such that $\sigma_{k-1}(P_k)=P_{k-1}$ and $Q_j=\Cp(P_j)$ the following 
relations hold,
\begin{align}
&\vh(Q_m)-\vh(P_m)=\ve_m\notag\\
&\vh(Q_k)-\vh(P_k)=0 \qquad \qquad\text{for $\ell<k<m$}\\
&\vh(Q_{\ell})-\vh(P_{\ell})=-\ve_{\ell-1}\notag\\
\intertext{and if $m=\ell$,}
&\vh(Q_m)-\vh(P_m)=\ve_m-\ve_{m-1}.
\end{align}
Here $\ve_m$ $(1\leq m\leq L-1)$ are the canonical basis vectors of 
$\Integer^{L-1}$ and $\ve_0=\ve_L=0$.
\end{lemma}

Thanks to equation~\eqref{sC} $\{Q_m,Q_{m-1},\ldots,Q_{\ell}\}$ is a subset of
$\O_{\Cp(P)}$. Defining $H'_{\gamma}(P):=\frac{1}{|\gamma|}\sum_{P'\in\gamma}
\sum_{i=1}^{L-1}(L-i)h_i(P')$ for a subset $\gamma\subset\O_P$,
lemma~\ref{lem_chain} ensures that 
\begin{equation} \label{H_chain}
H'_{\gamma}(P)-H'_{\Cp(\gamma)}(\Cp(P))=1
\end{equation}
for $\gamma=\{P_m\leadsto P_{m-1}\leadsto \cdots \leadsto P_{\ell}\}$
as long as $\ell>1$. For the case treated in ref.~\cite{LLT95}, where 
$[\omega(P)]\in\tab{\cdot}{\mu}$ is an ordinary Young tableau, $\ell$ is
always bigger than one when $P\neq P_{\min}$ and hence the proof of
theorem~\ref{theo_Hco} is complete in this case. For  
$[\omega(P)]\in\ttab{\cdot}{\mu}$, however, $\ell$ can take the value one
even if $P\neq P_{\min}$, due to point~\ref{rem_min} of remark~\ref{rem_diff}.
Hence~\eqref{H_chain} breaks down for $\ell=1$, i.e., when there is a 
$P'\in\gamma$ such that the letter $|\mu|$ is contained in the first step. 
However, in this case we are saved by the following lemma. Therein, the
height of a path $P=p_L\otimes\cdots\otimes p_1$ is defined as 
$\height(P):=\max_{1\leq i\leq L}\{\height(p_i)\}$.

\begin{lemma} \label{lem_equiv}
Let $P\in\Pb_{\mu}$ over $\{1,2,\ldots,|\mu|\}$. Then there exists 
a path $P'=p_L'\otimes\cdots\otimes p_1'\in\O_P$ such that $p_1'$ contains
the letter $|\mu|$ if and only if $\height(\omega(P))=\height(P)$.
\end{lemma}

\begin{proof}
Let us first show that the existence of $P'$ implies the condition on the 
height of $\omega(P)$. Since $p_1'$ contains $|\mu|$ the word $\omega(P')$ 
starts with $x_1^{(a_1)}$. By~\ref{i_2} of lemma~\ref{lem_om}
$\omega(P')$ is in row-representation. Hence the height of $\omega(P')$ 
equals the height of $p_1'$ and the first step is also (one of) the 
highest. Again by~\ref{i_2} of lemma~\ref{lem_om} $\omega(P)$ and 
$\omega(P')$ have the same shape so that the height of $\omega(P)$ equals
the height of $P$.

To prove the reverse, consider $P'\in\O_P$ such that the first step
is highest. Employing again~\ref{i_2} of lemma~\ref{lem_om} we see that
the height of $\omega(P')$ equals the height of the first step.
Now suppose that $p'_1$ does not contain $|\mu|$. This means that 
$\omega(P')=x_i^{(a_i)}u$ for some $u\in\W$ with $i>1$. Since $P'$ is
in row representation $x_i^{(a_i)}$ must be above $x_1^{(a_1)}$ in 
$[\omega(P')]$. This contradicts the fact that the height of $\omega(P')$ 
is the height of $p_1'$.
\end{proof}

The previous lemma shows that there exist chains $\gamma$ such that
$\gamma=\{P_m\leadsto\cdots\leadsto P_1\}$ (so that~\eqref{H_chain} 
is violated) if and only if the modified initial charge $\Cbp$ differs from 
$\S^{-1}\circ\Cp$. This is the case because the dropping and insertion 
operators $\Dp:=\omega^{-1}\circ\D\circ\omega$ and 
$\Up:=\omega^{-1}\circ\U\circ\omega$ in the relation 
$\Cbp=\Up\circ\S^{-1}\circ\Cp\circ\Dp$ only act non-trivially 
when the height of $\omega(P)$ equals the height of $P$, 
or equivalently by lemma~\ref{lem_equiv}, when there exists a chain 
$\gamma=\{P_m\leadsto\cdots\leadsto P_1\}$.
The dropping operator, however, does not change the weight of a path
as shown in the following lemma.

\begin{lemma} \label{lem_drop}
For $\mu\in\R^L$ let $P\in\Pb_{\mu}$ such that 
$\height(P)=\height(\omega(P))$. Then $H'(P)=H'(\Dp(P))$.
\end{lemma}

Lemmas~\ref{lem_chain}--\ref{lem_drop} imply~\eqref{HC} and hence 
theorem~\ref{theo_Hco} for the following reason. For $P\in\Pb_{\mu}$, 
the path $\Dp(P)$ does not contain any chains 
$\gamma=\{P_m\leadsto\cdots\leadsto P_1\}$ thanks to lemma~\ref{lem_equiv}. 
Hence $H'(P)=H'(\Dp(P))=H'(\Cp\circ\Dp(P))+1$.
On the other hand $H'(\Cp\circ \Dp(P))=H'(\Cbp(P))$ since $\S^{-1}$ does
not change the energy of a path and because of remark~\ref{rem_U} 
and lemma~\ref{lem_drop}. This proves~\eqref{HC}.

\begin{proof}[Proof of lemma~\ref{lem_drop}]
For a path in $\P_{\lm}$ we refer to $\mu$ as its content.
Now suppose the path $P$ of lemma~\ref{lem_drop} has $k$ steps of 
shape $\nu\in\R$ where $k\geq 1$ and $\height(\nu)=\height(P)$.
Then all $P'\in\O_P$ have $k$ steps of shape $\nu$ and, by~\ref{i_2} of 
lemma~\ref{lem_om}, $\height(\omega(P'))=\height(\nu)$.
Define $\D_{\nu}(P')$ as the path obtained from $P'$ by dropping all steps
of shape $\nu$.
Let $\eta$ be the content of $\D_{\nu}(P)$. Then for each permutation 
$\tilde{\eta}$ of $\eta$ define the suborbit $\S_{\tilde{\eta}}\subset\O_P$ as
\begin{equation}
\S_{\tilde{\eta}}=\{\Pt\in\O_P|~\text{content$(\D_{\nu}(\Pt))=
\tilde{\eta}$}\}.
\end{equation}
Then clearly $\O_P$ is the disjoint union of $\S_{\tilde{\eta}}$ over all
permutations $\tilde{\eta}$ of $\eta$, and $|\S_{\tilde{\eta}}|=\bin{L}{k}$.

Let us now show that for any $Q\in\S_{\tilde{\eta}}$ 
\begin{equation}\label{drop}
\sum_{\Pt\in\S_{\tilde{\eta}}}\sum_{i=1}^{L-1}ih_{L-i}(\Pt)
=\bin{L}{k}\sum_{i=1}^{L-k-1}ih_{L-k-i}(\D_{\nu}(Q)).
\end{equation}
Since $\Dp$ is a composition of $\D_{\nu}$'s, equation~\eqref{drop} clearly
implies $H'(P)=H'(\Dp(P))$.

To prove~\eqref{drop} we first study some properties of the energy
$h_i(\Pt)$ for $\Pt\in\S_{\tilde{\eta}}$. For $\Pt=\pt_L\otimes\cdots
\otimes \pt_1\in\S_{\tilde{\eta}}$, let
$\tilde{\mu}=(\tilde{\mu}_1,\ldots,\tilde{\mu}_L)$ be the content of $\Pt$
and define $L\geq m_1>\cdots>m_k\geq 1$ such that $\tilde{\mu}_{m_i}=\nu$.
Then $\Pt'=\sigma_{m_i+1}\circ\sigma_{m_i}(\Pt)$ is also in 
$\S_{\tilde{\eta}}$ and, for $1\leq i\leq k$,
\begin{align}
\label{h1} &h_{m_i}(\Pt)=h_{m_i-1}(\Pt)=0,\\
\label{h2} &h_{m_i}(\Pt')=h_{m_i+1}(\Pt).
\end{align}

The proof of~\eqref{h1} and~\eqref{h2} makes extensive use of
lemma~\ref{lem_om}. For two steps $p\in\B_{\la}$ and $p'\in\B_{\la'}$
let us call the shape $\la+\la'$ minimal because $h(p\otimes p')=0$
if $\shape(p\cdot p')=\la+\la'$.
Equation~\eqref{h1} states that 
$\pt_{m_i+1}\cdot \pt_{m_i}$ and $\pt_{m_i}\cdot \pt_{m_i-1}$
have minimal shape, or equivalently by~\ref{i_1} of lemma~\ref{lem_om},
that $\omega(\pt_{m_i+1}\otimes \pt_{m_i})$ and 
$\omega(\pt_{m_i}\otimes \pt_{m_i-1})$ have shapes $\tilde{\mu}_{m_i+1}+\nu$ 
and $\nu+\tilde{\mu}_{m_i-1}$, respectively. But since the height of 
$\omega(\Pt)$ is the height of $\nu$, the heights of 
$\omega(\pt_{m_i+1}\otimes \pt_{m_i})$ and $\omega(\pt_{m_i}\otimes 
\pt_{m_i-1})$ equal the height of $\nu$, and hence their shape has to be 
minimal.

We now turn to the proof of~\eqref{h2}. Denote $\Pt'=\pt'_L\otimes\cdots
\otimes\pt'_1$. Since $\Pt'=\sigma_{m_i+1}\circ\sigma_{m_i}(\Pt)$ we know 
by~\ref{i_2} of lemma~\ref{lem_om} that 
$\omega(\pt_{m_i+2}\otimes \pt_{m_i+1}\otimes \pt_{m_i})$ and
$\omega(\pt'_{m_i+2}\otimes \pt'_{m_i+1}\otimes \pt'_{m_i})$ 
have the same shape. But since by $\omega(\pt_{m_i+1}\otimes \pt_{m_i})$ and 
$\omega(\pt'_{m_i+2}\otimes \pt'_{m_i+1})$ have minimal shape by~\eqref{h1}
we can conclude that $\omega(\pt_{m_i+2}\otimes \pt_{m_i+1})$ and 
$\omega(\pt'_{m_i+1}\otimes \pt'_{m_i})$ have the same shape.
Hence by~\ref{i_1} of lemma~\ref{lem_om} also $\pt_{m_i+2}\cdot \pt_{m_i+1}$
and $\pt'_{m_i+1}\cdot \pt'_{m_i}$ have the same shape which 
implies~\eqref{h2}.

Analogous to the proof of~\eqref{h2} we find that for $\Pt'=\sigma_{m_i}(\Pt)$
the tableaux $\pt_{m_i+1}\cdot \pt_{m_i-1}$ and $\pt'_{m_i}\cdot \pt'_{m_i-1}$
have the same shape. Setting $\Pt$ to $Q$ in this argument shows that 
$h_{L-k-i}(\D_{\nu}(Q))$ is independent of $Q\in\S_{\tilde{\eta}}$.
Hence we can restrict our attention to $Q\in\S_{\tilde{\eta}}$
with steps 1 to $k$ of shape $\nu$ in the following.

If $k=L$ or $L-1$ the right-hand side of~\eqref{drop} is zero due to the
empty sum. Equation~\eqref{h1} ensures that the left-hand side is zero
as well. If $1\leq k\leq L-2$ set $X_i:=h_{L-k-i}(\D_{\nu}(Q))$ for 
$1\leq i< L-k$. Define $r_j$ as $r_j=L+1-m_j-j$ for $1\leq j\leq k$
and $r_0=0$, $r_{k+1}=L-k$ for a given $\Pt$ where, as above, the $m_i$
are the positions of the steps of shape $\nu$. Treating $X_i$ as an 
indeterminate we see from~\eqref{h1} and~\eqref{h2} that the contribution 
to $X_i$ from $\sum_{i=1}^{L-1}ih_{L-i}(\Pt)$ is given by
\begin{alignat*}{2}
&(i+j)& \quad & \text{for $r_j<i<r_{j+1}$ and $0\leq j\leq k$},\\
&0& \quad & \text{for $i=r_j$ and $1\leq j\leq k$}.
\end{alignat*}
Summing over all $\Pt\in\S_{\tilde{\eta}}$ or, equivalently, over all possible 
$r_i$ we find that
\begin{equation*}
\begin{split}
\sum_{\tilde{P}\in\S_{\tilde{\mu}}}\sum_{i=1}^{L-1}ih_{L-i}(\tilde{P})
=&\sum_{i=1}^{L-k-1}X_i
\sum_{j=0}^k (i+j)\sum_{r_0\leq\cdots\leq r_j<i<r_{j+1}\leq 
\cdots\leq r_{k+1}}1\\
=&\sum_{i=1}^{L-k-1}X_i
\sum_{j=0}^k (i+j)\bin{i+j-1}{j}\bin{L-i-j-1}{k-j}\\
=&\bin{L}{k}\sum_{i=1}^{L-k-1}iX_i
\end{split}
\end{equation*}
where the last step follows from (a special case of) the $_2F_1$ Gau{\ss} sum.
Recalling that $X_i=h_{L-k-i}(\D_{\nu}(Q))$ this proves equation~\eqref{drop} 
and hence lemma~\ref{lem_drop}.
\end{proof}

\section{The poset structure on $\ttab{\cdot}{\mu}$}
\label{sec_poset}

As shown in theorem~\ref{theo_Hco}, the weight and the cocharge are related 
as $H(P)=\co(\Om\circ\omega(P))$. Since $H(P)$ is by its 
definition~\ref{def_weight} finite and for each LR tableau 
$T\in\ttab{\cdot}{\mu}$ there exists a path $P\in\P_{\cdot\mu}$ such that 
$T=[\omega(P)]$, theorem~\ref{theo_Hco} immediately implies that $\co(T)$ is 
finite for all $T\in\ttab{\cdot}{\mu}$.
This means that each $T\in\ttab{\cdot}{\mu}$ reaches the minimal LR tableau
$T_{\min}$ after a finite number of applications of $\Cb$.
This, in turn, ensures the following corollary.
\begin{corollary}
\label{cor_finite}
The modified initial cyclage $\Cb$ induces a ranked poset structure
on $\ttab{\cdot}{\mu}$.
\end{corollary}
The statement of this corollary can be extended to more general cyclages
which generalize the $\la$-cyclages of Lascoux and 
Sch\"utzenberger~\cite{LS81,L89} for Young tableaux $T\in\tab{\cdot}{\mu}$. 
We define $\la$-cyclages for LR tableaux $T\in\ttab{\cdot}{\mu}$ in 
section~\ref{sec_lac} and prove the analogue of corollary~\ref{cor_finite} in
theorem~\ref{theo_poset}. In section~\ref{sec_prop} we deduce 
several important properties of the charge and cocharge which are 
needed in section~\ref{sec_rec} to prove recurrences for the 
A$_{n-1}$ supernomials and generalized Kostka polynomials.

\subsection{The $\la$-cyclage and $\la$-cocyclage}
\label{sec_lac}

The $\la$-cyclage is a generalization of the initial cyclage $\C$.
Let us first define the cyclage operator $\Z$ on words 
$w=x_i^{(a_i)}u\in\W_{\mu}$ as $\Z(w)=w^{(1)}$ where $w^{(1)}$ is defined
as in~\eqref{chain} by dropping the restriction that $w$ is in 
row-representation. The only additional requirement is that 
all $w'$ in the orbit of $w$ (i.e., all $w'$ such that 
$\omega^{-1}(w')\in\O_{\omega^{-1}(w)}$ with $\omega^{-1}$ as defined 
in~\eqref{omega_inv}) start with a letter different from $x_1^{(a_1)}$. 
If $w$ violates this condition we set $\Z(w)=0$.

Similarly one can define a cocyclage $\Z^{-1}$ on a word $w=ux_i^{(1)}\in
\W_{\mu}$. This time $\Z^{-1}(w)=0$ if there exists a $w'$ in the orbit of 
$w$ ending with $x_1^{(1)}$. If not, set $w^{(1)}=w$ and construct the 
chain of transformations
\begin{equation}
w^{(1)}\to w^{(2)}\to \cdots \to w^{(a_i)},
\end{equation}
where $w^{(j+1)}$ is obtained from $w^{(j)}$ by exchanging the last letter 
$x_i^{(j)}$ with its inverted partner $x_i^{(j+1)}$ in the subword of
$w^{(j)}$ consisting of the letters $x_i^{(j)}$ and $x_i^{(j+1)}$ only. 
Then the cocyclage $\Z^{-1}(w)$ is obtained by cycling the last
letter $x_i^{(a_i)}$ in $w^{(a_i)}$ to the front of the word. 
Obviously, $\Z^{-1}$ is the inverse of $\Z$.

The cyclage $\Z$ and cocyclage $\Z^{-1}$ defined on the set of words have 
analogues on the set of LR tableaux.
A $\la$-cyclage $\Z_{\la}$ on an LR tableau $T\in\ttab{\cdot}{\mu}$ is defined 
as $\Z_{\la}(T)=\Z(w)$ where $w=x_i^{(a_i)}u\in\W_{\mu}$ is a word
such that $T=[w]$ and $\shape(u)=\la$. Similarly, a $\la$-cocyclage 
$\Z^{-1}_{\la}$ is defined as $\Z_{\la}^{-1}(T)=\Z^{-1}(w)$ if 
$w=ux_i^{(1)}$ is a word such that $T=[w]$ and $\shape(u)=\la$. In both cases 
if no such $w$ exists we set $\Z_{\la}(T)=0$ and $\Z^{-1}_{\la}(T)=0$, 
respectively. The $\la$-cyclage and $\la$-cocyclage obey
\begin{equation}\label{ZZ}
\begin{aligned}
\Z_{\la}\circ\Z_{\la}^{-1}(T)&=T \qquad\text{if $\Z_{\la}^{-1}(T)\neq 0$},\\
\Z_{\la}^{-1}\circ\Z_{\la}(T)&=T \qquad\text{if $\Z_{\la}(T)\neq 0$}.
\end{aligned}
\end{equation}

In analogy with~\eqref{mic} we also define the modified $\la$-cyclage as
\begin{equation}
\Zb_{\la}=\U\circ\Z_{\la}\circ\D
\end{equation}
where $\D(T)$ was defined in section~\ref{sec_in_cyc} by successively 
dropping all $x_i^{(j)}$ $(1\leq j\leq a_i)$ from $T$ if
$\height(T)=\height(\mu_i)$ (recall that $\height(\mu_i)=a_i$)
and $\U$ was defined by reinserting all $x_i^{(j)}$ dropped by $\D$ in row 
$j$ such that the Young tableau conditions are satisfied. Similarly we set
\begin{equation}
\Zb'_{\la}=\U'\circ\Z^{-1}_{\la}\circ\D',
\end{equation}
where $\D'(T)$ is the tableau obtained by dropping all $x_i^{(j)}$
$(1\leq j\leq a_i)$ if $\width(T)=\width(\mu_i)$ and $\U'$ reinserts
all boxes dropped by $\D'$ such that there is one $x_i^{(j)}$ in each
column and the resulting object is again in $\ttab{\cdot}{\mu}$. 

The initial cyclage $\C$ is a special $\la$-cyclage since for each 
$T\in\ttab{\cdot}{\mu}$ there always exists a partition
$\la$ such that $\C(T)=\Z_{\la}(T)$. Hence also $\Cb(T)=\Zb_{\la}(T)$
for this $\la$.

For ordinary Young tableaux $T\in\tab{\cdot}{\mu}$, the cyclages 
$\Z_{\la}$ and $\Z^{-1}_{\la}$ have been considered in refs.~\cite{LS81,L89}. 
It was shown in ref.~\cite{LS81} that the $\la$-cyclages induce 
a ranked poset structure on $\tab{\cdot}{\mu}$ with the cocharge of a 
Young tableau being its rank and the minimal element being 
$T_{\min}=[x_1^{\mu_1}\cdots x_L^{\mu_L}]$. Thanks to
$\Zb_{\la}(T)=\Z_{\la}(T)$ for $T\in\tab{\cdot}{\mu}$ also $\Zb_{\la}$
induces a ranked poset structure on $\tab{\cdot}{\mu}$. Note, however,
that $\Zb'_{\la}\neq \Z_{\la}^{-1}$ even on $\tab{\cdot}{\mu}$
since for example $\Z_{\la}^{-1}(T)=0$, but $\Zb'_{\la}(T)=[3211]$
for $T=[2311]$ and $\la=(1)$.

{}From $\Zb_{\la}$ and $\Zb'_{\la}$ we may now define a cyclage- and
cocyclage-graph.
\begin{definition}[Cyclage- and cocyclage-graph]
For $\mu\in\R^L$, the cyclage-graph $\T(\mu)$ is defined by connecting all 
$T,U\in\ttab{\cdot}{\mu}$ as $T\to U$ if there exists a partition $\la$
such that $U=\Zb_{\la}(T)$. Similarly, the cocyclage-graph $\T'(\mu)$
is obtained by connecting $T\to U$ if there exists a $\la$ such that
$U=\Zb'_{\la}(T)$.
\end{definition}

An example of a cyclage-graph is given in appendix~\ref{sec_C}.
The cyclage- and cocyclage-graphs are related by an involution
\begin{equation*}
\La:\ttab{\la}{\mu}\to \ttab{\la^{\top}}{\mu^{\star}}
\end{equation*}
defined as follows.
Let $T=[w_1\cdots w_{\ell}]\in\ttab{\la}{\mu}$ with $w_i\in X^a$. Then 
$\La(T)=[w'_{\ell}\cdots w'_1]$ where $w'_m=x_i^{(k)}$ if $w_m=x_i^{(j)}$
is the $k$th $x_i^{(j)}$ in $T$ from the left. One may easily check that
$\La$ respects the Knuth equivalence relations~\eqref{knuth} 
and is therefore indeed a function on LR tableaux.
The $i$th row of $T$ gets mapped to the $i$th column in 
$\La(T)$ and hence the shape of $\La(T)$ is indeed $\la^{\top}$.

For example
\begin{equation*}
T=\; \raisebox{-0.6cm}{\epsfxsize=1.5 cm \epsffile{fig_ex7.ps}}
\; \overset{\La}{\mapsto}\;
\raisebox{-0.6cm}{\epsfxsize=1.5 cm \epsffile{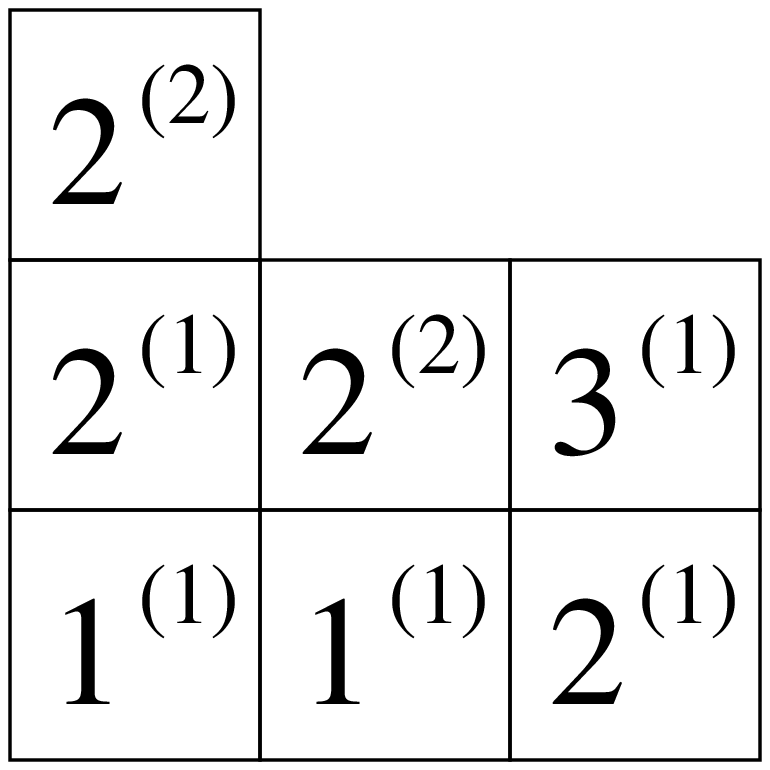}}\; .
\end{equation*}

\begin{lemma}\label{lem_TT'}
For $\mu\in\R^L$, $\T(\mu)=\La\T'(\mu^{\star})$ or, equivalently,
$\T'(\mu)=\La\T(\mu^{\star})$.
\end{lemma}

\begin{proof}
Observe that $\D'=\La\circ\D\circ\La$, $\U'=\La\circ\U\circ\La$
and $\Z_{\la}=\La\circ\Z_{\la^{\top}}^{-1}\circ\La$. This implies
\begin{equation}
\Zb_{\la}=\La\circ\Zb'_{\la^{\top}}\circ\La.
\end{equation}
Hence, for $T,T'\in\ttab{\cdot}{\mu}$ such that $T'=\Zb_{\la}(T)$ one finds
$\La(T')=\Zb'_{\la^{\top}}\circ\La(T)$ which proves the lemma.
\end{proof}

We now wish to show that both $\T(\mu)$ and $\T'(\mu)$ induce a ranked poset 
structure on the set of LR tableaux $\ttab{\cdot}{\mu}$.
To prove this we extend the standardization embedding
\begin{equation}\label{standLS}
\theta:\T(\mu)\hookrightarrow\T((1^{|\mu|})),
\end{equation}
of Lascoux and Sch\"utzenberger~\cite{LS81,L89} (see also chapter 2.6 of 
ref.~\cite{B94}) when $\mu$ is a partition to the case when 
$\mu\in\R^L$. Define the map $\phi$ on LR tableaux as follows:
\begin{center}
change the rightmost $x_1^{(j)}$ to $x_2^{(j)}$ for all $1\leq j\leq a_1$.
\end{center}
If $\height(\mu_1)=\height(\mu_2)$ and $\width(\mu_1)>\width(\mu_2)$ or
$\mu_2=0$ then $\phi(T)$ is an LR tableau of the same shape as $T$ and of
content $\mu'=(\mu_1-(1^{a_1}),\mu_2+(1^{a_1}),\mu_3,\ldots,\mu_L)$.
Denote by $\phi'$ the map $\phi$ restricted to the case when $\mu_2=0$.
One can show that $\Zb_{\la}\circ\phi'(T)=0$ if and only if
$\Zb_{\la}(T)=0$, and furthermore $[\phi',\Zb_{\la}]=0$.
(These statements can, for example, be proven by going over to paths using 
the map $\omega$ and noting that $\omega^{-1}\circ\phi'\circ\omega$
only acts on steps one and two --which is empty-- and 
$S\circ\omega^{-1}\circ\Zb_{\la}\circ\omega$ only acts on the step containing
the biggest entry in analogy to lemma~\ref{lem_Cp}. For the first statement 
it is sufficient to consider a path of length three for which it can be 
explicitly varified. Assuming $\Zb_{\la}(T)\neq 0$ the second statement
then follows trivially since the two operators act on different steps in
the path).

Denote by $\G$ the group spanned by $\omega\circ\sigma_i\circ\omega^{-1}$,
where $\sigma_i$ is the isomorphism of definition~\ref{def_iso} and
$\omega$ and $\omega^{-1}$ are defined in~\eqref{omega} 
and~\eqref{omega_inv}, respectively. Then, in analogy to~\eqref{standLS},
there exists an embedding
\begin{equation}\label{stand}
\theta:\T(\mu)\hookrightarrow\T(\nu) \qquad\text{$\nu^{\star}$ a partition.}
\end{equation}
for $\mu\in\R^L$ by combining $\phi'$ with the action of $\G$.
Since both $\tau\in\G$ and $\phi'$ are compatible with the cyclages
(the proof of the first statement is analogous to the proof of~\eqref{sC}), 
we find that
\begin{equation}\label{thZ}
[\theta,\Zb_{\la}]=0.
\end{equation}

\begin{example}
If $T$ is the LR tableau of equation~\eqref{gyt} then under $\theta$ its content
\begin{equation*}
\mu=\bigl(
\;\raisebox{-0.15cm}{\epsfxsize=0.3 cm \epsffile{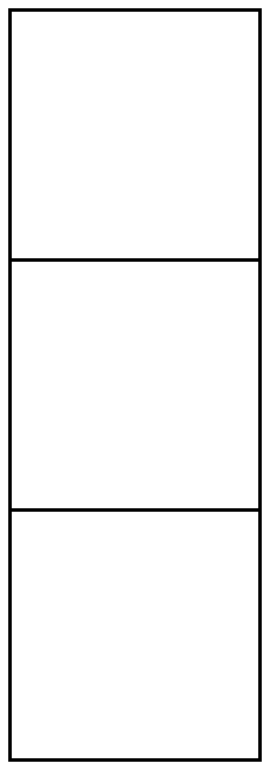}}\;,
\;\raisebox{-0.15cm}{\epsfxsize=0.6 cm \epsffile{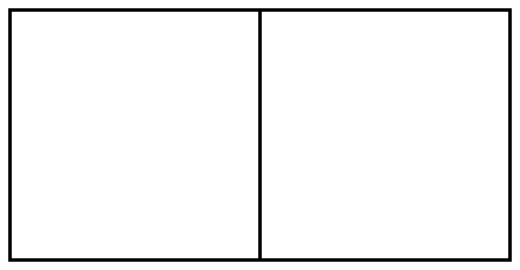}}\;,
\;\raisebox{-0.15cm}{\epsfxsize=0.6 cm \epsffile{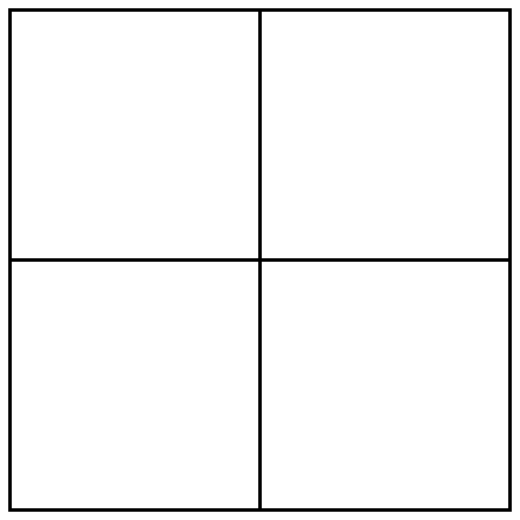}}\;
\bigr)\quad\text{will be changed to}\quad
\nu=\bigl(
\;\raisebox{-0.15cm}{\epsfxsize=0.3 cm \epsffile{fig_ex15.ps}}\;,
\;\raisebox{-0.15cm}{\epsfxsize=0.3 cm \epsffile{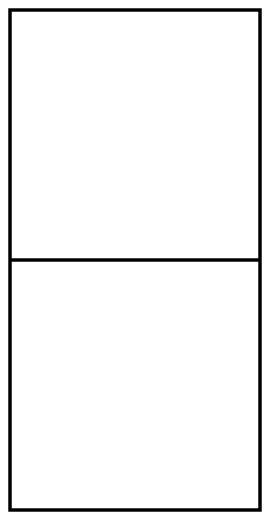}}\;,
\;\raisebox{-0.15cm}{\epsfxsize=0.3 cm \epsffile{fig_ex18.ps}}\;,
\;\raisebox{-0.15cm}{\epsfxsize=0.3 cm \epsffile{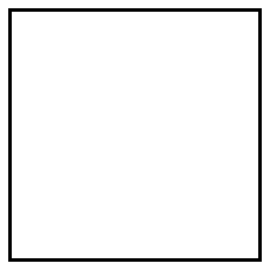}}\;,
\;\raisebox{-0.15cm}{\epsfxsize=0.3 cm \epsffile{fig_ex19.ps}}\;
\bigr).
\end{equation*}
The standardization $\theta(T)$ can be determined from
\begin{equation*}
T\xrightarrow{\tau_1}
\;\raisebox{-0.5cm}{\epsfxsize=1.8 cm \epsffile{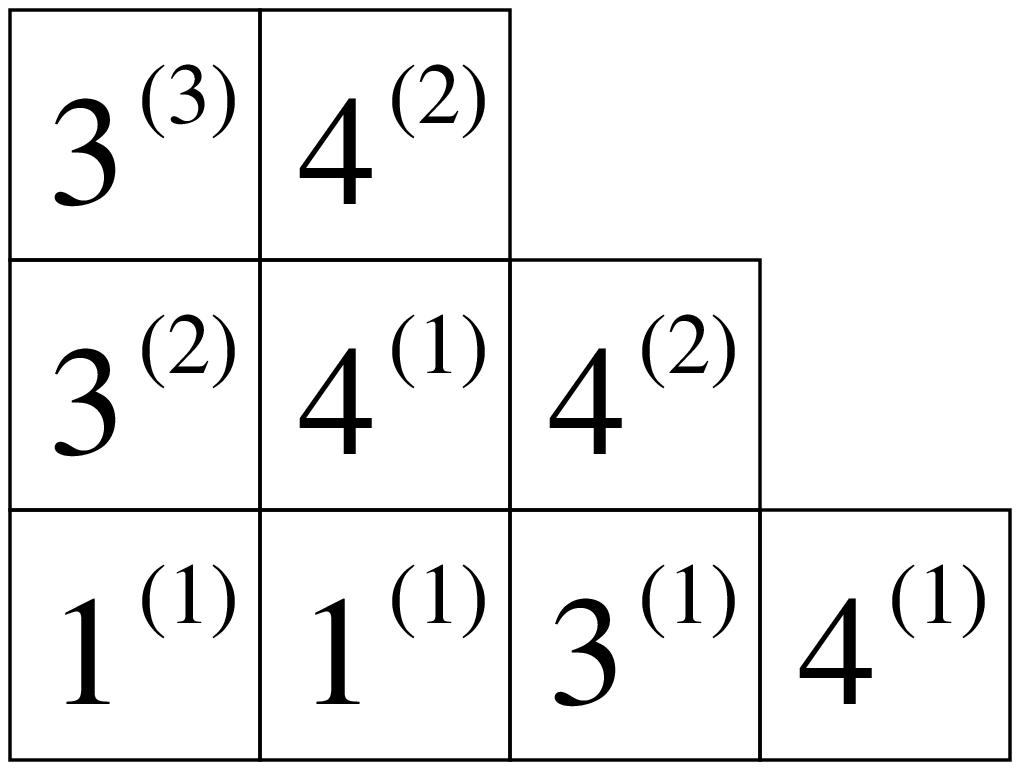}}\;
\overset{\phi}{\hookrightarrow}
\;\raisebox{-0.5cm}{\epsfxsize=1.8 cm \epsffile{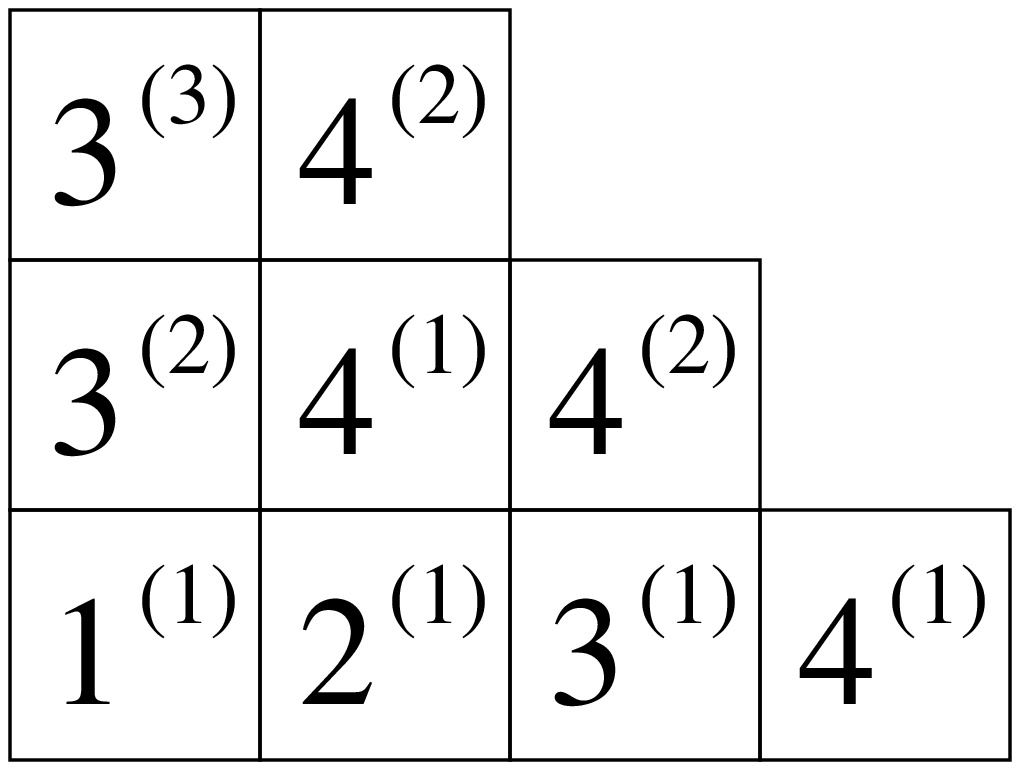}}\;
\xrightarrow{\tau_2}
\;\raisebox{-0.5cm}{\epsfxsize=1.8 cm \epsffile{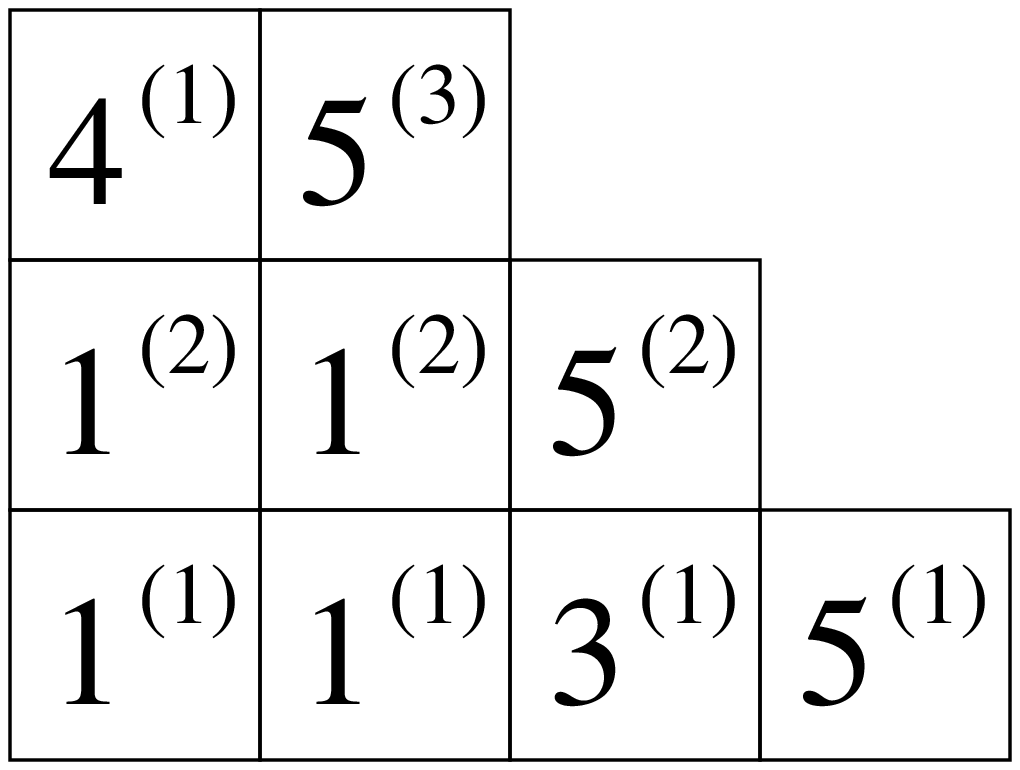}}\;
\overset{\phi}{\hookrightarrow}
\;\raisebox{-0.5cm}{\epsfxsize=1.8 cm \epsffile{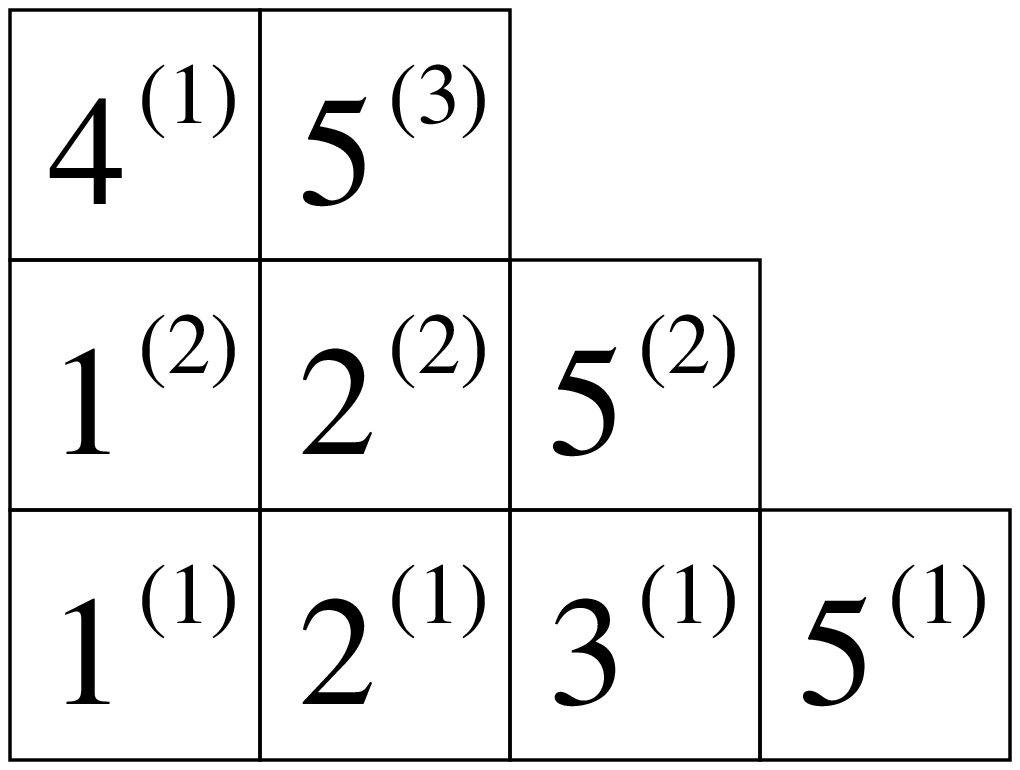}}\;
\xrightarrow{\tau_3}
\;\raisebox{-0.5cm}{\epsfxsize=1.8 cm \epsffile{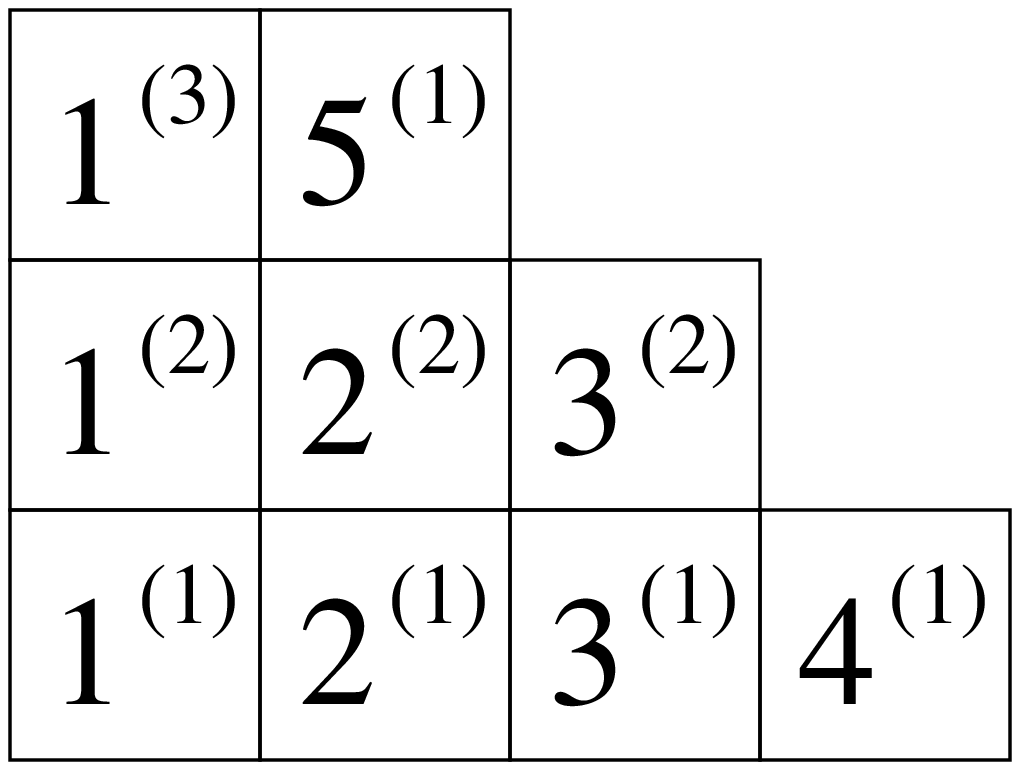}}\;
=\theta(T),
\end{equation*}
where $\tau_1=\omega\circ\sigma_2\circ\sigma_3\circ\sigma_1\circ\omega^{-1}$,
$\tau_2=\omega\circ\sigma_2\circ\sigma_3\circ\sigma_4\circ\sigma_1\circ
\sigma_2\circ\sigma_3\circ\omega^{-1}$ and $\tau_3=\omega\circ\sigma_2\circ
\sigma_3\circ\sigma_1\circ\sigma_2\circ\omega^{-1}$.
\end{example}

\begin{theorem}
\label{theo_poset}
Let $\mu\in\R^L$. Then the cyclage-graph $\T(\mu)$ imposes a ranked poset 
structure on $\ttab{\cdot}{\mu}$ with minimal element $T_{\min}=
[x_1^{\mu_1}\cdots x_L^{\mu_L}]$.
Similarly, $\T'(\mu)$ imposes a ranked poset structure on $\ttab{\cdot}{\mu}$
with minimal element $T_{\max}=[x_L^{\mu_L}\cdots x_1^{\mu_1}]$.
\end{theorem}

\begin{proof}
Let us first consider $\mu$ to be a partition and show that in this case
$\T'(\mu)$ is a ranked poset. For every $T\in\tab{\cdot}{\mu}$ with 
$T\neq T_{\max}$ there exists at least one partition $\la$ such that 
$\Zb'_{\la}(T)\neq 0$ and one can show that
\begin{equation}\label{cT}
c(\Zb'_{\la}(T))=c(T)-1.
\end{equation}
Namely, if $\D'(T)=T$, i.e., $\Zb'_{\la}(T)=\Z_{\la}^{-1}(T)$ then 
$c(\Z_{\la}^{-1}(T))=\|\mu\|-\co(\Z_{\la}^{-1}(T))=\|\mu\|-\co(T)-1=c(T)-1$
by equation~\eqref{ZZ} and the fact that the cocharge is the rank of $\T(\mu)$
for a partition $\mu$ as shown by Lascoux and Sch\"utzenberger~\cite{LS81}.
{}From the explicit prescription for calculating the charge of a Young
tableau $T\in\tab{\cdot}{\mu}$ via indices (see for example~\cite{M95}
page 242 or~\cite{B94} page 111) one may easily check that
$c(T)=c(\D'(T))$ which proves~\eqref{cT}. This shows that for a partition 
$\mu$, $\T'(\mu)$ is a poset ranked by the charge with minimal element 
$T_{\max}$.

{}From lemma~\ref{lem_TT'} and equations~\eqref{stand} and~\eqref{thZ} we 
deduce that also $\T(\mu)$ with $\mu\in\R^L$ is a ranked poset. Since for 
$T_{\max}=\bigl[x_L^{\mu^{\top}_L}\cdots x_1^{\mu^{\top}_1}\bigr]$ with 
$\mu\in\R^L$
\begin{equation}\label{TminTmax}
\La(T_{\max})=T_{\min}=\bigl[x_1^{\mu_1}\cdots x_L^{\mu_L}\bigr],
\end{equation}
the minimal element of $\T(\mu)$ is $T_{\min}$. According to 
lemma~\ref{lem_TT'} also $\T'(\mu)$ is a ranked poset for all $\mu\in\R^L$
with minimal element equal to $T_{\max}$.
\end{proof}

The standardization embedding~\eqref{stand} can be refined by combining
$\phi$ with the action of $\G$ to obtain
\begin{equation}\label{psi}
\psi_{\nu\mu}:\T(\nu)\hookrightarrow\T(\mu), \qquad \nu\geq \mu
\end{equation}
for $\mu,\nu\in\R^L$ with the ordering $\nu\geq \mu$ as defined in
section~\ref{sec_review}. Similar to~\eqref{thZ}
\begin{equation}\label{pZ}
[\psi_{\nu\mu},\Zb_{\la}]=0.
\end{equation}
Certainly, $[\psi_{\nu\mu},\Cb]=0$ thanks to~\eqref{sC} and $[\phi,\Cb]=0$
which can be varified explicitly. To establish~\eqref{pZ} for general
$\Zb_{\la}$ we are left to show $[\phi,\Zb_{\la}]=0$. Let us briefly sketch 
the proof here. Firstly, $\psi_{\nu\mu}$ only depends on $\nu$ and $\mu$, 
but not on its explicit composition in terms of $\phi$ and $\sigma_i$'s. This 
can be shown by induction on the cocharge using $[\psi_{\nu\mu},\Cb]=0$. 
Secondly, $\Zb_{\la}(T)=0$ if and only if $\Zb_{\la}\circ\phi(T)=0$. This can 
be seen as follows. For every LR tableaux there exists a 
standardization composed only of $\phi'$ and $\sigma_i$'s. Denote by 
$\theta_1$ and $\theta_2$ such standardizations for $T$ and $\phi(T)$, 
respectively. Since the standardization is independent of the composition 
of $\phi$ and $\sigma_i$'s we conclude $\theta_1(T)=\theta_2\circ\phi(T)$. 
Thanks to~\eqref{thZ} this means that
\begin{equation*}
\theta_1\circ\Zb_{\la}(T)=\Zb_{\la}\circ\theta_1(T)=\Zb_{\la}\circ
\theta_2\circ\phi(T)=\theta_2\circ\Zb_{\la}\circ\phi(T)
\end{equation*}
which proves the assertion. When $\Zb_{\la}(T)\neq 0$ the commutation
relation $[\Zb_{\la},\phi]=0$ can again be explicitly shown on paths
using the maps $\omega$ and $\omega^{-1}$.

\subsection{Properties of charge and cocharge}
\label{sec_prop}

In this section we establish some properties of the charge and cocharge
for LR tableaux.
The cocharge of an LR tableau $T$ is its rank in the poset
induced by the modified initial cyclage $\Cb$. Since the initial cyclage
is a special $\la$-cyclage the cocharge is also the rank of the poset
$\T(\mu)$. In definition~\ref{def_charge} the charge was defined as 
$c(T)=\|\mu\|-\co(T)$ where $\|\mu\|=\sum_{i<j}|\mu_i\cap\mu_j|$. 
We show now that $\|\mu\|$ is in fact the cocharge of 
$T_{\max}=\bigl[x_L^{\mu_L}\cdots x_1^{\mu_1}\bigr]$ and that the charge of 
an LR tableau is equal to its rank in the poset $\T'(\mu)$.

\begin{lemma}
\label{lem_mu}
For $L\in\Int$ and $\mu\in\R^L$, $\|\mu\|=\co(T_{\max})$.
\end{lemma}

\begin{proof}
Set $P_{\max}:=\omega^{-1}(w_{T{\max}})$ with $\omega^{-1}$ as defined
in~\eqref{omega_inv}. Using theorem~\ref{theo_Hco}, $\|\mu\|=\co(T_{\max})$
is equivalent to $\|\mu\|=H'(P_{\max})$ with $H'$ as in~\eqref{Hp}.
Defining $L_i^{(a)}$ as the number of components of $\mu$ equal to $(i^a)$
and setting $X_{ij}^{ab}:=|(i^a)\cap (j^b)|$, we may rewrite $\|\mu\|$ as
\begin{equation}
\label{X_ijab}
\|\mu\|=\frac{1}{2}\sum_{\substack{1\leq i,j\leq N\\ 1\leq a,b\leq n}} 
L_i^{(a)}(L_j^{(b)}-\delta_{ij}\delta_{ab}) X_{ij}^{ab},
\end{equation}
where $\delta_{ij}=1$ if $i=j$ and zero otherwise.

Let $P\in\O_{P_{\max}}$. Then $P\in\P_{\la\tilde{\mu}}$ where $\tilde{\mu}$
is a permutation of $\mu$. One may easily show that 
$h(p_{k+1}\otimes p_k)=|\tilde{\mu}_k\cap \tilde{\mu}_{k+1}|$ for all 
$P\in\O_{P_{\max}}$. Hence
\begin{equation*}
H'(P_{\max})=\frac{1}{|\O_{P_{\max}}|}\sum_{\tilde{\mu}}
\sum_{k=1}^{L-1}(L-k)|\tilde{\mu}_k\cap \tilde{\mu}_{k+1}|,
\end{equation*}
where the first sum is over all permutations of $\mu$.
Notice that $|\tilde{\mu}_k\cap\tilde{\mu}_{k+1}|=X_{ij}^{ab}$ if
$\tilde{\mu}_k=(i^a)$ and $\tilde{\mu}_{k+1}=(j^b)$.
We now wish to determine the coefficient of $X_{ij}^{ab}$ in $H'(P_{\max})$.
Since 
\begin{equation*}
|\O_{P_{\max}}|=\text{number of permutations of $\mu$}
=\frac{L!}{\prod_{c,k\geq 1} L_k^{(c)}!},
\end{equation*}
the contribution to $X_{ij}^{ab}$ in $H'(P_{\max})$ is
\begin{equation*}
\frac{1}{|\O_{P_{\max}}|}\cdot
\frac{(L-2)!}{
\prod_{c,k\geq 1} (L_k^{(c)}-\delta_{ac}\delta_{ik}-\delta_{bc}\delta_{jk})!}
\sum_{k=1}^{L-1}(L-k)
=\frac{1}{2}L_i^{(a)}(L_j^{(b)}-\delta_{ij}\delta_{ab}).
\end{equation*}
Summing over all $i,j,a,b$ yields~\eqref{X_ijab} which completes the proof.
\end{proof}

The charge and cocharge are dual in the following sense.
\begin{lemma}\label{lem_dual}
For $T\in\ttab{\cdot}{\mu}$, $c(T)=\co(\La(T))$.
\end{lemma}

\begin{proof}
Since $c(T)=\|\mu\|-\co(T)$, the lemma is equivalent to
\begin{equation}\label{co_sum}
\co(\La(T))+\co(T)=\|\mu\|.
\end{equation}
For $T_{\min}=[x_1^{\mu_1}\cdots x_L^{\mu_L}]$ we have $\La(T_{\min})=T_{\max}
=\bigl[x_L^{\mu_L^{\top}}\cdots x_1^{\mu_1^{\top}}\bigr]$. Since 
$\co(T_{\min})=0$ and $\co(T_{\max})=\|\mu^{\star}\|=\|\mu\|$ by 
lemma~\ref{lem_mu}, equation~\eqref{co_sum} holds for $T=T_{\min}$ and 
$T=T_{\max}$. 

Now assume that~\eqref{co_sum} holds for some $T\in\ttab{\cdot}{\mu}$
so that $\D(T)=T$. Then~\eqref{co_sum} also holds for 
$\Zb_{\la}(T)=\Z_{\la}(T)$ if we can show that 
$\co(\La(T))=\co(\La\circ\Z_{\la}(T))-1$ because $\co(T)=\co(\Z_{\la}(T))+1$.
Since $\La\circ\Z_{\la}=\Z_{\la^{\top}}^{-1}\circ\La$ and
$\D\circ\Z_{\la^{\top}}^{-1}=\Z_{\la^{\top}}^{-1}$ this is fulfilled.

If on the other hand $\D'(T)=T$ or equivalently $\D\circ\La(T)=\La(T)$
then~\eqref{co_sum} also holds for $\Zb'_{\la}(T)=\Z_{\la}^{-1}(T)$
because $\co(T)=\co(\Z_{\la}^{-1}(T))-1$ and
$\co(\La(T))=\co(\La\circ\Z_{\la}^{-1}(T))+1$ 
thanks to $\La\circ\Z_{\la}^{-1}=\Z_{\la^{\top}}\circ\La$.

Since $\D(T)=T$ if $\D'(T)\neq T$ and vice versa, unless $T$ is
equal to both $T_{\min}$ and $T_{\max}$, this proves~\eqref{co_sum}
for all $T\in\ttab{\cdot}{\mu}$.
\end{proof}

As argued before the cocharge is the rank of the poset $\T(\mu)$ since
the initial cyclage is a special $\la$-cyclage. Lemmas~\ref{lem_dual} 
and~\ref{lem_TT'} show that the charge is the rank of the poset $\T'(\mu)$.
This is summarized in the following corollary.
\begin{corollary}
For $\mu\in\R^L$, the cocharge is the rank of the poset $\T(\mu)$
and the charge is the rank of the poset $\T'(\mu)$. In addition
\begin{equation*}
0\leq \co(T)\leq \|\mu\| \quad\text{and}\quad 0\leq c(T)\leq \|\mu\|,
\end{equation*}
with $\co(T_{\min})=c(T_{\max})=0$ and $\co(T_{\max})=c(T_{\min})=\|\mu\|$.
\end{corollary}

\section{Properties of the supernomials and generalized Kostka polynomials}
\label{sec_propSK}

Several interesting properties of the supernomials~\eqref{super} and 
generalized Kostka polynomials~\eqref{gK} are stated. In section~\ref{sec_rel}
a duality formula for the generalized Kostka polynomials as well as relations
between the supernomials and the generalized (cocharge) Kostka polynomial
are given.
Recurrences for the A$_{n-1}$ supernomials and the generalized Kostka 
polynomials are established in section~\ref{sec_rec}. These will be used in 
section~\ref{sec_fermi} to obtain a representation of the generalized Kostka 
polynomials of the Kirillov--Reshetikhin-type~\eqref{KRKostka}. 
In section~\ref{sec_A1} we treat the A$_1$ supernomials in more detail and
sketch an elementary proof of the Rogers--Ramanujan-type identities
of ref.~\cite{SW97}.

\subsection{General properties}
\label{sec_rel}

The results of the previous section imply the following
duality formula for the generalized Kostka polynomials.
\begin{theorem}
\label{theo_Kinv}
For $\la$ a partition and $\mu\in\R^L$,
\begin{equation}\label{Kinv}
K_{\la\mu}(q)=q^{\|\mu\|}K_{\la^{\top}\mu^{\star}}(1/q).
\end{equation}
\end{theorem}

\begin{proof}
This follows from the charge representation of the generalized
Kostka polynomials of corollary~\ref{cor_SK}, lemma~\ref{lem_dual}
and $c(T)=\|\mu\|-\co(T)$.
\end{proof}

The supernomial $S_{\lm}(q)$ and the generalized cocharge Kostka polynomial
$\Sb_{\lm}(q)$ satisfy linear relations as follows.
\begin{theorem}
\label{theo_SasSb}
For $\la\in\Int^n$ and $\mu\in\R^L$,
\begin{equation}\label{SasSb}
S_{\lm}(q)=\sum_{\eta\vdash |\la|} K_{\eta\la}\Sb_{\eta\mu}(q),
\end{equation}
where $K_{\eta\la}=K_{\eta\la}(1)$ is the Kostka number.
\end{theorem}

\begin{proof}
By definition the supernomial $S_{\lm}(q)$ is the generating function over
all paths $P\in\P_{\lm}$ weighted by $H(P)$ and by~\eqref{Sb_int} 
$\Sb_{\eta\mu}(q)$ is the generating function over all 
LR tableaux $T\in\ttab{\eta}{\mu}$ with cocharge statistic. 
Hence, since $[\omega(P)]\in\ttab{\cdot}{\mu}$ and
$H(P)=\co(\Om\circ \omega(P))$ by theorem~\ref{theo_Hco} for $P\in\P_{\lm}$, 
equation~\eqref{SasSb} is proven if we can show that for all partitions
$\eta$ of $|\la|$ and $T\in\ttab{\eta}{\mu}$ there are $K_{\eta\la}$
paths such that $[\omega(P)]=T$.

To this end let us show that for all partitions $\eta$ of $|\la|$ with 
$\eta\geq \la$ a pair $(T,t)$ with $T\in\ttab{\eta}{\mu}$ and 
$t\in\tab{\eta}{\la}$ uniquely specifies a path 
$P=p_L\otimes\cdots\otimes p_1\in\P_{\lm}$ by requiring that $p_L\cdot \ldots
\cdot p_1=t$ and $[\omega(P)]=T$. Firstly, by point~\ref{i_1} of 
lemma~\ref{lem_om} indeed $\shape(p_L\cdot\ldots\cdot p_1)
=\shape([\omega(P)])$.
Let us now construct $P\in\P_{\lm}$ from a given pair $(T,t)$. Set $a_i=
\height(\mu_i)$ and define $p_i^{(k)}$ and $t_i^{(k)}$
$(1\leq i\leq L; 1\leq k\leq a_i)$ recursively as follows. Set 
$t_{L+1}^{(1)}=t$ and decompose for $1\leq i\leq L$
\begin{equation}
\label{decom}
t_{i+1}^{(1)}=p_i^{(a_i)}\cdot t_i^{(a_i)}\quad\text{and}
\quad t_i^{(k+1)}=p_i^{(k)}\cdot t_i^{(k)}\quad (1\leq k<a_i)
\end{equation}
such that $\shape(p_i^{(k)})=(\width(\mu_i))$ and $\shape(t_i^{(k)})
=\shape(T_i^{(k)})$, where $T_i^{(k)}$ is obtained from $T$ by dropping all
letters $x\geq x_i^{(k)}$. The decompositions in~\eqref{decom} are unique
by the Pieri formula. The desired path is $P=p_L\otimes\cdots\otimes p_1$ 
where $p_i:=p_i^{(a_i)}\cdot \ldots \cdot p_i^{(1)}$ $(1\leq i\leq L)$ 
because $p_i$ has shape $\mu_i$ since $T\in\ttab{\eta}{\mu}$,
$p_L\cdot \ldots \cdot p_1=t$ and $[\omega(P)]=T$ by construction.
\end{proof}

{}From equation~\eqref{SasSb} one can infer that the special cases 
of the supernomials for which $\mu$ or $\mu^{\star}$ is a partition
have previously occurred in the literature. In the study of finite abelian 
subgroups, Butler~\cite{B87}-\cite{B94} defines polynomials 
$\alpha_{\mu}(S;q)$, 
where $\mu$ is a partition and $S=\{a_1<\cdots<a_{n-1}\}$ an ordered set 
of $n-1$ integers such that $a_{n-1}<|\mu|$, and shows that they satisfy
\begin{equation}
h_{a_1}(x)h_{a_2-a_1}(x)\cdots h_{m-a_{n-1}}(x)
=\sum_{\mu\vdash m} \alpha_{\mu}(S;q^{-1})q^{\|\mu\|} P_{\mu}(x;q).
\end{equation} 
Here $h_k(x)$ is the $k$th homogeneous symmetric function and
$P_{\mu}(x;q)$ is the Hall-Littlewood polynomial.
Using $h_{\la_1}(x)\cdots h_{\la_n}(x)=
\sum_{\eta} K_{\eta \la} s_{\eta}(x)$ and equations \eqref{HL}
and \eqref{SasSb} 
immediately yields that $\alpha_{\mu}(S;q)=S_{\lm}(q)$ where
$\la=(a_1,a_2-a_1,\ldots,|\mu|-a_{n-1})$.
When $\mu^{\star}$ is a partition the supernomial has been studied by
Hatayama {\em et al.}~\cite{HKKOTY98}.

An immediate consequence of theorem~\ref{theo_SasSb} is the inverse of
relation~\eqref{SasSb}.
\begin{corollary}
For $\la$ a partition with $\height(\la)\leq n$ and $\mu\in\R^L$ 
\begin{equation}\label{SbasS}
\Sb_{\lm}(q)=\sum_{\tau\in S_n}\epsilon(\tau)
S_{(\la_1+\tau_1-1,\ldots,\la_n+\tau_n-n)\mu}(q),
\end{equation}
where $S_n$ is the permutation group on $1,2,\ldots,n$ and $\epsilon(\tau)$
is the sign of $\tau$.
\end{corollary}
\begin{proof}
Substitute \eqref{SasSb} into the right-hand side of \eqref{SbasS}
and use (see page~76 of ref.~\cite{F97})
\begin{equation*}
\sum_{\tau\in S_n}\epsilon(\tau)
K_{\eta(\la_1+\tau_1-1,\ldots,\la_n+\tau_n-n)}=\delta_{\eta\la}.
\end{equation*}
\end{proof}

\subsection{Recurrences of the A$_{n-1}$ supernomials and generalized Kostka
polynomials}
\label{sec_rec}

We have seen in equations~\eqref{step_order} and~\eqref{step_orderK} that the 
supernomials and generalized Kostka polynomials are independent of the 
ordering of $\mu$. We may therefore
label the supernomials and generalized Kostka polynomials by a matrix 
$\vL$ with component $L_i^{(a)}$ in row $a$ and column $i$ where 
\begin{equation}\label{Lmu}
L_i^{(a)}=L_i^{(a)}(\mu):=
\text{number of components of $\mu$ equal to $(i^a)$}.
\end{equation} 
If $N:=\max\{\width(\mu_k)\}$ then $\vL$ is an $n\times N$ matrix.
We denote the supernomials and generalized Kostka polynomials with the
label $\vL$ by $S(\vL,\la)$ and $K(\vL,\la)$, respectively, and 
from now on we identify $S(\vL,\la)$ and $S_{\lm}(q)$ (similarly
$K(\vL,\la)$ and $K_{\lm}(q)$) if $\mu$ and $\vL$ are related as 
in~\eqref{Lmu}. Define $\ve_i^{(a)}$ as the $n\times N$ matrix with the
only non-zero element in row $a$ and column $i$ equal to 1 and 
furthermore set $L=\sum_{i=1}^N\sum_{a=1}^n L_i^{(a)}$,
\begin{equation}
\label{l}
\l_i^{(a)}=\sum_{j=1}^N\min\{i,j\}L_j^{(a)}\quad\text{and}\quad
\lb_i^{(a)}=\sum_{b=1}^n\min\{a,b\}L_i^{(b)}.
\end{equation}
With this notation we can state the following recurrence relations.

\begin{theorem}[Recurrences]
\label{theo_rec}
Let $i,a,N,n\in\Int$ such that $1\leq i< N$ and $1\leq a< n$. Let $\vL$ be 
an $n\times N$ matrix with non-negative integer
components such that $L_i^{(a)}\geq 2$. Then for $\la\in\Int^n$ 
\begin{align}
\label{recS}
S(\vL,\la)&=S(\vL+\ve_{i-1}^{(a)}-2\ve_i^{(a)}+\ve_{i+1}^{(a)},\la)
+q^{\l_i^{(a)}-i}S(\vL+\ve_i^{(a-1)}-2\ve_i^{(a)}+\ve_i^{(a+1)},\la)\\
\intertext{and for $\la$ a partition}
\label{recK}
K(\vL,\la)&=q^{\lb_i^{(a)}-a}
K(\vL+\ve_{i-1}^{(a)}-2\ve_i^{(a)}+\ve_{i+1}^{(a)},\la)
+K(\vL+\ve_i^{(a-1)}-2\ve_i^{(a)}+\ve_i^{(a+1)},\la).
\end{align}
\end{theorem}

\begin{proof}
First we prove~\eqref{recS}. Take $\mu\in\R^L$ corresponding to $\vL$
such that $\mu_{L-1}=\mu_L=(i^a)$ (which is possible since $L_i^{(a)}\geq 2$).
Define $\mu'$ and $\mu''$ by $\mu_L'=((i+1)^a)$, $\mu_{L-1}'=((i-1)^a)$,
$\mu_L''=(i^{a+1})$, $\mu_{L-1}''=(i^{a-1})$ and $\mu_j'=\mu_j''=\mu_j$ for
$1\leq j\leq L-2$. Recalling definition~\ref{def_super} of the supernomials,
it is obvious that $\P_{\lm'}$ and $\P_{\lm''}$ are the sets of 
paths underlying the two terms on the right-hand side of~\eqref{recS}.
Furthermore $\P_{\lm'}$ and $\P_{\lm''}$ are disjoint. We now wish to 
establish a bijection between $\P_{\lm}$ and $\P_{\lm'}\cup\P_{\lm''}$.
To this end define $\tau(p_L\otimes p_{L-1})=\tilde{p}_L\otimes 
\tilde{p}_{L-1}$ for $p_{L-1},p_L\in\B_{\mu_L}$ such that 
\begin{equation}
\label{prod}
\tilde{p}_L\cdot \tilde{p}_{L-1}=p_L\cdot p_{L-1}
\end{equation}
and either (a) $\tilde{p}_{L-1}\in\B_{\mu_{L-1}'}$, 
$\tilde{p}_L\in\B_{\mu_L'}$ if $\nu\cap ((i+1)^a)=((i+1)^a)$ or 
(b) $\tilde{p}_{L-1}\in\B_{\mu_{L-1}''}$, $\tilde{p}_L\in\B_{\mu_L''}$
if $\nu\cap (i^{a+1})=(i^{a+1})$ where $\nu=\shape(p_L\cdot p_{L-1})$.
Indeed these conditions are mutually excluding and determine 
$\tilde{p}_{L-1}$ and $\tilde{p}_L$ uniquely, i.e.,
the Littlewood-Richardson coefficient $c_{\mu_{L-1}'\mu_L'}^{\nu}=1$ if and 
only if $c_{\mu_{L-1}''\mu_L''}^{\nu}=0$ and vice versa.
Conversely, if $\tilde{p}_{L-1}\in\B_{\mu_{L-1}'}$, $\tilde{p}_L\in\B_{\mu_L'}$
(or $\tilde{p}_{L-1}\in\B_{\mu_{L-1}''}$, $\tilde{p}_L\in\B_{\mu_L''}$)
one can find unique $p_L\otimes p_{L-1}=\tau^{-1}(\tilde{p}_L\otimes 
\tilde{p}_{L-1})$ with $p_{L-1},p_L\in\B_{\mu_L}$ by requiring~\eqref{prod}.
Hence $\tau:\P_{\lm}\to\P_{\lm'}\cup\P_{\lm''}$ with
$\tau(P):=\tau(p_L\otimes p_{L-1})\otimes p_{L-2}\otimes\cdots\otimes p_1$
for each path $P=p_L\otimes\cdots\otimes p_1\in\P_{\lm}$, is the 
desired bijection. This proves~\eqref{recS} at $q=1$.

To prove~\eqref{recS} at arbitrary base $q$ notice that if 
$\tau(P)\in\P_{\lm'}$, then the LR tableaux $T=[\Om\circ\omega(P)]$ and
$T'=[\Om\circ\omega\circ\tau(P)]$ are related as
\begin{equation*}
T'=\psi_{\mu\mu'}(T),
\end{equation*}
with $\psi_{\mu\mu'}$ defined in~\eqref{psi}. Because of~\eqref{pZ}
we have $\co(T)=\co(T')$. Hence theorem~\ref{theo_Hco} implies that also 
$H(P)=H(\tau(P))$ for all $P$ such that $\tau(P)\in\P_{\lm'}$. Therefore, 
the term $S(\vL+\ve_{i-1}^{(a)}-2\ve_i^{(a)}+\ve_{i+1}^{(a)},\la)$ 
in~\eqref{recS} comes without a power of $q$.

Similarly, if $\tau(P)\in\P_{\lm''}$ then the LR tableaux 
$T=[\Om\circ\omega(P)]$ and $T''=[\Om\circ\omega\circ\tau(P)]$ are related 
as $\La(T'')=\psi_{\mu\mu''}\circ\La(T)$ which implies
\begin{equation}\label{coLa}
\co(\La(T))=\co(\La(T'')).
\end{equation}
Therefore, pulling all strings in our register, we derive
\begin{align*}
H(P)-H(\tau(P))&=\co(T)-\co(T'') &
&\text{by theorem~\ref{theo_Hco}}\\
&=\|\mu\|-\co(\La(T))-\|\mu''\|+\co(\La(T'')) & 
&\text{by lemma~\ref{lem_dual}}\\
&=\|\mu\|-\|\mu''\| & &\text{by equation~\eqref{coLa}}\\
&=\l_i^{(a)}-i &
&\text{recalling $\|\mu\|=\sum_{j<k}|\mu_j\cap\mu_k|$}
\end{align*}
which is the power of $q$ in front of the second term in~\eqref{recS}.
This concludes the proof of~\eqref{recS}.

To prove~\eqref{recK}, recall definition~\ref{def_gKp} of the generalized
Kostka polynomials. The generalized cocharge Kostka polynomials~\eqref{Sb} 
obey the same recurrences~\eqref{recS} as the supernomials. This follows 
from the fact that $\Pb_{\lm}\subset\P_{\lm}$ and that $\tau(P)$ is in 
$\Pb_{\lm'}$ or $\Pb_{\lm''}$ if $P\in\Pb_{\lm}$ by the same arguments as 
in the proof of point~\ref{i_2} of lemma~\ref{lem_om}.
Using~\eqref{gK}, $\|\mu\|-\|\mu'\|=\lb_i^{(a)}-a$ and
$\|\mu\|-\|\mu''\|=\l_i^{(a)}-i$ yields~\eqref{recK}.
\end{proof}

Rectangular Young tableaux over the alphabet $\{1,2,\ldots,n\}$ of height
$n$ are often identified with the empty tableau.
When this identification is made for the steps of the paths in the generating
functions defining the supernomials and generalized Kostka polynomials
one obtains the following properties.

\begin{lemma}
\label{lem_rel}
Let $n\geq 2$, $N,i\geq 1$ be integers and $\vL$ an $n\times N$ matrix with
non-negative entries. Then for $\la\in\Int^n$
\begin{align}
\label{Sn}
S(\vL+\ve_i^{(n)},\la+(i^n))&=S(\vL,\la),\\
\intertext{and for $\la$ a partition with at most $n$ parts}
\label{Kn}
K(\vL+\ve_i^{(n)},\la+(i^n))&=q^{\sum_a a\l_i^{(a)}}K(\vL,\la).
\end{align}
\end{lemma}

\begin{proof}
Writing $S(\vL+\ve_i^{(n)},\la+(i^n))$ as a generating function over paths
as in~\eqref{super}, each path $P$ in the sum has at least one step of
height $n$. Hence $\height(\omega(P))=n$ by lemma~\ref{lem_om}.
Denoting by $P'$ the path obtained from $P$ by dropping the step $p_k$ of
shape $(i^n)$, we find from lemma~\ref{lem_drop} and $\Om(p_k)=p_k$ that
$H(P)=H(P')$. This proves~\eqref{Sn}.

The generalized cocharge Kostka polynomials $\Sb$ obey the same 
relation~\eqref{Sn} as $S$. 
Let $\mu$ and $\mu'$ be any of the arrays of rectangular partitions 
corresponding to $\vL+\ve_i^{(n)}$ and $\vL$, respectively, by~\eqref{Lmu}.
Then using $\|\mu\|-\|\mu'\|=\sum_a a\l_i^{(a)}$ and recalling~\eqref{gK}
one finds~\eqref{Kn}.
\end{proof}

\subsection{The A$_1$ supernomials}
\label{sec_A1}

The A$_1$ supernomials are given by specializing definition~\ref{def_super}
to $n=2$, i.e., $\la=(\la_1,\la_2)\in\Int^2$. By lemma~\ref{lem_rel}
it is sufficient to label all A$_1$ supernomials by a vector $\vL\in\Int^N$ 
instead of a two-row matrix. Recall that the supernomials vanish unless
$\la_1+\la_2=|\mu|=\ell_N$ where $\ell_i=\sum_{j=1}^N \min\{i,j\} L_j$. 
We therefore set
\begin{equation}\label{S_A1}
S_1(\vL,a):=S_{\lm}(q)=\sum_{P\in \P_{\lm}} q^{H(P)}
\end{equation}
where $\vL\in\Int^N$, $a+\frac{1}{2}\ell_N=0,1,\ldots,\ell_N$ and
\begin{equation}
\label{lm}
\begin{aligned}
\mu&=(1^{L_1}2^{L_2}\cdots N^{L_N}),\\
\la&=(\tfrac{1}{2}\ell_N+a,\tfrac{1}{2}\ell_N-a).
\end{aligned}
\end{equation}

For $S_1(\vL,a)$ the recurrences~\eqref{recS} read
\begin{equation}
\label{rec_S1a}
S_1(\vL,a)=S_1(\vL+\ve_{i-1}-2\ve_i+\ve_{i+1},a)
+q^{\ell_i-i}S_1(\vL-2\ve_i,a),
\end{equation}
for $1\leq i<N$. The $\ve_i$ are the canonical basis vectors of $\Integer^N$
and $\ve_0=0$. The above equation is in fact part of a larger family of 
recursion relations.

\begin{lemma}
\label{lem_rec_S1}
Let $A,B,N$ be integers such that $1\leq A\leq B<N$ and let $\vL\in\Int^N$
such that $L_1=\cdots=L_{B-1}=0$ if $A<B$. Then for
$a+\tfrac{1}{2}\ell_N=0,1,\ldots,\ell_N$
\begin{equation}
\label{rec_S1}
S_1(\vL+\ve_A+\ve_B,a)=S_1(\vL+\ve_{A-1}+\ve_{B+1},a)+q^{\ell_A+A}
S_1(\vL+\ve_{B-A},a).
\end{equation}
\end{lemma}
\begin{proof}
When $A=B$ equation~\eqref{rec_S1} reduces to~\eqref{rec_S1a} with $\vL$
replaced by $\vL+2\ve_i$. For $A<B$
equation~\eqref{rec_S1} follows from the recurrences
\begin{equation}\label{rec_int}
S(\vL+\ve_A^{(1)}+\ve_B^{(1)},\la)=S(\vL+\ve_{A-1}^{(1)}+\ve_{B+1}^{(1)},\la)
+q^{\ell_A^{(1)}+A}S(\vL+\ve_{B-A}^{(1)}+\ve_A^{(2)},\la)
\end{equation}
where $n=2$ and $\vL$ is a $2\times N$ matrix such that 
$L_1^{(1)}=\cdots=L_{B-1}^{(1)}=0$. This can be seen by dropping
all entries of the matrix in the second row using lemma~\ref{lem_rel} 
and then replacing the matrix $\vL$ by a vector $\vL$.

Equation~\eqref{rec_int} can be proven in complete analogy to the proof of 
theorem~\ref{theo_rec} as follows. Take $\mu$ corresponding to 
$\vL+\ve_A^{(1)}+\ve_B^{(1)}$, replace the components 
$\mu_{L-1},\ldots,\mu_L''$ in the proof of theorem~\ref{theo_rec} by
$\mu_{L-1}=(B),\mu_L=(A),\mu_{L-1}'=(B+1),\mu_L'=(A-1),\mu_{L-1}''=(B-A),
\mu_L''=(A,A)$ and set $\tau(p_L\otimes p_{L-1})=\tilde{p}_L\otimes
\tilde{p}_{L-1}$ where again $\tilde{p}_L\cdot\tilde{p}_{L-1}
=p_L\cdot p_{L-1}$ and now (a) $\tilde{p}_L\in\B_{\mu_L'}$,
$\tilde{p}_{L-1}\in\B_{\mu_{L-1}'}$ if $\shape(p_L\cdot p_{L-1})\neq (B,A)$
and (b) $\tilde{p}_L\in\B_{\mu_L''}$,
$\tilde{p}_{L-1}\in\B_{\mu_{L-1}''}$ if $\shape(p_L\cdot p_{L-1})=(B,A)$.
Note that we have used here that $n=2$ which ensures that the shape
of the product of two steps has at most height 2. One may explicitly 
check that $\co(T)=\co(T')$ and $\co(T)=\co(T'')+\ell_A^{(1)}+A$ with
$T,T',T''$ as defined in the proof of theorem~\ref{theo_rec} which 
proves~\eqref{rec_int}.
\end{proof}

Using $S_{\lm}(q)=\alpha_{\mu}(S;q)$ for $\mu$ a partition (see
the discussion after theorem~\ref{theo_SasSb})
and the explicit representations for 
$\alpha_{\mu}(S;q)$ in \cite{B87}-\cite{B94}, one finds that
\begin{equation}\label{super_rel}
S_1(\vL,a)=\qbin{\vL}{a}
\end{equation}
with $\qbin{\vL}{a}$ given in equation~\eqref{qsuper}. 

We now recall some identities of ref.~\cite{SW97} involving the A$_1$ 
supernomial and show how the recurrences of lemma \ref{lem_rec_S1} yield an
elementary proof. 
The identities unify and extend many
of the known Bose--Fermi or Rogers--Ramanujan-type identities for
one-dimensional configuration sums of solvable lattice models.
Below we only quote the result of~\cite{SW97} corresponding to
the Andrews--Baxter--Forrester model and its fusion hierarchy.
Set
\begin{equation}
\label{T1}
T_1(\vL,a)=q^{\frac{1}{4}\sum_{i=1}^NL_i\ell_i-\frac{a^2}{N}}
\qbin{\vL}{a}_{1/q}.
\end{equation}

\begin{theorem}\label{theo_unitary}
Let $a,b,p,N$ be integers such that $N<p-2$, $1\leq a \leq p-1$ and 
$1\leq b\leq p-N-1$ and let $\vL\in\Int^N$. Then
\begin{multline}\label{unitary}
\sum_{j=-\infty}^{\infty}\Bigl\{
q^{\frac{j}{N}\left(p(p-N)j+p b-(p-N)a\right)}
T_1\bigl(\vL,\tfrac{b-a}{2}+pj\bigr)
-q^{\frac{1}{N}(pj+a)((p-N)j+b)}
T_1\bigl(\vL,\tfrac{b+a}{2}+pj\bigr)\Bigr\} \\
=q^{\frac{1}{4N}(b-a)(a-b-N)} \hspace{-3mm}
\sum_{\substack{\vm\in \Int^{\,p-3} \\
\vm \equiv \vQ_{ab}~(\bmod{2})}} \hspace{-3mm}
q^{\frac{1}{4}\vm C\vm-\frac{1}{2}m_{a-1}}
\prod_{j=1}^{p-3} \qbin{m_j+n_j}{m_j},
\end{multline}
where $C$ is the Cartan matrix of A$_{p-3}$ and 
$\vQ_{ab}=\vQ^{(a-1)}+\vQ^{(p-b-1)}+\vQ^{(p-2)}+\sum_{i=2}^N L_i\vQ^{(i)}$ with
$\vQ^{(i)}=\ve_{i-1}+\ve_{i-3}+\cdots$. The expression 
$\vm\equiv \vQ~(\bmod 2)$ stands for $m_i\equiv Q_i~(\bmod 2)$ and
$\vm C\vm=\sum_{i,j=1}^{p-3}C_{ij}m_im_j$.
The variable $m_0=0$ and $\vn$ is determined by
\begin{equation*}
\vn=\frac{1}{2}(\ve_{a-1}+\ve_{p-b-1}+\sum_{i=1}^N L_i\ve_i-C\vm).
\end{equation*}
\end{theorem}
\noindent
For $\vL=L\ve_{N'}$, $1\leq N'\leq N$, in the limit $L\to \infty$ the 
identity~\eqref{unitary} yields an identity for branching functions of 
A$_1^{(1)}$ cosets.

If we can show that the (suitably normalized) $q\to 1/q$ forms of both sides 
of \eqref{unitary} satisfy the recurrence
\begin{equation}\label{rec_X}
X(\vL+\ve_A+\ve_B)=X(\vL+\ve_{A-1}+\ve_{B+1})+q^{\ell_A+A}X(\vL+\ve_{B-A})
\end{equation}
then the identity is proven if it holds
for the trivial initial conditions $\vL=\ve_i$ ($i=0,1,\ldots,N$).
The $q\to 1/q$ version of the left-hand side of \eqref{unitary}
satisfies \eqref{rec_X} by lemma \ref{lem_rec_S1}. For the right-hand
side of \eqref{unitary} with $q\to 1/q$ it is readily shown \cite{SW97} that
\eqref{rec_X} holds for $A=B$ and all $\vL\in\Integer^N$ by using
modified $q$-binomials obtained by extending the range of $\la_2$ in 
the top line of \eqref{qbin} to $\la_2\in\Integer$. This implies
\eqref{rec_X} thanks to the following lemma.

\begin{lemma}\label{lem_X}
Let $X$ be a function of $\vL\in\Integer^N$ satisfying the recurrences
\begin{equation}\label{rec_X1}
X(\vL+2\ve_A)=X(\vL+\ve_{A-1}+\ve_{A+1})+q^{\ell_A+A}X(\vL)
\end{equation}
for all $A=1,2,\ldots,N-1$ and $\vL\in\Integer^N$. Then $X$ fulfills
the more general recurrences
\begin{equation}\label{rec_X2}
X(\vL+\ve_A+\ve_B)=X(\vL+\ve_{A-1}+\ve_{B+1})+q^{\ell_A+A}X(\vL+\ve_{B-A})
\end{equation}
for all $1\leq A\leq B<N$ and $\vL\in\Integer^N$ such that 
$L_1=\cdots=L_{B-1}=0$ if $A<B$.
\end{lemma}

\begin{proof}
Assume that~\eqref{rec_X2} is proven for all $1\leq A'\leq B'<B$. This is 
certainly true for $A'=B'=1$ thanks to~\eqref{rec_X1}. Using~\eqref{rec_X1}
successively (with $A$ replaced by $i$) for $i=B,B-1,\ldots,A$ yields
\begin{equation}
\label{X_int1}
X(\vL+\ve_A+\ve_B)=X(\vL+\ve_{A-1}+\ve_{B+1})
+\sum_{i=A}^B q^{\ell_i+A}X(\vL+\ve_A-\ve_i-\ve_{i+1}+\ve_{B+1}).
\end{equation}
Applying~\eqref{rec_X1} with $A$ replaced by $B$ to the second term 
on the right-hand side in~\eqref{X_int1} one obtains
\begin{multline}\label{X_int}
\sum_{i=A}^{B-1}q^{\ell_i+A}\bigl\{
X(\vL+\ve_A-\ve_i-\ve_{i+1}-\ve_{B-1}+2\ve_B)
-q^{\ell_B-2i+A}X(\vL+\ve_A-\ve_i-\ve_{i+1}-\ve_{B-1})\bigr\}\\
+q^{\ell_B+A}X(\vL+\ve_A-\ve_B).
\end{multline}
Telescoping the last term with the negative terms in the sum at 
$i=B-1,B-2,\ldots,A$ using~\eqref{rec_X1} with $A\to i-1$ yields
$q^{\ell_B+A}X(\vL+\ve_{A-1}-\ve_{B-1})$. The positive terms in the
sum in~\eqref{X_int} can be simplified to
$q^{\ell_A+A}X(\vL+\ve_{B-A-1}-\ve_{B-1}+\ve_B)$ by combining successively
the term $i=A$ with $i=A+1,\ldots,B-1$ using~\eqref{rec_X2}
with $A\to i-A$ and $B\to i$. Therefore~\eqref{X_int1} becomes
\begin{multline*}
X(\vL+\ve_A+\ve_B)=X(\vL+\ve_{A-1}+\ve_{B+1})\\
+q^{\ell_A+A}X(\vL+\ve_{B-A-1}-\ve_{B-1}+\ve_B)+q^{\ell_B+A}
X(\vL+\ve_{A-1}-\ve_{B-1}).
\end{multline*}
The last two terms can be combined to $q^{\ell_A+A}X(\vL+\ve_{B-A})$
employing~\eqref{rec_X2} with $A\to B-A$ and $B\to B-1$ and using
that $L_1=\cdots =L_{B-1}=0$. This yields~\eqref{rec_X2}.
\end{proof}

Let us make some comments about the above outlined proof of theorem
\ref{theo_unitary}. First, lemma \ref{lem_rec_S1} requires $N>1$. 
However, thanks to 
$T_1\bigl((L_1,\ldots,L_N,0,\ldots,0),a\bigr)=q^{\frac{M-N}{MN}a^2}
T_1\bigl((L_1,\ldots,L_N),a\bigr)$,
where the dimension of the vector on the left-hand side is $M$, one can
derive the identities~\eqref{unitary} for all $N\geq 1$ except when $p=4$.
Second, we note that for $\vL\in\Int^N$ the polynomials on the right-hand side
of \eqref{unitary} indeed remain unchanged by replacing the $q$-binomial
with the modified $q$-binomial.
Finally, the proof given in~\cite{SW97} used the identities at $\vL=L\ve_1$ 
as initial conditions. The knowledge of these non-trivial
identities is not necessary in the above proof.

In the discussion section we will conjecture higher-rank analogues 
of~\eqref{unitary}.

\section{Fermionic representation of the generalized Kostka polynomials}
\label{sec_fermi}

In this section we give a fermionic representation of the generalized Kostka
polynomials generalizing the Kirillov--Reshetikhin expression~\eqref{KRKostka}.
Recalling the definitions~\eqref{l} we introduce the following 
function.
\begin{definition}\label{def_F}
Let $n\geq 2$ and $N\geq 1$ be integers, $\la$ a partition with 
$\height(\la)\leq n$ and $\vL$ an $n\times N$ matrix with entry $L_i^{(a)}\in
\Int$ in row $a$ and column $i$. Then set $F(\vL,\la)=0$ if $|\la|\neq
\sum_{a,i\geq 1}ai L_i^{(a)}$ and otherwise
\begin{equation}\label{F}
F(\vL,\la)=\sum_{\alpha} q^{C(\alpha)} \prod_{a,i\geq 1}
\qbin{P_i^{(a)}+\alpha_i^{(a)}-\alpha_{i+1}^{(a)}}{\alpha_i^{(a)}
-\alpha_{i+1}^{(a)}},
\end{equation}
where the sum is over sequences $\alpha=(\alpha^{(1)},\alpha^{(2)},\ldots)$
of partitions such that $|\alpha^{(a)}|=\sum_{j\geq 1}j\lb_j^{(a)}-(\la_1
+\cdots+\la_a)$. Furthermore, with the convention that $\alpha^{(0)}_i=0$,
\begin{align}
P_i^{(a)}&=\sum_{k=1}^i(\alpha_k^{(a-1)}-2\alpha_k^{(a)}+\alpha_k^{(a+1)})
+\l_i^{(a)},\\
C(\alpha)&=\sum_{a,i\geq 1}\binom{A_i^{(a)}+\alpha_i^{(a-1)}
-\alpha_i^{(a)}}{2},
\quad A_i^{(a)}=\sum_{\substack{k\geq i\\ b\geq a}} L_k^{(b)}.
\end{align}
\end{definition}

Recalling that $K_{\lm}(q)=0$ unless $|\la|=|\mu|=\sum_{a,i\geq 1}aiL_i^{(a)}$
we find that $F(\vL,\la)=K(\vL,\la)$ if $L_i^{(a)}=0$ for $a>1$ by
comparing~\eqref{F} with~\eqref{KRKostka}. We wish to show that $F(\vL,\la)$
equals the generalized Kostka polynomial $K(\vL,\la)$ for more general $\vL$. 
We begin by showing that $F$ obeys the same recurrence relation as $K$.

\begin{lemma}\label{lem_recF}
Let $i,a,N,n\in\Int$ such that $1\leq i< N$ and $1\leq a< n$ and let 
$\la$ be a partition with $\height(\la)\leq n$. Let $\vL$ be an $n\times N$ 
matrix with non-negative integer entries such that $L_i^{(a)}\geq 2$. Then
\begin{equation}
\label{recF}
F(\vL,\la)=q^{\lb_i^{(a)}-a}
F(\vL+\ve_{i-1}^{(a)}-2\ve_i^{(a)}+\ve_{i+1}^{(a)},\la)
+F(\vL+\ve_i^{(a-1)}-2\ve_i^{(a)}+\ve_i^{(a+1)},\la).
\end{equation}
\end{lemma}

\begin{proof}
Under the substitution 
$\vL\to \vL+\ve_{i-1}^{(a)}-2\ve_i^{(a)}+\ve_{i+1}^{(a)}$
the variable $P_j^{(b)}$ and the function $C(\alpha)$ transform as
\begin{equation}\label{trans1}
\begin{aligned}
P_j^{(b)}&\to P_j^{(b)}-\delta_{ij}\delta_{ab},\\
C(\alpha)&\to C(\alpha)-\lb_i^{(a)}+a+\alpha_i^{(a)}-\alpha_{i+1}^{(a)}.
\end{aligned}
\end{equation}
On the other hand, replacing 
$\vL\to\vL+\ve_i^{(a-1)}-2\ve_i^{(a)}+\ve_i^{(a+1)}$ induces the changes
\begin{equation}\label{trans2}
\begin{aligned}
P_j^{(b)}&\to P_j^{(b)}
+\min\{i,j\}(\delta_{a-1,b}-2\delta_{ab}+\delta_{a+1,b}),\\
C(\alpha)&\to C(\alpha)-m_i^{(a)}+i.
\end{aligned}
\end{equation}

Now apply the $q$-binomial recurrence
\begin{equation}\label{recqbin}
\qbin{m+n}{n}=q^n\qbin{m+n-1}{n}+\qbin{m+n-1}{n-1},
\end{equation}
to the $(a,i)$th term in the product in~\eqref{F} (this term cannot be
$\qbin{0}{0}$ because of the condition $L_i^{(a)}\geq 2$).
Thanks to~\eqref{trans1} one can immediately recognize the first term
of the resulting expression as
$q^{\lb_i^{(a)}-a}F(\vL+\ve_{i-1}^{(a)}-2\ve_i^{(a)}+\ve_{i+1}^{(a)},\la)$.
In the second term we perform the variable change
$\alpha_j^{(b)}\to \alpha_j^{(b)}+\chi(j\leq i)\delta_{ab}$ where recall
that $\chi({\rm true})=1$ and $\chi({\rm false})=0$.
Since this leads to exactly the same change in $P_j^{(b)}$ and $C(\alpha)$ as 
in~\eqref{trans2}, the second term indeed yields
$F(\vL+\ve_i^{(a-1)}-2\ve_i^{(a)}+\ve_i^{(a+1)},\la)$.
\end{proof}

\begin{theorem}\label{theo_KF}
Let $N\geq 1,n\geq 2$ be integers, $\la$ a partition with 
$\height(\la)\leq n$ and $\vL$ an $n\times N$ matrix with components 
$L_i^{(a)}\in\Int$ in row $a$ and column $i$. If either
$L_i^{(a)}\geq L_i^{(a+2)}$ for all $1\leq a\leq n-2$
and $1\leq i\leq N$ or $L_i^{(a)}\geq L_{i+2}^{(a)}$ for all $1\leq a\leq n$
and $1\leq i\leq N-2$, then $F(\vL,\la)=K(\vL,\la)$.
\end{theorem}

\begin{proof}
We use $F(\vL,\la)=K(\vL,\la)$ for $\vL$ such that $L_i^{(a)}=0$ when $a>1$ 
as initial condition. Since $K$ and $F$ both satisfy the recurrences
\begin{equation}\label{KFrec}
X(\vL+\ve_{i+1}^{(a+1)})=X(\vL-\ve_i^{(a-1)}+2\ve_i^{(a)})
-q^{\lb_i^{(a)}+1}X(\vL-\ve_i^{(a-1)}+\ve_{i-1}^{(a)}+\ve_{i+1}^{(a)})
\end{equation}
(compare with~\eqref{recK} and~\eqref{recF}, respectively)
the theorem follows immediately for the first set of restrictions on $\vL$.
The second set of restrictions comes about by using the symmetry~\eqref{Kinv}
of the generalized Kostka polynomials.
\end{proof}

The recurrences~\eqref{KFrec} are not sufficient to prove 
theorem~\ref{theo_KF} for $\vL$ with arbitrary entries $L_i^{(a)}\in\Int$. 
However, we nevertheless believe the theorem to be true for this case as well.

\begin{conjecture}
Let $N\geq 1,n\geq 2$ be integers, $\la$ a partition with 
$\height(\la)\leq n$ and $\vL$ an $n\times N$ matrix with nonnegative integer
entries. Then $F(\vL,\la)=K(\vL,\la)$.
\end{conjecture}

\section{Discussion}
\label{sec_dis}

We believe that there exist many further results for the generalized Kostka 
polynomials and supernomials. For example, \eqref{unitary} admits higher-rank 
analogues in terms of
\begin{equation}
\label{T_super}
T(\vL,\la)=q^{\frac{1}{2}\sum_{i=1}^N\sum_{a,b=1}^{n-1}L_i^{(a)}
C^{-1}_{ab}\l_i^{(b)}-\frac{1}{2N}\sum_{i=1}^n (\la_i-\frac{1}{n}|\la|)^2}
S(\vL,\la)|_{1/q},
\end{equation}
where $C^{-1}$ is the inverse of the Cartan matrix of A$_{n-1}$.
For integers $n\ge 2, N,p\ge 1$ such that $N<p-n$ and any 
$n\times N$ matrix $\vL$ with non-negative integer entries such that 
$\sum_{b=1}^{n-1}C_{ab}^{-1}L_i^{(b)}\in\Integer$ for all $i$ and $a$, we 
conjecture
\begin{equation}
\label{AnRR}
\sum_{k_1+\cdots+k_n=0} 
\sum_{\tau\in S_n}  \epsilon(\tau)
q^{\sum_{i=1}^n\{
\frac{1}{2N}(pk_i+\tau_i-i)^2-\frac{p}{2}k_i^2+ik_i\}}
T(\vL,\la(\vk,\tau))
=\sum_{\vm} q^{\frac{1}{2}\vm(C^{-1}\otimes C)\vm}
\qbin{\vm+\vn}{\vm},
\end{equation}
where the following notation has been used. On the left-hand side the 
components of $\la(\vk,\tau)$ are given by 
$\la_j(\vk,\tau)=\frac{1}{n}\sum_{a,i\geq 1}aiL_i^{(a)}+pk_j+\tau_j-j$.
On the right-hand side the sum runs over all 
$\vm=\sum_{a=1}^{n-1}\sum_{i=1}^{p-n-1}m_i^{(a)} (\ve_a\otimes \ve_i)$
with $m_i^{(a)}\in\Int$ such that
$\sum_{b=1}^{n-1} C_{ab}^{-1}m_i^{(b)}\in\Integer$ for all
$a=1,\dots,n-1$ and $i=1,\dots,p-n-1$. 
The variable $\vn$ is fixed by
\begin{equation}
(C\otimes I)\vn+(I\otimes C)\vm=\sum_{a=1}^{n-1}\sum_{i=1}^N
L_i^{(a)}(\ve_a\otimes \ve_i)
\end{equation}
where $I$ is the identity matrix and $C$ is the Cartan matrix of an A-type
Lie algebra. The dimension of the first space in the tensor product is $n-1$
and that of the second space is $p-n-1$. Finally we used the notation
\begin{align*}
(A\otimes B)\vm &=\sum_{a,b=1}^{n-1}\sum_{i,j=1}^{p-n-1}
A_{ab}B_{ij} m_j^{(b)}(\ve_a\otimes \ve_i),\\
\vn(A\otimes B)\vm &=\sum_{a,b=1}^{n-1}\sum_{i,j=1}^{p-n-1}A_{ab}
B_{ij} n_i^{(a)}m_j^{(b)}\\
\intertext{and}
\qbin{\vm+\vn}{\vm}&=\prod_{a=1}^{n-1}\prod_{i=1}^{p-n-1}
\qbin{m_i^{(a)}+n_i^{(a)}}{m_i^{(a)}}.
\end{align*}

The identities~\eqref{AnRR} are polynomial analogues of branching function 
identities of the Rogers--Ramanujan type for A$^{(1)}_{n-1}$ cosets. 
For $n=2$ they follow from theorem~\ref{theo_unitary} with $a=b=1$ and
for $\vL=L(\ve_1\otimes \ve_1)$ they were claimed in~\cite{D97}
\footnote{The proof in~\cite{D97} seems to be incomplete.}.
Unfortunately, the recurrences of theorem~\ref{theo_rec} are not
sufficient to prove~\eqref{AnRR} for general $n$ and $\vL$. 
A proof would require a more complete set of recurrences for the A$_{n-1}$ 
supernomials analogous to those stated in lemma~\ref{lem_rec_S1} for $n=2$.

The left-hand side of equation~\eqref{AnRR} can be interpreted in terms of 
paths of a level-$(p-n)$ A$_{n-1}^{(1)}$ lattice model. Denote by $\La_k$ 
($0\leq k\leq n-1$) the dominant integral weights of A$_{n-1}^{(1)}$. Then
the states $a$ of the lattice model underlying~\eqref{AnRR} are given by the 
level-$(p-n)$ dominant integral weights, i.e., $a=\sum_{k=0}^{n-1} a_k\La_k$ 
such that $\sum_{k=0}^{n-1}a_k=p-n$. Define the adjacency matrices 
$\A$ labelled by two states $a,b$ and a Young tableau as
$\A^{\emptyset}_{a,b}=\chi(a=b)$, $\A^{[i]}_{a,b}=\chi(b=a+\La_i-\La_{i-1})$ 
$(i=1,\ldots,n;\La_n=\La_0)$ and recursively
$\sum_{b}\A^{T}_{a,b}\A^{[i]}_{b,c}=\A^{T\cdot [i]}_{a,c}$. 
Call a path $P=p_L\otimes\cdots\otimes p_1\in\P_{\lm}$ admissible with initial
state $a^{(1)}$ if $\prod_{i=1}^L\A_{a^{(i)},a^{(i+1)}}^{p_i}=1$ where 
$a^{(i+1)}=a^{(i)}+\sum_{k=0}^{n-1}\La_k(\la^{(i)}_k-\la^{(i)}_{k+1})$ for 
$i=1,\ldots,L$ and $\la^{(i)}=\content(p_i)$.
Then, up to an overall factor,~\eqref{AnRR} is an identity for the generating
function of admissible paths $P\in\P_{\lm}$ starting at $a^{(1)}=(p-n)\La_0$ 
with $\mu$ and $\vL$ related as in~\eqref{Lmu} and 
$\la=(\frac{|\mu|}{n},\cdots,\frac{|\mu|}{n})$. The weights of the paths are
given by $-H(P)$ with $H$ as defined in~\eqref{H}. 

Our initial motivation for studying the A$_{n-1}$ supernomials is
their apparent relevance to a higher-rank generalization of Bailey's 
lemma~\cite{B49}. Indeed, a Bailey-type lemma involving 
the supernomials $S_{\lm}(q)$ such that $\mu^{\star}$ (or any 
permutation thereof) is a partition can be formulated. 
Here we briefly sketch some of our findings. Further details about a 
Bailey lemma and Bailey chain for A$_2$ supernomials are given
in~\cite{ASW98}, whereas we hope the report more on the general A$_{n-1}$ case
in a future publication.

Set
\begin{equation}
f(\vL,\la)=S_{\lm}(q)/(q)_{\vL}
\end{equation} 
for $\vL\in\Int^{n-1}$ and $\la\in\Int^n$ and zero otherwise 
where $(q)_{\vL}=(q)_{L_1}\cdots (q)_{L_{n-1}}$ and
$\mu=((1^{n-1})^{L_{n-1}},\ldots,(1)^{L_1})$. Here 
$(1^i)^{L_i}$ denotes $L_i$ components $(1^i)$.
Note that $\mu^{\star}=(1^{L_1}2^{L_2}\cdots (n-1)^{L_{n-1}})$ is indeed a
partition.

Let $\vL=(L_1,\dots,L_{n-1})$ and $\vk=(k_1,\dots,k_n)$ such that
$k_1+\cdots+k_n=0$ denote arrays of integers and let
$\alpha=\{\alpha_{\vk}\}_{k_1\geq \dots \geq k_n}$,
$\gamma=\{\gamma_{\vk}\}_{k_1\geq \dots \geq k_n}$,
$\beta=\{\beta_{\vL}\}$ and $\delta=\{\delta_{\vL}\}$ be sequences.
Then $(\alpha,\beta)$ and $(\gamma,\delta)$ such that 
\begin{align}
\label{ab}
\beta_{\vL}&=\sum_{\substack{k_1+\dots+k_n=0\\k_1\geq\dots\geq k_n}}
\alpha_{\vk}f(C\vL+\ell\ve_1,L_{n-1}\vrho-\vk+\ell\ve_n),\\
\intertext{and}
\label{cd}
\gamma_{\vk}&=\sum_{\vL\in\Integer^{n-1}} \delta_{\vL}
f(C\vL+\ell\ve_1,L_{n-1}\vrho-\vk+\ell\ve_n)
\end{align}
are called an A$_{n-1}$ Bailey pair relative to $q^{\ell}$ and an A$_{n-1}$
conjugate Bailey pair relative to $q^{\ell}$, respectively.
Here $\ell\in\Int$, $C$ is the Cartan matrix of A$_{n-1}$ and
$\vrho$ is the $n$-dimensional Weyl vector $\vrho=\ve_1+\cdots+\ve_n$.

When $n=2$, $f(L,\la)=1/(q)_{\la_1}(q)_{\la_2}$ for $L=\la_1+\la_2$ and zero
otherwise, and \eqref{ab} and~\eqref{cd} reduce 
(up to factors of $(q)_{\ell}$) to the usual
definition~\cite{B49} of a Bailey pair and conjugate Bailey pair (after
identifying $\vk=(k_1,k_2)=(k,-k)$),
\begin{equation*}
\beta_L=\sum_{k=0}^L\frac{\alpha_k}{(q)_{L-k}(q)_{L+k+\ell}}
\quad\text{and}\quad
\gamma_k=\sum_{L=k}^{\infty}\frac{\delta_L}{(q)_{L-k}(q)_{L+k+\ell}}.
\end{equation*}

Analogous to the A$_1$ case the A$_{n-1}$ Bailey pair and conjugate Bailey 
pair satisfy
\begin{equation}\label{abcd}
\sum_{\substack{k_1+\dots+k_n=0\\k_1\geq\dots\geq k_n}}
\alpha_{\vk}\gamma_{\vk}=
\sum_{\vL\in\Integer^{n-1}} \beta_{\vL}\delta_{\vL}.
\end{equation}
For $n \geq 2$, $N\geq 1$ we now claim the following A$_{n-1}$ 
conjugate Bailey pair relative to $q^{\ell}$.
Choose integers $\la_j^{(a)}\geq 0$ $(a=1,\dots,n-1,~j=1,\dots N-1)$ and 
$\sigma$ such that
\begin{equation}
\ell-\sum_{a=1}^{n-1}\sum_{i=1}^{N-1}ai\la_i^{(a)}+\sigma N\equiv 
0 \pmod{n}.
\end{equation}
Setting $\la=\sum_{a=1}^{n-1}\sum_{i=1}^{N-1}\la_i^{(a)}(\ve_a\otimes\ve_i)$
and $\vk=\vk(\vL)$, such that $k_i(\vL)=L_i-L_{i+1}$ ($L_n=0$, $L_{n+1}=L_1$,
so that $\sum_{i=1}^n k_i=0$) the $(\gamma,\delta)$ pair
\begin{equation}
\label{gd}
\begin{aligned}
\gamma_{\vk(\vL)}&=\frac{q^{\frac{1}{2N}(\vL C\vL+2\ell L_1)}}
{(q)_{\infty}^{n-1}}\sum_{\vn}
\frac{q^{\frac{1}{2}\vn(C\otimes C^{-1})\vn-\vn(I\otimes C^{-1})\lambda}}
{(q)_{\vn}}\\[2mm]
\delta_{\vL}&=q^{\frac{1}{2N}(\vL C\vL+2\ell L_1)}
\sum_{\vn} q^{\frac{1}{2} \vn (C\otimes C^{-1}) \vn-
\vn(I\otimes C^{-1})\lambda}\qbin{\vm+\vn}{\vn}
\end{aligned}
\end{equation}
satisfies~\eqref{cd}. The summations in~\eqref{gd} run over all 
$\vn=\sum_{a=1}^{n-1}\sum_{i=1}^{N-1}n_i^{(a)} (\ve_a\otimes\ve_i)$ such that
\begin{equation}
\frac{L_a+\ell C^{-1}_{a1}-
\sum_{b=1}^{n-1}\sum_{k=1}^{N-1}k C^{-1}_{ab}\lambda_k^{(b)}}{N}-
\sum_{j=1}^{N-1}C^{-1}_{1j} n^{(a)}_j \in \Integer + \frac{a\sigma}{n}
\qquad a=1,\dots,n-1.
\end{equation}
In the expression for $\delta_{\vL}$ the variable $\vm$ is related to $\vn$ by
\begin{equation}
(C\otimes I)\vn+(I\otimes C)\vm=(C\vL+\ell\ve_1)\otimes\ve_{N-1}+\lambda.
\end{equation}

Inserting~\eqref{gd} into~\eqref{abcd} yields a rank $n-1$ and level $N$ 
version of Bailey's lemma. Indeed, when $\la=\ve_a\otimes \ve_i$, 
$\gamma_{\vk}$ is proportional to the level-$N$ A$_{n-1}^{(1)}$ string
function in the representation given by Georgiev~\cite{G95}.
When $n=2$ the pair $(\gamma,\delta)$ of equation~\eqref{gd} reduces to the 
conjugate Bailey pair of refs.~\cite{SW97a,SW96}.
The identities in~\eqref{AnRR} provide A$_{n-1}$ Bailey pairs relative to 1.
We remark that Milne and Lilly~\cite{ML92,ML95} also considered higher-rank 
generalizations of Bailey's lemma. However, their definition of an
A$_{n-1}$ Bailey is different from the one above, and
in particular we note that the function $f$ is not $q$-hypergeometric
for $n>2$.

\section*{Acknowledgements}
We would like to thank M. Okado for generously sharing some of his 
unpublished notes on energy functions of lattice models and
for drawing our attention to ref.~\cite{NY95}. Many thanks to A. Lascoux 
for providing us with copies of his papers and to D. Dei Cont for travelling 
all the way to Italy to get us a copy of ref.~\cite{LS81}.
Furthermore, we would like to acknowledge useful discussions with
P. Bouwknegt, O. Foda, R. Kedem, A. Kuniba, B. McCoy and A. Nakayashiki.
AS has been supported by the ``Stichting Fundamenteel Onderzoek der
Materie'' which is financially supported by the Dutch foundation for
scientific research NWO. SOW has been supported by a fellowship of the
Royal Netherlands Academy of Arts and Sciences.

\section*{Note added}
After submission, several papers~\cite{K98,KS98}~\cite{S98a}-\cite{ShWe98} 
with considerable overlap with this work have appeared. 
The generalized Kostka polynomials studied in this paper were also
introduced in~\cite{ShWe98} as special types of Poincar\'e polynomials
and further studied in \cite{S98a}-\cite{S98c}.
In refs.~\cite{K98,KS98} it was conjectured that the generalized Kostka
polyonomials coincide with special cases of spin generating functions
of ribbon tableaux~\cite{LLT97} and that the fermionic 
representation~\eqref{F} of the generalized Kostka polynomials 
is the generating function of rigged configurations. This last conjecture
has now been established in~\cite{KSS}.

We are indebted to Mark Shimozono for his questions and comments which led 
to several refinements of the paper.
We also thank him for pointing out that the analogues of the recurrences of 
lemma~\ref{lem_rec_S1} for the (cocharge) Kostka polynomials have occurred 
in~\cite{R88}.

\appendix

\section*{Appendix}

\section{Proof of lemma~\ref{lem_chain}}
\label{sec_A}

Obviously, the following lemma implies lemma~\ref{lem_chain}.
\begin{lemma}
\label{lem_h}
Let $P\in\Pb_{\mu}$ be a path over $\{1,2,\ldots,M\}$ where 
$M=|\mu|$ and set $Q=\Cp(P)$. Suppose $M$ is contained in 
step $p_i$ of $P=p_L\otimes\cdots\otimes p_1$.
\begin{enum}
\item If $M$ is contained in the $(i-1)$th step of $\sigma_{i-1}(P)$
then $h_{i-1}(Q)-h_{i-1}(P)=0$.
\item If $M$ is contained in the $i$th step of $\sigma_{i-1}(P)$
then $h_{i-1}(Q)-h_{i-1}(P)=-1$.
\item If $M$ is contained in the $(i+1)$th step of $\sigma_i(P)$
then $h_i(Q)-h_i(P)=0$.
\item If $M$ is contained in the $i$th step of $\sigma_i(P)$
then $h_i(Q)-h_i(P)=1$.
\end{enum}
\end{lemma}

This lemma, in turn, follows from the next lemma. The height of an entry
$M$ in a Young tableau is defined to be $i$ if $M$ is in the $i$th row from 
the bottom.

\begin{lemma}\label{lem_bump}
For $\mu,\mu'\in\R$, let $p\in\B_{\mu}$ and $p'\in\B_{\mu'}$ such that
each entry in $p\cdot p'$ occurs at most once and $p$ contains the largest 
entry $M$. Then
\begin{enum}
\item $\shape(\Cp(p)\cdot p')\neq \shape(p\cdot p')$ if and only if
the height of $M$ in $p\cdot p'$ is bigger than $\height(p)$, and
\item $0\leq h(p\otimes p')-h(\Cp(p)\otimes p')\leq 1$.
\end{enum}
\end{lemma}

Before we prove lemma~\ref{lem_bump} let us first show that it indeed
implies lemma~\ref{lem_h}. 

\begin{proof}[Proof of lemma~\ref{lem_h}]~\newline\indent
$(i)$ Let $\tilde{p}_i\otimes \tilde{p}_{i-1}:=\sigma(p_i\otimes p_{i-1})$.
The steps $p_i$ and $\tilde{p}_{i-1}$ have the same shape and by
assumption they both contain $M$.
Since $p_i\cdot p_{i-1}=\tilde{p}_i\cdot \tilde{p}_{i-1}$ we conclude
that the height of $M$ in $p_i\cdot p_{i-1}$ has to be $\height(p_i)$.
By $(i)$ of lemma~\ref{lem_bump} it follows that $\shape(\Cp(p_i)\cdot
p_{i-1})=\shape(p_i\cdot p_{i-1})$ which proves $h_{i-1}(Q)-h_{i-1}(P)=0$.

$(ii)$ Again we denote $\tilde{p}_i\otimes \tilde{p}_{i-1}:=
\sigma(p_i\otimes p_{i-1})$. By assumption $p_i$ and $\tilde{p}_i$ contain
$M$. Equation $p_i\cdot p_{i-1}=\tilde{p}_i\cdot \tilde{p}_{i-1}$ can only
hold if the box with entry $M$ has been bumped at least once. But this
implies that the height of $M$ in $p_i\cdot p_{i-1}$ is bigger than 
$\height(p_i)$. By $(i)$ of lemma~\ref{lem_bump} this means that
$h_{i-1}(Q)\neq h_{i-1}(P)$, and by $(ii)$ of lemma~\ref{lem_bump} the
difference has to be $-1$.

$(iii)$ This point can be proven analogous to $(i)$.

$(iv)$ Let us show that this case follows from $(ii)$ by considering 
$P'=\Cp^{-1}\circ\Omp(P)$. The path $P'$ satisfies the conditions of 
case~$(ii)$ with $i\to L+1-i$ since $\sigma_i$ commutes with $\Cp^{-1}$ due 
to~\eqref{sC} and~\eqref{CO} and since~\eqref{sO} holds. Hence
$h_{L-i}(\Cp(P'))-h_{L-i}(P')=-1$ which is equivalent to 
$h_{L-i}(\Omp(P))-h_{L-i}(\Omp\circ\Cp(P))=-1$ by inserting the definition 
of $P'$ and using~\eqref{CO} and $\Omp^2=\Id$.
Finally employing~\eqref{hO} proves $(iv)$.
\end{proof}

\begin{proof}[Proof of lemma~\ref{lem_bump}]
Let $p'=[w]$ with $w=w_N\cdots w_1$ in row-representation.
Define $p^{(0)}=p$ and $p^{(i+1)}=p^{(i)}\cdot [w_{N-i}]$ for 
$i=0,1,\ldots,N-1$. Then obviously $p^{(N)}=p\cdot p'$.

We will show inductively that either $M$ got bumped in $p^{(i)}$ (which 
implies that the height of $M$ is bigger than $\height(p)$) or the action 
of $\Cp$ on $p^{(i)}$ is still described by the inverse sliding mechanism 
starting at the largest element $M$ and ending in the bottom left corner.

We prove this claim by induction on $i$.
The initial condition is satisfied since $\Cp$ acts on $p^{(0)}=p$
by the inverse sliding mechanism by definition.
To prove the induction step suppose that $M$ did not yet get bumped in 
$p^{(i)}$ (if it has been bumped in $p^{(i)}$ then this is also true for 
$p^{(k)}$ with $i\leq k\leq N$ and we are finished). By the induction 
hypothesis the action of $\Cp$ on $p^{(i)}$ is still given by the inverse 
sliding mechanism. It is useful to denote the boxes
in $p^{(i)}$ affected by the inverse sliding pictorially by drawing 
$\;\raisebox{-0.3cm}{\epsfxsize=0.6 cm \epsffile{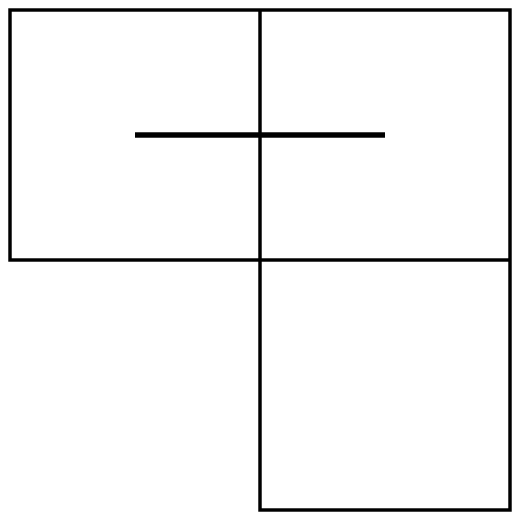}}\;$ 
and $\;\raisebox{-0.3cm}{\epsfxsize=0.6 cm \epsffile{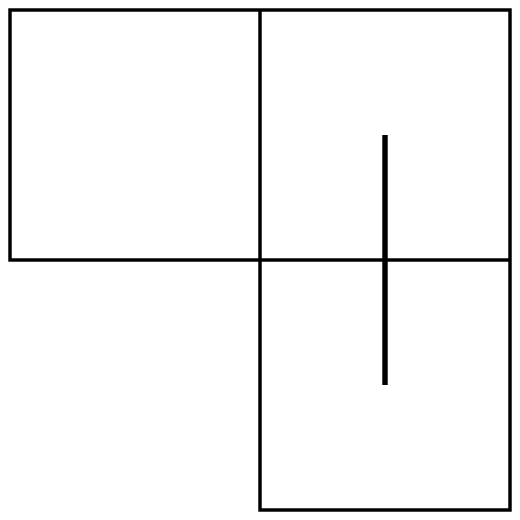}}\;$ 
for $a>b$ and $a<b$, respectively, if the corresponding boxes of $p^{(i)}$ are 
$\;\raisebox{-0.3cm}{\epsfxsize=0.6 cm \epsffile{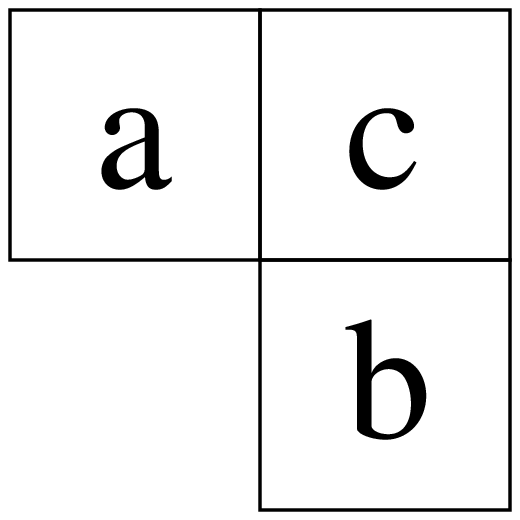}}\;$.
For example
\begin{equation*}
p^{(i)}=\;\raisebox{-0.7cm}{\epsfxsize=3 cm \epsffile{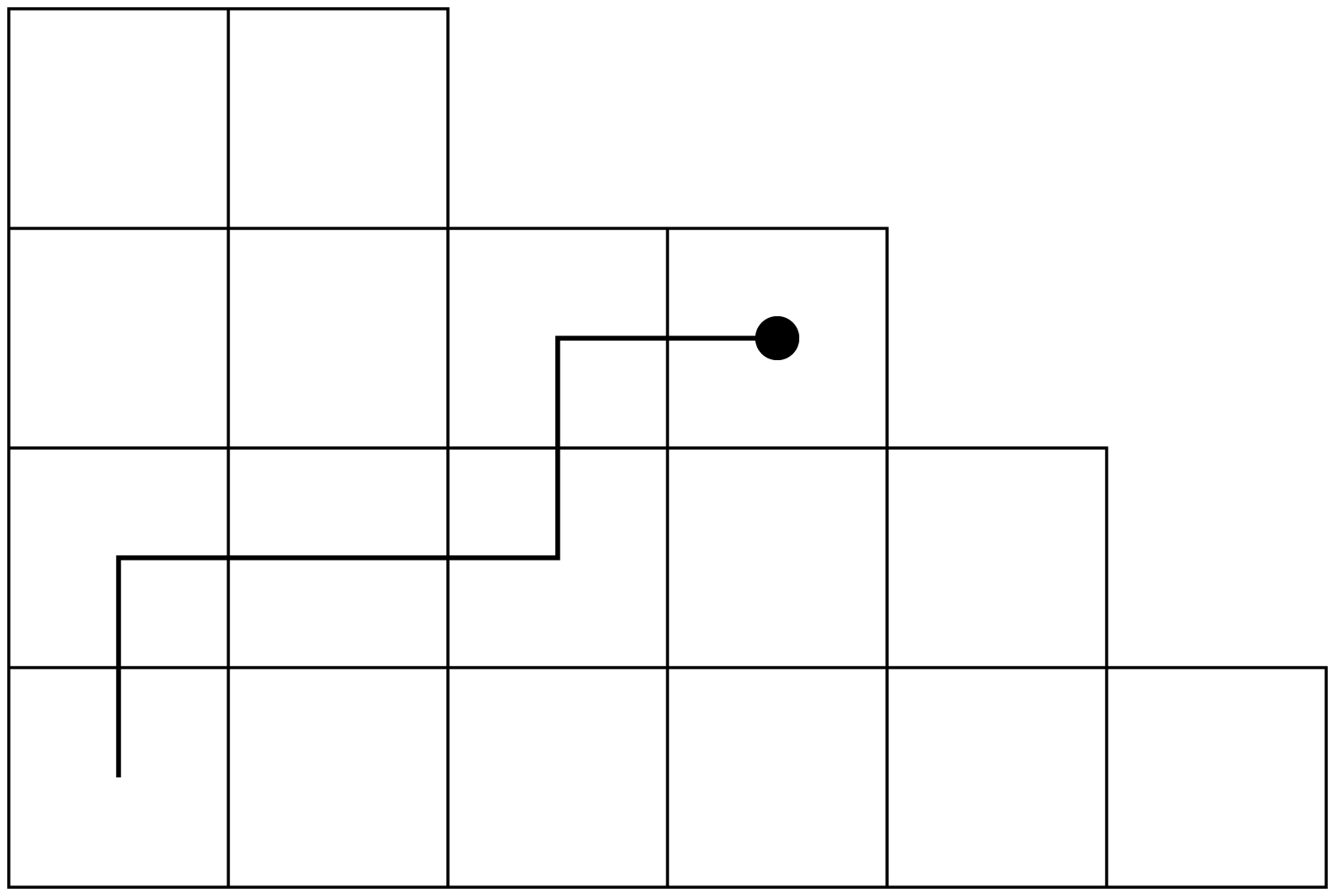}}\;\, .
\end{equation*}
The dot indicates the position of $M$ in $p^{(i)}$.
Comparing with figure~\ref{fig_Islide} we see that the line traces exactly
the movement of an empty box under the inverse sliding mechanism.
We now wish to insert $[w_{N-i}]$ by the Schensted bumping algorithm to obtain 
$p^{(i+1)}$, i.e., $[w_{N-i}]$ gets inserted in the first row of $p^{(i)}$ and 
possibly bumps another box to the second row and so on. Let us label the boxes
of $p^{(i)}$ which get bumped when inserting $[w_{N-i}]$ by a cross. 
Two things may happen:
\begin{enum}
\item[(1)] 
None of the boxes $\; \epsfxsize=0.45 cm \epsffile{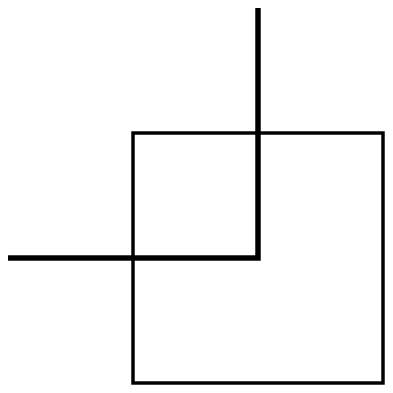}\; $ and
$\; \raisebox{-0.15 cm}{\epsfxsize=0.3 cm \epsffile{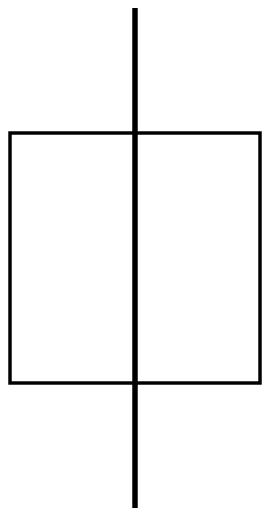}}\; $ contain
a cross. (We include 
$\; \raisebox{-0.15 cm}{\epsfxsize=0.3 cm \epsffile{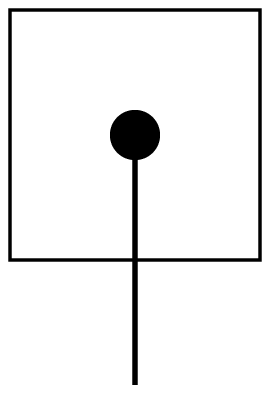}}\; $
in the set of boxes depicted by
$\; \raisebox{-0.15 cm}{\epsfxsize=0.3 cm \epsffile{fig_C6.ps}}\; $).
\item[(2)]
There are boxes $\; \epsfxsize=0.45 cm \epsffile{fig_C5.ps}\; $ or
$\; \raisebox{-0.15 cm}{\epsfxsize=0.3 cm \epsffile{fig_C6.ps}}\; $ which
contain a cross. 
\end{enum}
If (1) occurs there can be at most one box containing both a line and
a cross. One may easily see that the line of $p^{(i)}$ also describes the 
route of the inverse sliding mechanism in $p^{(i+1)}$ and that $M$ does not 
get bumped. Hence we are finished in this case.
If (2) occurs all boxes vertically above
$\; \epsfxsize=0.45 cm \epsffile{fig_C5.ps}\; $ or
$\; \raisebox{-0.15 cm}{\epsfxsize=0.3 cm \epsffile{fig_C6.ps}}\; $
up to and including either ({\em a})
$\; \raisebox{-0.15 cm}{\epsfxsize=0.45 cm \epsffile{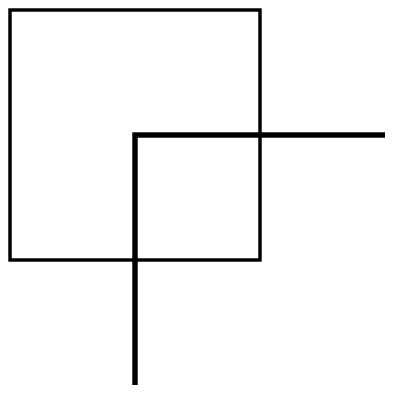}}\; $ or 
({\em b}) $\; \raisebox{-0.15 cm}{\epsfxsize=0.3 cm \epsffile{fig_C8.ps}}\; $
must also contain a cross.
In case ({\em a}) the line indicating the inverse sliding changes from
$p^{(i)}$ to $p^{(i+1)}$ as illustrated in figure~\ref{fig_change}.
\begin{figure}
\centering
\begin{tabular}{c}
\epsfxsize=1.8 cm \epsffile{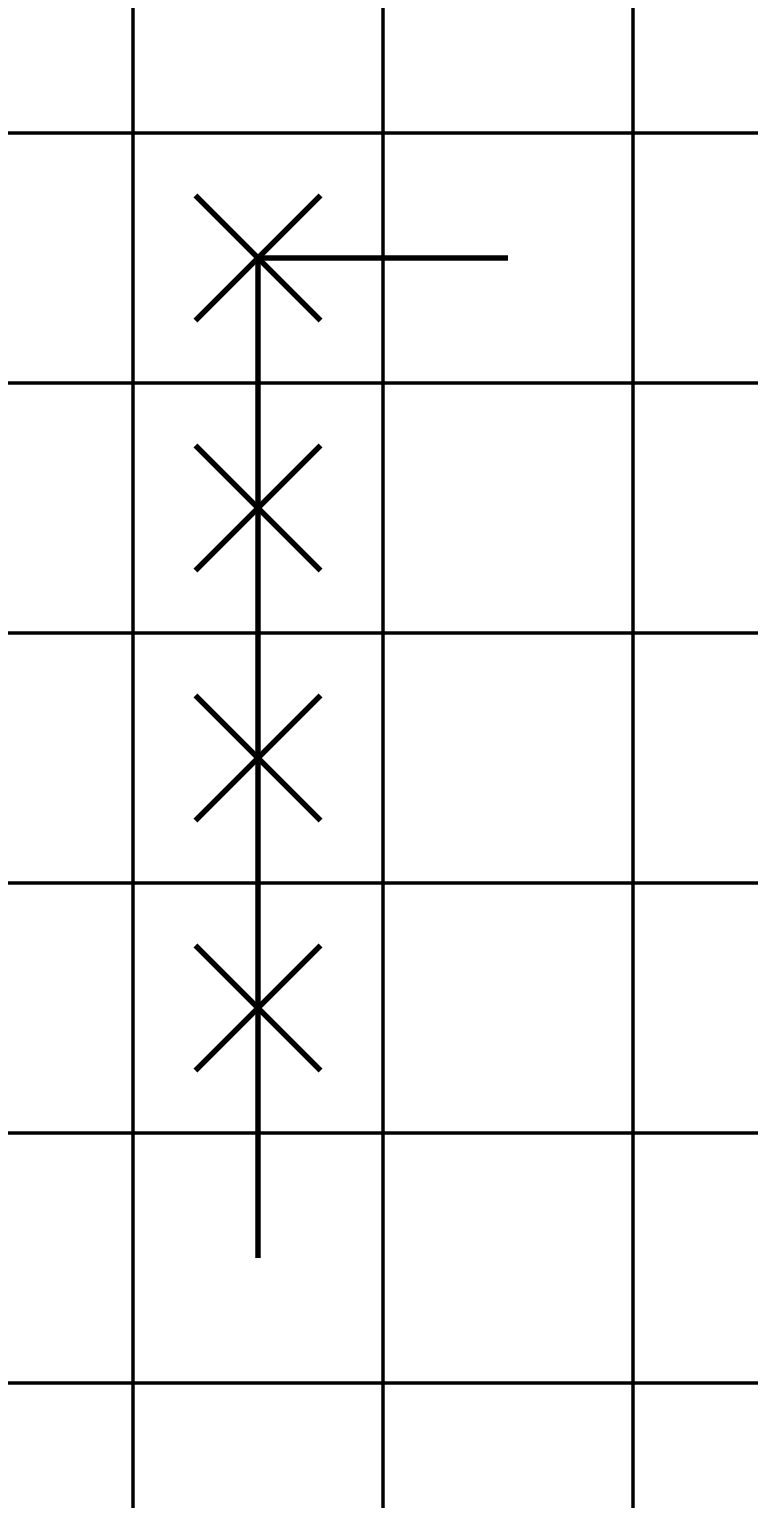} \\[2mm]
in $p^{(i)}$ 
\end{tabular}
\raisebox{0.3cm}{$\;\longrightarrow\;$}
\begin{tabular}{c}
\epsfxsize=1.8 cm \epsffile{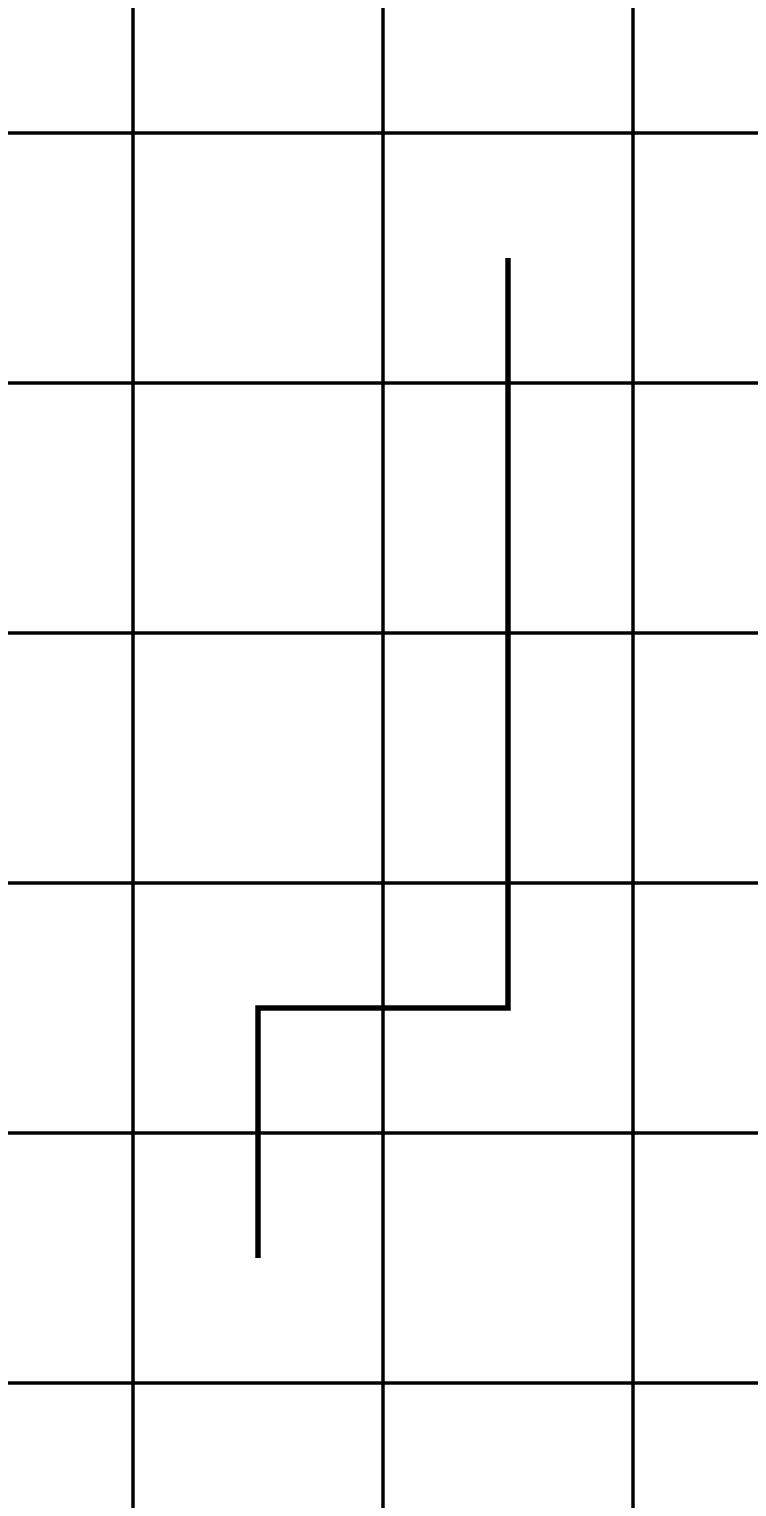}\\[2mm]
in $p^{(i+1)}$ 
\end{tabular}
\hspace{1.5cm}
\begin{tabular}{c}
\epsfxsize=1.8 cm \epsffile{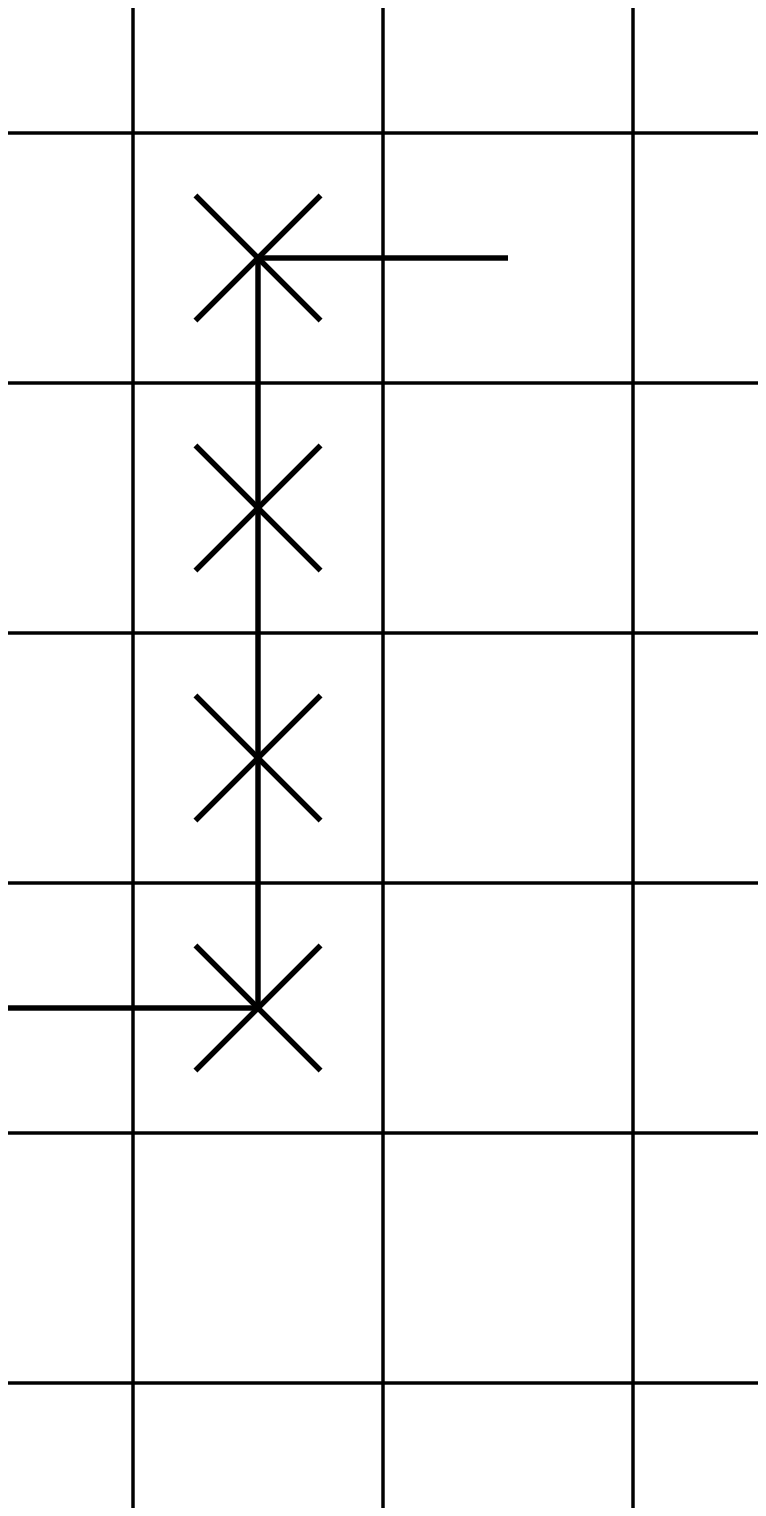} \\[2mm]
in $p^{(i)}$ 
\end{tabular}
\raisebox{0.3cm}{$\;\longrightarrow\;$}
\begin{tabular}{c}
\epsfxsize=1.8 cm \epsffile{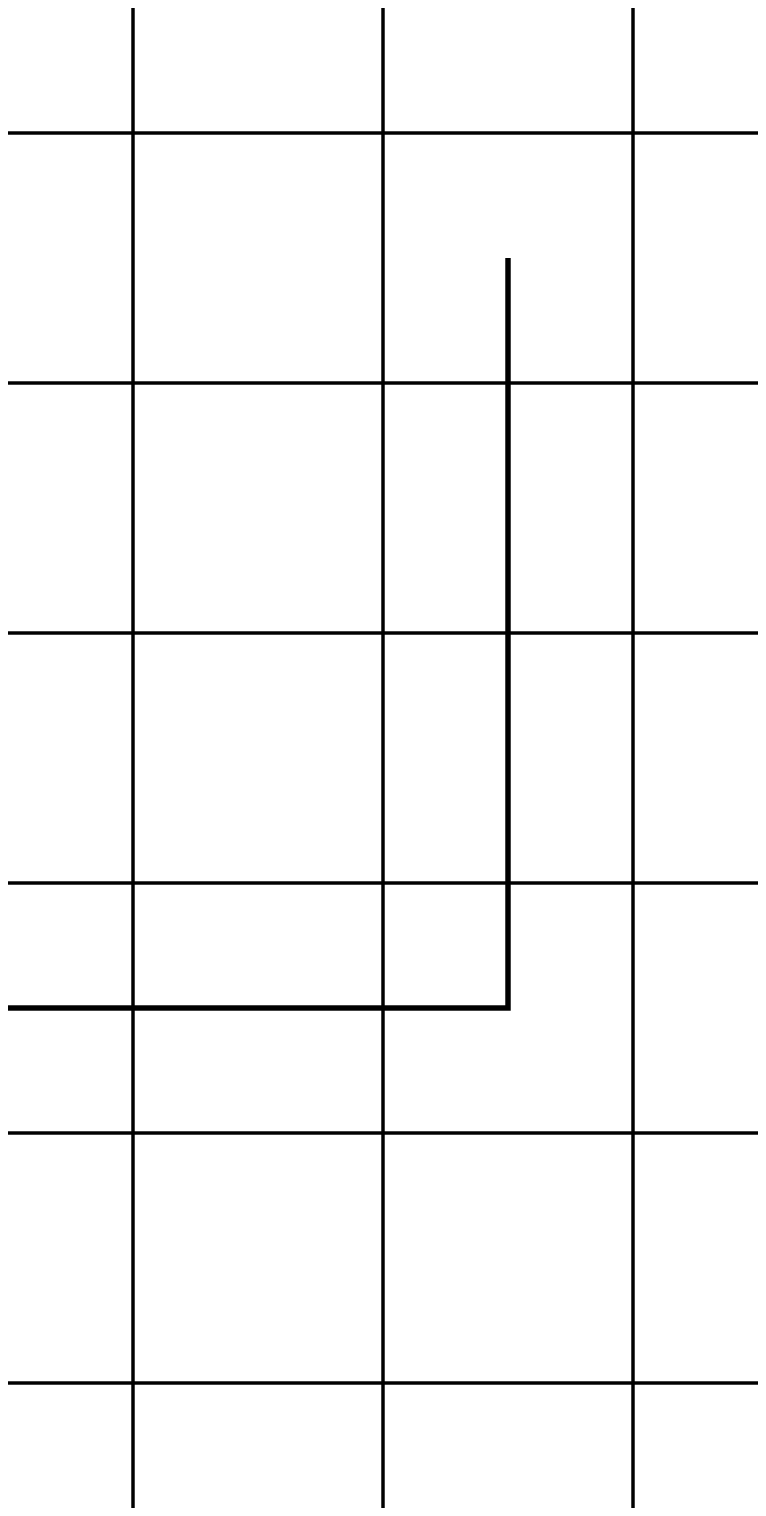}\\[2mm]
in $p^{(i+1)}$ 
\end{tabular}
\caption{Change of inverse sliding route from $p^{(i)}$ to $p^{(i+1)}$
\label{fig_change}}
\end{figure}
In case ({\em b}) 
$\; \raisebox{-0.15 cm}{\epsfxsize=0.3 cm \epsffile{fig_C8.ps}}\; $
contains a cross and hence $M$ got bumped. This concludes the proof of the 
claim.

Observe that, as long as cases (1) or (2{\em a}) occur, $M$ does not get 
bumped and $\shape(p^{(i+1)})=\shape(\Cp(p^{(i+1)})$ since $\Cp$ is still 
described by the inverse sliding mechanism.
If, however, case (2{\em b}) occurs for $p^{(i)}$ which implies that $M$ got 
bumped in $p^{(i+1)}$, then $\shape(p^{(i+1)})$ and $\shape(\Cp(p^{(i+1)}))$
differ. This is so since $p^{(i)}$ must contain
\begin{equation*}
\epsfxsize=2 cm \epsffile{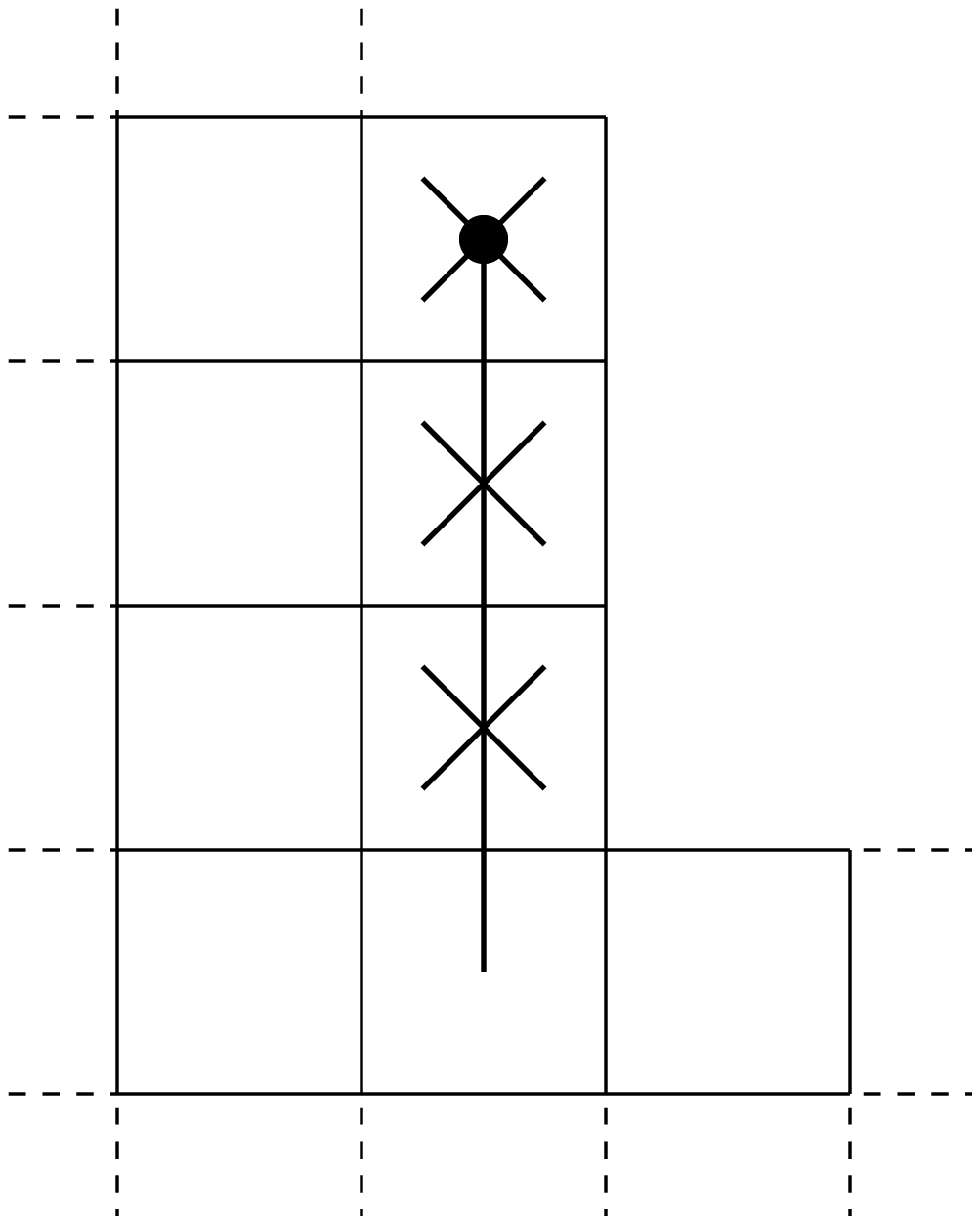}
\quad \raisebox{1cm}{\text{or}}\quad
\epsfxsize=2 cm \epsffile{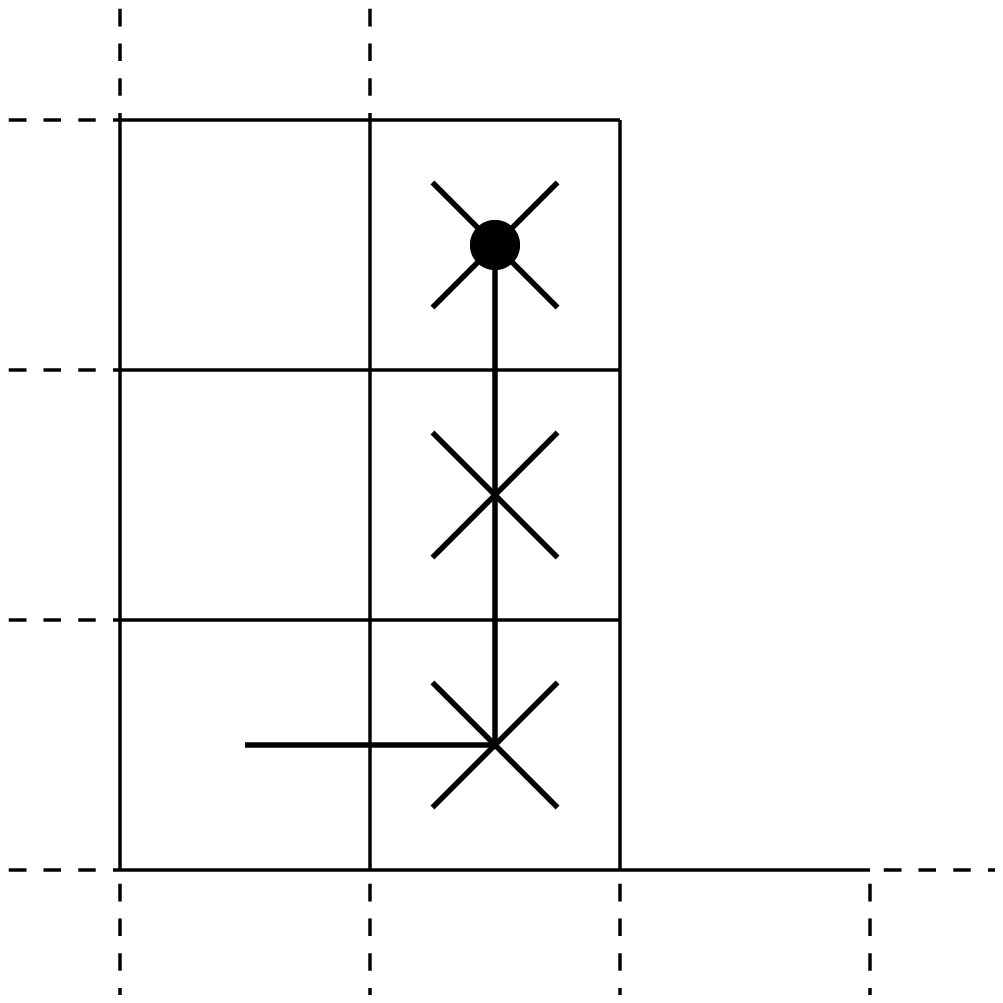}
\end{equation*}
where the dashed lines indicate possible other boxes and the number of 
the vertically aligned crossed boxes may of course vary (but at least 
the box containing the dot must also contain a cross).
Suppose the lowest crossed box in the vertical line below 
$\; \raisebox{-0.15 cm}{\epsfxsize=0.3 cm \epsffile{fig_C8.ps}}\; $ is in 
row $k$. In comparison with $p^{(i)}$, the shape of $p^{(i+1)}$ has one more 
box above the $\height(p)$. In the shape of 
$\Cp(p^{(i+1)})=\Cp(p^{(i)})\cdot [w_{N-i}]$, 
this box has been moved to row $k$. One may also easily see that
$\shape(p^{(j)})$ and $\shape(\Cp(p^{(j)}))$ for $j>i$ differ by moving
exactly one box from above $\height(p)$ to the $k$th row. This proves
lemma~\ref{lem_bump}.
\end{proof}

\newpage
\section{Example of a cyclage-graph}
\label{sec_C}

\begin{figure}[h]
\centering
\epsfysize=16 cm \epsffile{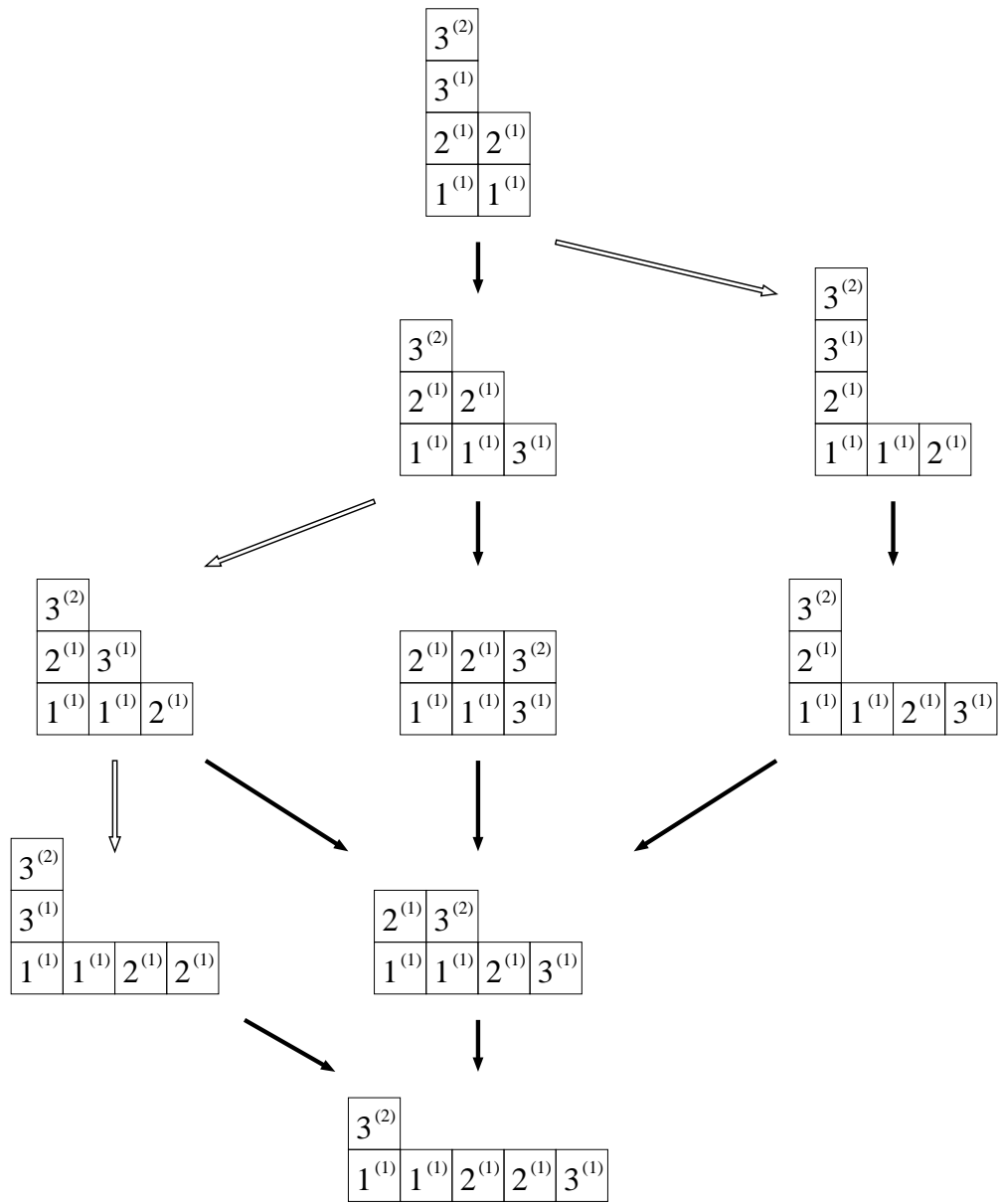}
\caption{The cyclage-graph $\T(\mu)$ for $\mu=((2),(2),(1^2))$
}\label{fig_poset}
\end{figure}

Figure~\ref{fig_poset} shows the poset structure of $\ttab{\cdot}{\mu}$ for 
$\mu=((2),(2),(1^2))$. A black arrow
from LR tableau $T$ to LR tableau $T'$ means $T'=\Cb(T)$. A white arrow indicates
that $T'$ and $T$ are related by a modified $\la$-cyclage (as defined in
section~\ref{sec_poset}) other than the modified initial cyclage, i.e.
$T'=\Zb_{\la}(T)$ for some shape $\la$ but $T'\neq \Cb(T)$.

\newpage


\begin{thebibliography}{99}

\bibitem{A84}
G.~E.~Andrews,
{\em Multiple series Rogers--Ramanujan type identities},
Pacific J. Math. {\bf 114} (1984), 267--283.

\bibitem{A94}
G.~E.~Andrews,
{\em Schur's theorem, Capparelli's conjecture and $q$-trinomial
coefficients},
Contemp. Math. {\bf 166} (1994), 141--154.

\bibitem{AB87}
G.~E.~Andrews and R.~J.~Baxter,
{\em Lattice gas generalization of the hard
hexagon model. III. $q$-trinomial coefficients},
J. Stat. Phys. {\bf 47} (1987), 297--330.

\bibitem{ABF84}
G.~E.~Andrews, R.~J.~Baxter and P.~J.~Forrester,
{\em Eight-vertex SOS model and generalized 
Rogers--Ramanujan--type identities},
J. Stat. Phys. {\bf 35} (1984), 193--266.

\bibitem{ASW98}
G.~E.~Andrews, A.~Schilling and S.~O.~Warnaar, 
{\em An A$_2$ Bailey lemma and Rogers--Ramanujan-type identities},
math.QA/9807125.

\bibitem{B49}
W.~N.~Bailey,
{\em Identities of the Rogers-Ramanujan type},
Proc. London Math. Soc. (2) {\bf 50} (1949), 1--10.

\bibitem{B87}
L.~M.~Butler,
{\em A unimodality result in the enumeration of subgroups of a finite 
abelian group},
Proc. Amer. Math. Soc. {\bf 101} (1987), 771--775.

\bibitem{B91}
L.~M.~Butler,
{\em Generalized flags in finite abelian $p$-groups},
Discrete Appl. Math. {\bf 34} (1991), 67--81.

\bibitem{B94}
L.~M.~Butler,
{\em Subgroup lattices and symmetric functions},
Memoirs of the Amer. Math. Soc., no. 539, vol. {\bf 112} (1994).

\bibitem{D97}
S.~Dasmahapatra,
{\em On the combinatorics of row and corner transfer matrices of the
A$_{n-1}^{(1)}$ restricted face models},
Int. J. Mod. Phys. A {\bf 12} (1997), 3551--3586.

\bibitem{DF96}
S.~Dasmahapatra and O.~Foda,
{\em Strings, paths, and standard tableaux},
Int. J. Mod. Phys. A {\bf 13} (1998), 501--522.

\bibitem{DJKMO89}
E.~Date, M.~Jimbo, A.~Kuniba, T.~Miwa and M.~Okado,
{\em Paths, Maya diagrams and representations of $\widehat{sl}(r,C)$},
Adv. Stud. Pure Math. {\bf 19} (1989), 149--191.

\bibitem{F97}
W.~Fulton,
{\em Young tableaux: with applications to representation theory and geometry},
London Math. Soc. student texts {\bf 35}, Cambridge University Press (1997).

\bibitem{G95}
G.~Georgiev,
{\em Combinatorial constructions of modules for infinite-dimensional Lie 
algebras, II. Parafermionic space},
q-alg/9504024.

\bibitem{HKKOTY98}
G.~Hatayama, A.~N.~Kirillov, A.~Kuniba, M.~Okado, T.~Takagi and Y.~Yamada,
{\em Character formulae of $\widehat{sl}_n$-modules and inhomogeneous paths},
Nucl. Phys. B {\bf 536} [PM] (1998), 575--616.

\bibitem{JMO88}
M.~Jimbo, T.~Miwa and M.~Okado,
{\em Local state probabilities of solvable lattice models:
An A$^{(1)}_{n-1}$ family},
Nucl. Phys. B {\bf 300 [FS22]} (1988), 74--108.

\bibitem{KKMMNN92a}
S.-J. Kang, M. Kashiwara, K.~C. Misra, T. Miwa, T. Nakashima, A. Nakayashiki,
{\em Affine crystals and vertex models},
Int. J. Mod. Phys. A Suppl. {\bf 1A} (1992), 449--484.

\bibitem{KKMMNN92b}
S.-J. Kang, M. Kashiwara, K.~C. Misra, T. Miwa, T. Nakashima, A. Nakayashiki,
{\em Perfect crystals of quantum affine Lie algebras},
Duke Math. J. {\bf 68} (1992), 499--607.

\bibitem{K91}
M.~Kashiwara, 
{\em On crystal bases of the $q$-analogue of universal enveloping algebras},
Duke Math. J. {\bf 63} (1991), 465--516.

\bibitem{K94}
M.~Kashiwara, 
{\em Crystal bases of modified quantized enveloping algebras},
Duke Math. J. {\bf 73} (1994), 383--413.

\bibitem{KN94}
M.~Kashiwara and T.~Nakashima,
{\em Crystal graph for representations of the $q$-analogue of classical
Lie algebras},
J. Alg. {\bf 165} (1994), 295--345.

\bibitem{KKMM93a}
R.~Kedem, T.~R.~Klassen, B.~M.~McCoy and E.~Melzer,
{\em Fermionic quasi-particle representations for characters of
$(G^{(1)})_1\times (G^{(1)})_1/(G^{(1)})_2$},
Phys. Lett. B {\bf 304} (1993), 263--270.

\bibitem{KKMM93b}
R.~Kedem, T.~R.~Klassen, B.~M.~McCoy and E.~Melzer,
{\em Fermionic sum representations for conformal field theory characters},
Phys. Lett. B {\bf 307} (1993), 68--76.

\bibitem{K95}
A.~N.~Kirillov,
{\em Dilogarithm identities},
Prog. Theor. Phys. Suppl. {\bf 118} (1995), 61--142.

\bibitem{K98}
A.~N.~Kirillov,
{\em New combinatorial formula for modified Hall--Littlewood polynomials},
math.QA/9803006.

\bibitem{KKN97}
A.~N.~Kirillov, A.~Kuniba and T.~Nakanishi,
{\em Skew Young diagram method in spectral decomposition of integrable
lattice models II: Higher levels},
q-alg/9711009.

\bibitem{KR88}
A.~N.~Kirillov and N.~Yu.~Reshetikhin,
{\em The Bethe Ansatz and the combinatorics of Young tableaux},
J. Soviet Math. {\bf 41} (1988), 925--955.

\bibitem{KSS}
A.~N.~Kirillov, A.~Schilling and M.~Shimozono,
{\em A bijection between Littlewood--Richardson tableaux and rigged
configurations}, preprint.

\bibitem{KS98}
A.~N.~Kirillov and M.~Shimozono,
{\em A generalization of the Kostka-Foulkes polynomials},
math.QA/9803062.

\bibitem{K70}
D.~E.~Knuth,
{\em Permutations, matrices and generalized Young tableaux},
Pacific J. Math. {\bf 34} (1970), 709--727.

\bibitem{KMOTU96}
A.~Kuniba, K.~C.~Misra, M.~Okado, T.~Takagi and J.~Uchiyama,
{\em Paths, Demazure crystals and symmetric functions},
to appear in Nankai-CRM proceedings 
{\em Extended and quantum algebras and their applications to physics},
q-alg/9612018.

\bibitem{L89}
A.~Lascoux,
{\em Cyclic permutations on words, tableaux and harmonic polynomials},
Proc. Hyderabad Conference on Algebraic Groups 1989, Manoj Prakashan, Madras
(1991), 323--347.

\bibitem{LLT95}
A.~Lascoux, B.~Leclerc and J.-Y.~Thibon,
{\em Crystal graphs and $q$-analogues of weight multiplicities for the
root system $A_n$},
Lett. Math. Phys. {\bf 35} (1995), 359--374.

\bibitem{LLT97}
A.~Lascoux, B.~Leclerc and J.-Y.~Thibon,
{\em Ribbon tableaux, Hall--Littlewood functions, quantum affine
algebras, and unipotent varieties},
J. Math. Phys. {\bf 38} (1997), 1041--1068.

\bibitem{LS78}
A.~Lascoux and M.~P.~Sch\"utzenberger,
{\em Sur une conjecture de H.O. Foulkes},
CR Acad. Sci. Paris {\bf 286A} (1978), 323--324.

\bibitem{LS81}
A.~Lascoux and M.~P.~Sch\"utzenberger,
{\em Le monoid plaxique},
Quaderni della Ricerca scientifica {\bf 109} (1981), 129--156.

\bibitem{M95}
I.~G.~Macdonald,
{\em Symmetric functions and Hall polynomials},
Oxford University Press, second edition (1995).

\bibitem{ML92}
S.~C.~Milne and G.~M.~Lilly,
{\em The $A_{\ell}$ and $C_{\ell}$ Bailey transform and lemma},
Bull. Amer. Math. Soc. (N.S.) {\bf 26} (1992), 258--263.

\bibitem{ML95}
S.~C.~Milne and G.~M.~Lilly,
{\em Consequences of the $A_{\ell}$ and $C_{\ell}$ Bailey transform and
Bailey lemma},
Discrete Math. {\bf 139} (1995), 319--346.

\bibitem{NY95}
A.~Nakayashiki and Y.~Yamada,
{\em Kostka polynomials and energy functions in solvable lattice models},
Selecta Math. (N.S.) 3 (1997), 547--599.

\bibitem{O97}
M.~Okado, private communication.

\bibitem{R88}
F.~Regonati, 
{\em Sui numeri dei sottogruppi di dato ordine dei $p$-gruppi abeliani finiti},
Istit. Lombardo (Rend. Sc.) A {\bf 122} (1988), 369--380.

\bibitem{S61}
C.~Schensted,
{\em Longest increasing and decreasing subsequences},
Canad. J. Math. {\bf 13} (1961), 179--191.

\bibitem{S96}
A.~Schilling,
{\em Multinomials and polynomial bosonic forms for the branching functions
of the $\widehat{su}(2)_M \times \widehat{su}(2)_N / \widehat{su}(2)_{N+M}$
conformal coset models},
Nucl. Phys. B {\bf 467} (1996), 247--271.

\bibitem{SW97a}
A.~Schilling and S.~O.~Warnaar,
{\em A higher-level Bailey lemma},
Int. J. Mod. Phys. B {\bf 11} (1997), 189--195.

\bibitem{SW96}
A.~Schilling and S.~O.~Warnaar,
{\em A higher-level Bailey lemma: proof and application},
The Ramanujan Journal {\bf 2} (1998), 327--349.

\bibitem{SW97}
A.~Schilling and S.~O.~Warnaar,
{\em Supernomial coefficients, polynomial identities and $q$-series},
The Ramanujan Journal {\bf 2} (1998), 459--494.

\bibitem{S63}
M.~P.~Sch{\"u}tzenberger,
{\em Quelques remarques sur une construction de Schensted},
Math. Scand. {\bf 12} (1963), 117--128.

\bibitem{S98a}
M.~Shimozono,
{\em A cyclage poset structure for Littlewood--Richardson tableaux},
math.QA/9804037.

\bibitem{S98b}
M.~Shimozono,
{\em Multi-atoms and monotonicity of generalized Kostka polynomials},
math.QA/9804038.

\bibitem{S98c}
M.~Shimozono,
{\em Affine type A crystal structure on tensor products of rectangles, 
Demazure characters, and nilpotent varieties},
math.QA/9804039.

\bibitem{ShWe98}
M.~Shimozono and J.~Weyman,
{\em Graded characters of modules supported in the closure of a nilpotent 
conjugacy class},
math.QA/9804036.

\bibitem{W97}
S.~O.~Warnaar,
{\em The Andrews--Gordon identities and $q$-multinomial coefficients},
Commun. Math. Phys. {\bf 184} (1997), 203--232.

\end{thebibliography}
\end{document}